\AtBeginDocument{%
  \paperwidth=\dimexpr
    1in + \oddsidemargin
    + \textwidth
    + 1in + \oddsidemargin
  \relax
  \paperheight=\dimexpr
    1in + \topmargin
    + \headheight + \headsep
    + \textheight
    + 1in + \topmargin
  \relax
  \usepackage[pass]{geometry}\relax
}
\RequirePackage{fix-cm}
\PassOptionsToPackage{numbers,sort&compress,square}{natbib}
\documentclass[smallcondensed,natbib]{svjour3}     %
\usepackage[utf8]{inputenc}
\usepackage{bm}
\usepackage{bbm}
\makeatletter
\if@twocolumn
  \renewcommand\normalsize{%
   \@setfontsize\normalsize\@xpt{12.5pt}%
   \abovedisplayskip=3 mm plus6pt minus 4pt
   \belowdisplayskip=3 mm plus6pt minus 4pt
   \abovedisplayshortskip=0.0 mm plus6pt
   \belowdisplayshortskip=2 mm plus4pt minus 4pt
   \let\@listi\@listI}%

  \renewcommand\small{%
   \@setfontsize\small{8.5pt}\@xpt
   \abovedisplayskip 8.5\p@ \@plus3\p@ \@minus4\p@
   \abovedisplayshortskip \z@ \@plus2\p@
   \belowdisplayshortskip 4\p@ \@plus2\p@ \@minus2\p@
   \def\@listi{\leftmargin\leftmargini
               \parsep 0\p@ \@plus1\p@ \@minus\p@
               \topsep 4\p@ \@plus2\p@ \@minus4\p@
               \itemsep0\p@}%
   \belowdisplayskip \abovedisplayskip}

\else
  \if@smallext
   \renewcommand\normalsize{%
   \@setfontsize\normalsize\@xpt\@xiipt
   \abovedisplayskip=3 mm plus6pt minus 4pt
   \belowdisplayskip=3 mm plus6pt minus 4pt
   \abovedisplayshortskip=0.0 mm plus6pt
   \belowdisplayshortskip=2 mm plus4pt minus 4pt
   \let\@listi\@listI}%

  \renewcommand\small{%
   \@setfontsize\small\@viiipt{9.5pt}%
   \abovedisplayskip 8.5\p@ \@plus3\p@ \@minus4\p@
   \abovedisplayshortskip \z@ \@plus2\p@
   \belowdisplayshortskip 4\p@ \@plus2\p@ \@minus2\p@
   \def\@listi{\leftmargin\leftmargini
               \parsep 0\p@ \@plus1\p@ \@minus\p@
               \topsep 4\p@ \@plus2\p@ \@minus4\p@
               \itemsep0\p@}%
   \belowdisplayskip \abovedisplayskip}
 \else
  \renewcommand\normalsize{%
   \@setfontsize\normalsize{9.5pt}{11.5pt}%
   \abovedisplayskip=3 mm plus6pt minus 4pt
   \belowdisplayskip=3 mm plus6pt minus 4pt
   \abovedisplayshortskip=0.0 mm plus6pt
   \belowdisplayshortskip=2 mm plus4pt minus 4pt
   \let\@listi\@listI}%

  \renewcommand\small{%
   \@setfontsize\small\@viiipt{9.25pt}%
   \abovedisplayskip 8.5\p@ \@plus3\p@ \@minus4\p@
   \abovedisplayshortskip \z@ \@plus2\p@
   \belowdisplayshortskip 4\p@ \@plus2\p@ \@minus2\p@
   \def\@listi{\leftmargin\leftmargini
               \parsep 0\p@ \@plus1\p@ \@minus\p@
               \topsep 4\p@ \@plus2\p@ \@minus4\p@
               \itemsep0\p@}%
   \belowdisplayskip \abovedisplayskip}
  \fi
\fi

\makeatother

\usepackage{pgfplots}
\pgfplotsset{compat=1.13}
\usepackage{subcaption}

\definecolor{red}{rgb}{0.7592, 0.3137, 0.3020}
\definecolor{blu}{rgb}{0.3098, 0.5059, 0.7412}
\definecolor{grn}{rgb}{0,0.4,0}
\definecolor{org}{rgb}{1.0,0.5,0.0}
\definecolor{prp}{rgb}{0.4,0.35,0.8}

\usepackage{mathtools}
\DeclarePairedDelimiter{\ceil}{\lceil}{\rceil}
\usepackage{stmaryrd}
\usepackage{amsfonts}
\usepackage[colorlinks=true]{hyperref}
\usepackage{siunitx}
\usepackage{booktabs}
\usepackage{authblk}
\smartqed  %
\DeclareMathOperator{\arcsinh}{arcsinh}
\DeclareMathOperator{\atan2}{atan2}
\DeclareMathOperator{\sgn}{sgn}

\makeatletter
\newcommand{\footnoteref}[1]{%
\ltx@ifpackageloaded{hyperref}{%
  \ifHy@hyperfootnotes%
    \hbox{\hyperref[#1]{%
            \@textsuperscript {\normalfont \ref*{#1}}}}%
  \else%
    \hbox{\@textsuperscript {\normalfont \ref*{#1}}}%
  \fi%
}{%
    \hbox{\@textsuperscript {\normalfont \ref{#1}}}%
 }%
}
\makeatother

\SetMathAlphabet{\mathrm}{normal}{\encodingdefault}{cmss}{\mddefault}{n}
\SetMathAlphabet{\mathrm}{bold}{\encodingdefault}{cmss}{\bfdefault}{n}

\newcommand{\sref}[1]{Section~\ref{#1}}
\newcommand{\aref}[1]{Appendix~\ref{#1}}

\renewcommand{\vec}[1] {\ensuremath{\bm{#1}}}
\newcommand{\dof}[2] {\ensuremath{{\left\{#1\right\}}_{#2}}}
\renewcommand{\Re}{{\mathbb{R}}}
\newcommand{\BB}{{\mathcal{B}(\Omega)}}

\newcommand{\vv}[1]{{\bm{ #1 } }}
\newcommand{\VV}[1]{{\bm{\tilde{ #1}} }}
\newcommand{\mm}[1]{{\bm{{ #1}} }}
\newcommand{\Mm}[1]{{\bm{{\bar{#1}}} }}
\newcommand{\MM}[1]{{\bm{\tilde { #1}} }}

\title{A Non-stiff Summation-By-Parts Finite Difference Method for
the Scalar Wave Equation in Second Order Form: Characteristic Boundary Conditions and
Nonlinear Interfaces
\thanks{%
  B.A.E. was supported by National Science Foundation Awards EAR-1547603 and EAR-1916992\\
  J.E.K. was supported by National Science Foundation Award EAR-1547596\\
  T.H. was supported by National Science Foundation EAR-1916992
  }
}
\titlerunning{Non-stiff SBP in Second-Order Form}

\makeatletter
\let\ORIGINAL@spythm\@spythm{}
\def\@spythm#1#2#3#4[#5]{%
  \NR@gettitle{#5}%
  \ORIGINAL@spythm{#1}{#2}{#3}{#4}[#5]%
}
\makeatother

\begin{document}

\author{%
Brittany A. Erickson
  \and
    Jeremy E. Kozdon
  \and
  Tobias Harvey
}

\institute{%
  B. A. Erickson\at%
  Department of Computer and Information Science \& Department Earth Sciences,\\
  1202 University of Oregon,\\
  1477 E. 13th Ave.\\
  Eugene, OR 97403--1202\\
  \email{bae@uoregon.edu}
  \and
    J. E. Kozdon \at%
  Department of Applied Mathematics,\\
  Naval Postgraduate School,\\
  833 Dyer Road,\\
  Monterey, CA 93943--5216\\
  \email{jekozdon@nps.edu}
  \and
  T. Harvey\at%
  Department of Computer and Information Science\\
  1202 University of Oregon\\
  1477 E. 13th Ave.\\
  Eugene, OR 97403--1202\\
  \email{tharvey2@uoregon.edu}
  \and
  \begin{center}
    The views expressed in this document are those of the authors and do
    not reflect the official policy or position of the Department of
    Defense or the U.S. Government.\\
    Approved for public release; distribution unlimited
  \end{center}
  }

\date{June 28, 2022}

\maketitle

\begin{abstract}
  Curvilinear, multiblock summation-by-parts finite difference operators with
  the simultaneous approximation term method provide a stable and accurate
  framework for solving the wave equation in second order form.
  That said, the standard method can become arbitrarily stiff when
  characteristic boundary conditions and nonlinear interface conditions are
  used.
  Here we propose a new technique that avoids this stiffness by using
  characteristic variables to ``upwind'' the boundary and interface treatment.
  This is done through the introduction of an additional block boundary
  displacement variable.
  Using a unified energy, which expresses both the standard as well as
  characteristic boundary and interface treatment, we show that the
  resulting scheme has semidiscrete energy stability for the scalar anisotropic
  wave equation.
  The theoretical stability results are confirmed with numerical experiments
  that also demonstrate the accuracy and robustness of the proposed scheme.
  The numerical results also show that the characteristic scheme has a time step
  restriction based on standard wave propagation considerations and not the
  boundary closure.
\end{abstract}

\section{Introduction}
Due to their superior dispersion properties, high-order methods are ideally
suited for wave-dominated partial differential equations (PDEs)
\citep{KreissOliger1972}.
That said, unless great care is taken in the treatment of boundary conditions,
interface couplings, and variable coefficients, high-order methods are often
less robust than their low-order counterparts.
An important tool in robust high-order methods is utilization of the
summation-by-parts (SBP) property \citep{KreissScherer1974, KreissScherer1977}.
SBP is the discrete analogue of integration by parts and allows the discrete
stability analysis to mimic the continuous well-posedness analysis
\citep{Nordstrom2017}.

When combined with multiblock domain decompositions and curvilinear coordinates,
SBP finite difference methods can be used to stably and accurately model complex
geometries and variable material parameters.
SBP finite difference methods use standard central difference stencils in the
interior of a domain and transition to one-sided stencils at boundaries and
interfaces in a manner that maintains the SBP property.
An important feature of SBP finite difference methods is the built-in norm
matrix, which is similar to the mass matrix in finite element methods.
A variety of SBP finite difference operators have been developed with the most
relevant to this work being the first and second derivative operators on
unstaggered grids \citep{KreissScherer1974, KreissScherer1977, Strand1994,
MattssonNordstrom2004, Mattsson2012}.
With SBP finite difference methods it is possible to either enforce boundary
conditions strongly \citep{Olsson1995a, Olsson1995b} or weakly
\citep{CarpenterGottliebAbarbanel1994, CarpenterNordstromGottlieb1999}; weak
enforcement of boundary conditions with SBP methods is often called the
simultaneous approximation term (SAT) method and is the approach taken here due
to its ability to handle non-compatible boundary and interface conditions.

We are primarily interested in the wave equation in second-order form, that is,
a displacement formulation of the wave equation as opposed to velocity-stress or
velocity-strain.
Our motivation for this is to address our ultimate goal of advancing simulations
of the earthquake cycle, though our work is applicable to other coupling
situations and whenever characteristic-based boundary conditions are natural.
In earthquake cycle models interseismic loading (decade long tectonic loading)
is coupled to dynamic rupture (earthquake rupture taking place over seconds to
minutes); the importance of this coupling has been recently highlighted in, for
example, \citet{Erickson2020}.
In the interseismic phase, a quasidynamic formulation is often used that
neglects inertial effects, e.g., acceleration, resulting in an elliptic PDE for
the displacement.
During active fault slip, the coseismic rupture phase, the equations have a
hyperbolic character due to the inclusion of inertial effects.

In order to avoid having to transition between displacements and velocity-stress
(or velocity-strain) it is desirable to use a displacement-based formulation for
the coseismic phase.
\citet{VirtaMattsson2014}, building on \citet{MattssonHamIaccarino2008,
MattssonHamIaccarino2009}, developed an SBP-SAT finite difference scheme for the
scalar second-order wave equation with variable coefficients on curved
geometries.
\citet{Duru2019} extended this scheme for use with nonlinear friction laws which
govern the sliding of fault interfaces in earthquake problems; nonlinear
friction laws relate the interface traction to the sliding velocity.
However, as noted in \citet{Duru2019}, the modified scheme that incorporates the
nonlinear friction law results in a numerically stiff system of ordinary
differential equations (ODEs) that prevents the use of, for instance, explicit
Runge-Kutta time stepping methods; in \citet{Duru2019} a custom second-order
accurate time stepping method is used.
Similar numerical stiffness is also seen in the velocity-stress formulation of
the wave equation for earthquake problems, though this can be circumvented by
rewriting the nonlinear friction law in terms of the characteristic variables,
\citep{KozdonDunhamNordstrom2012}.

The main contribution of this work is the use of a characteristic formulation of
boundary and interface conditions within a displacement-based scheme, namely
merging the ideas of \citet{VirtaMattsson2014, Duru2019} with those of
\citet{KozdonDunhamNordstrom2012}.
The key idea of the work is to track the evolution of the boundary and interface
displacements, which allows the use of a non-stiff characteristic formulation.

The paper is structured as follows:
\sref{sec:SBP} introduces our notation for SBP finite difference methods.
Since the main contribution of the work can be seen in a simplified setup,
\sref{sec:1d} presents the continuous and semidiscrete energy analysis for
the one-dimensional scalar wave equation.
Here the characteristic method is formulated for both boundary and
nonlinear interface conditions and compared against the standard approach.
\sref{sec:model:problem} introduces a multi-dimensional anisotropic, scalar wave
equation along with its continuous energy analysis.
The semidiscretization of the multi-dimensional problem is given in
\sref{sec:discretization}, where a unified energy analysis is presented that
incorporates both the proposed characteristic method and the previous
approach.
The paper concludes with multi-dimensional tests problems in
\sref{sec:numerical:experiments} and concluding remarks in
\sref{sec:conclusions}.

The Julia programming language \citep[v1.6.0]{bezanson2017julia} was used for
all simulations with the codes available at
\url{https://github.com/Thrase/sbp_waveprop_characteristic}.

\section{Summation-By-Parts Operators}\label{sec:SBP}
In this work summation-by-parts (SBP) finite difference operators are used to
approximate spatial derivatives.
In order to centralize the explanation of these operators, here we introduce the
one-dimensional operators and then generalize the operators to multiple
dimensions using tensor products.

Throughout we refer to the SBP operators as $2p$ where $p$ is the boundary
accuracy and $2p$ is the interior accuracy.
SBP methods for the wave equation in second order form typically see a global
convergence rate of $\min(2p, p + 2)$, i.e., two orders of accuracy greater than
the boundary accuracy except in the case of $2p = 2$.
For first derivatives we use the operators from \citet{Strand94}\footnote{The
free parameter $x_1=0.70127127127127$ is used for $2p = 6$.} and for second
derivatives the variable coefficient operators from \citet{Mattsson2012}.

\subsection{One Dimensional SBP Operators}
Let the domain $0 \le \xi \le 1$ be discretized with an $N + 1$ equally spaced
grid of points. The grid of points are represented as $\vv{\xi}$ with spacing $h
= 1 / N$ and points located at $\dof{\vv{\xi}}{k} = k h$ for $k = 0, 1, \dots,
N$.  Let $\vv{u}$ be the projection of $u$ onto the computational grid. We
define the operator $\vv{e}_{k}$ to be the grid basis functions, that is the
vector which is $1$ at grid point $k$ and zero at all other grid points.
Importantly $\vv{e}_{k}^{T}$ selects the value of a grid function $\vv{u}$ at
the point $k$, namely $\vv{e}_{k}^{T} \vv{u} = \dof{\vv{u}}{k}$.

Let the first and $C(\xi)$-weighted second derivatives of $u$ be approximated as
\begin{subequations}
  \begin{align}
    {\left.\partial_{1} u\right|}_{\xi_{1} = kh} \approx
    \dof{\mm{D}_{1}\vv{u}}{k},\\
    {\left.\partial_{1} C \partial_{1} u\right|}_{\xi_{1} = kh} \approx
    \dof{\mm{D}_{11}^{(C)}\vv{u}}{k}.
  \end{align}
\end{subequations}
The derivative approximations $\mm{D}_{1}$ and $\mm{D}^{(C)}_{11}$ are called
SBP if they satisfy the following definitions.
\begin{definition}{\bf (SBP First Derivative)}\label{def:sbp1}
  The operator $\mm{D}_{1}$ is called an SBP approximation if it can be
  decomposed as $\mm{H}_{1}\mm{D}_{1} = \mm{Q}_{1}$ with $\mm{H}_{1}$ being a
  symmetric positive definite matrix and
  \begin{equation}
    \vv{u}^{T} \left(\mm{Q}_{1} + \mm{Q}_{1}^{T}\right) \vv{v} =
    \vv{u}^{T} \vv{e}_{N} \vv{e}_{N}^{T} \vv{v} - \vv{u}^{T} \vv{e}_{0}
    \vv{e}_{0}^{T} \vv{v} =
    \dof{\vv{u}}{N} \dof{\vv{v}}{N} - \dof{\vv{u}}{0} \dof{\vv{v}}{0},
  \end{equation} for all vectors $\vv{u}$ and
  $\vv{v}$.
\end{definition}
\begin{definition}{\bf (SBP Second Derivative)}\label{def:sbp2}
  The operator $\mm{D}_{11}^{(C)}$ is called an SBP approximation if it can be
  decomposed as
  \begin{equation}
    \label{eqn:sbp:d2:1d}
    \mm{H}_{1}\mm{D}_{11}^{(C)} = -\mm{A}_{11}^{(C)} +
    \dof{\mm{C}}{N} \vv{e}_{N} \vv{b}_{N}^{T}
    -
    \dof{\mm{C}}{0} \vv{e}_{0} \vv{b}_{0}^{T},
  \end{equation}
  where $\mm{A}_{11}^{(C)}$ is a symmetric positive semidefinite matrix and
  $\vv{b}_{N}^{T} \vv{u}$ and $\vv{b}_{0}^{T} \vv{u}$ are accurate
  approximations of the first derivative of $u$ at the boundary points
  $\dof{\vv{\xi}}{N}$ and $\dof{\vv{\xi}}{0}$, respectively.
\end{definition}

In addition the derivative approximations are assumed to be compatible, namely
that $\mm{H}_{1}$ is the same for both the first and second derivative operators
and the weighting matrix $\mm{H}_{1}$ is diagonal.

\begin{remark}
  It is not assumed that the boundary derivative operators $\vv{b}_{0}^{T}$ and
  $\vv{b}_{N}^{T}$ are the first and last rows of $\mm{D}_{1}$, namely
  $\vv{b}_{0}^{T} \ne \vv{e}_{0} \mm{D}_{1}$ and $\vv{b}_{N}^{T} \ne \vv{e}_{N}
  \mm{D}_{1}$. That is, we do not assume that the operators are fully-compatible
  SBP operators \citep{MattssonParisi2010CICP}.
\end{remark}

The reason operators satisfying Definitions~\ref{def:sbp1} and~\ref{def:sbp2}
are called SBP is that the following identities
\begin{subequations}
  \begin{align}
    \vv{u}^{T}\mm{H}\mm{D}_{1}\vv{v} &=
    \dof{\vv{u}}{N} \dof{\vv{v}}{N} - \dof{\vv{u}}{0} \dof{\vv{v}}{0}
    -
    \vv{u}^{T}\mm{D}_{1}^{T}\mm{H}\vv{v},\\
    \vv{u}^{T}\mm{H}\mm{D}_{11}^{(C)}\vv{v} &=
    \dof{\mm{C}}{N}\dof{\vv{u}}{N} \vv{b}_{N}^{T}\vv{v} -
    \dof{\mm{C}}{0}\dof{\vv{u}}{0} \vv{b}_{0}^{T}\vv{v}
    -
    \vv{u}^{T}\mm{A}_{11}^{(C)}\vv{v},
  \end{align}
\end{subequations}
discretely mimic the continuous integration by parts identities
\begin{subequations}
  \begin{align}
    \int_{0}^{1} u \partial_{1} v &=
    {(uv)|}_{0}^{1} -
    \int_{0}^{1} (\partial_{1} u) v,\\
    \int_{0}^{1} u \partial_{1} C \partial_{1} v &=
    {(Cu\partial_{1}v)|}_{0}^{1} -
    \int_{0}^{1} (\partial_{1} u) C \partial_{1} v.
  \end{align}
\end{subequations}
It is useful to note that $\mm{H}_{1}$ and $\mm{A}^{(C)}_{11}$ lead to
quadrature approximations of the following integrals \citep{HickenZingg2013}:
\begin{subequations}
  \begin{align}
    \int_{0}^{1} u v
    &\approx
    \vv{u}^{T}\mm{H}_{1}\vv{v},\\
    \int_{0}^{1} (\partial_{1} u) C \partial_{1} v
    &\approx
    \vv{u}^{T}\mm{A}_{11}^{(C)}\vv{v}.
  \end{align}
\end{subequations}

\subsection{Multidimensional SBP operators}\label{sec:SBP:MD}
Multiple dimensional SBP operators can be constructed via tensor products. In
particular the one-dimensional operators are applied along the grid lines.
Derivative approximations are taken to be of the form
\begin{equation}
  \partial_{i} C \partial_{j} u \approx \MM{D}_{ij}^{(C)} \VV{u}.
\end{equation}
The variable coefficients $C$ present in the approximation make it cumbersome
to define the form of $\MM{D}_{ij}^{(C)} \VV{u}$, so here we outline some of the
important discrete properties of the operator; \aref{app:2D:sbp} presents the
tensor product construction of the operators in two spatial dimensions from
which the higher dimensional extensions can be generalized.

We define multidimensional SBP operators on the reference domain $\hat{B} = {[0,
1]}^d$ with $d$ being the number of spatial dimensions. A regular Cartesian grid
is used to discretize the reference domain with $N_{i} + 1$ grid points in each
direction and grid spacing $h_{i} = 1 / N_{i}$. A field is represented as a
vector with the leading dimension being the fastest index, i.e., column-major
order. So in two dimensions the grid function of $u(\xi_{1}, \xi_{2})$ is the
vector
\begin{equation}
  \label{eqn:dof:ordering}
  \VV{u} =
  \begin{bmatrix}
    \dof{\VV{u}}{00} & \dof{\VV{u}}{10} & \dots & \dof{\VV{u}}{N_{1}N_{2}}
  \end{bmatrix}^{T},
\end{equation}
where $\dof{\VV{u}}{ij} \approx u(i h_{1}, j h_{2})$.
On the reference domain the faces are numbered so that face $1$ is $\xi_{1} =
0$, face $2$ is $\xi_{1} = 1$, face $3$ is $\xi_{2} = 0$, etc..

Let $\MM{H}$ be the tensor product volume norm matrix,
\begin{equation}
  \MM{H} = \mm{H}_{1} \otimes \cdots \otimes \mm{H}_{d},
\end{equation}
which can be thought of as an approximation of the inner product
\begin{equation}
  \int_{\hat{B}} v u \approx \VV{v} \MM{H} \VV{u}.
\end{equation}
The tensor product derivative operators have the following SBP structure
\begin{equation}\label{eqn:multi:sbp:orig}
  \MM{H} \MM{D}_{ij}^{(C)} = -\MM{A}_{ij}^{(C)}
  + \sum_{f=2i-1}^{2i} \hat{n}^{f}_{i} {\left(\Mm{L}^{f}\right)}^{T}
  \mm{H}^{f}\mm{C}^{f}\Mm{B}^{f}_{j},
\end{equation}
where the two terms in the multidimensional SBP decomposition can be thought of
as approximations of the following volume and surface integrals:
\begin{subequations}
  \begin{align}
    \int_{\hat{B}} (\partial_{i} v) C (\partial_{j} u) &\approx
    \VV{v}^{T}\MM{A}_{ij}^{(C)}\VV{u}^{T},\\
    \int_{\partial \hat{B}^{f}}
    v \hat{n}_{i}^{f} C (\partial_{j} u) &\approx
    \hat{n}_{i} \VV{v} {\left(\Mm{L}^{f}\right)}^{T}
    \mm{H}^{f}\mm{C}^{f}\Mm{B}^{f}_{j} \VV{u}.
  \end{align}
\end{subequations}
If $C_{ij}$ defines a symmetric positive definite, spatially varying coefficient
matrix then the matrix $\MM{A}_{ij}^{(C_{ij})}$ (summation implied over $i$ and
$j$) is symmetric positive semidefinite; see \aref{app:2D:sbp}.

The scalar $\hat{n}_{i}^{f}$ is component $i$ of the outward pointing normal to
the reference domain along face $f$; the Cartesian nature of the reference
domain and the face numbering imply that
\begin{equation}
  \label{eqn:hat:n:kron}
  \hat{n}_{i}^{f} =
  \begin{cases}
    {(-1)}^{f}, & \mbox{ if } i = \ceil*{\frac{f}{2}},\\
    0, & \mbox{ otherwise}.
  \end{cases}
\end{equation}
The matrix $\Mm{L}^{f}$ selects the points from the volume vector along face
$f$.
The matrix $\Mm{B}^{f}_{j}$ computes the derivative approximation in the
direction $\xi_{j}$ and evaluates it along face $f$.
When $i = j$ in~\eqref{eqn:multi:sbp:orig} then $f \in (2j-1, 2j)$ and
$\Mm{B}^{f}_{j}$ is based on the boundary derivatives from the one-dimensional
second derivative SBP operator. When $i \ne j$ in~\eqref{eqn:multi:sbp:orig}
then $f \notin (2j-1, 2j)$ and $\Mm{B}^{f}_{j}$ is based on the first derivative
SBP operator.
The diagonal matrix $\mm{H}^{f}$ is the tensor product surface norm matrix,
which approximates
\begin{equation}
  \int_{\partial \hat{B}^{f}} v u \approx
    \VV{v} {\left(\Mm{L}^{f}\right)}^{T}
    \mm{H}^{f}\Mm{L}^{f} \VV{u},
\end{equation}
and the diagonal matrix $\mm{C}^{f}$ is the variable coefficient evaluated at the
points of face $f$.
Since $\hat{n}_{i}^{f} = 0$ for $f \notin (2i-1, 2i)$ the summation in SBP
decomposition~\eqref{eqn:multi:sbp:orig} can be extended to be a summation over
all faces,
\begin{equation}\label{eqn:multi:sbp}
  \MM{H} \MM{D}_{ij}^{(C)} = -\MM{A}_{ij}^{(C)}
  + \sum_{f=1}^{2d} \hat{n}^{f}_{i} {\left(\Mm{L}^{f}\right)}^{T}
  \mm{H}^{f}\mm{C}^{f}\Mm{B}^{f}_{j};
\end{equation}
this new form will be used to simplify the statement of the discretization.

As noted above, here we have only outlined our basic notation and more details
about the construction of the operators are given in \aref{app:2D:sbp}.

\section{One-Dimensional Example}\label{sec:1d}
To highlight the key contributions of the work we begin with the one dimensional
scalar wave equation with linear boundary and nonlinear interface conditions.
Even in these simplified problem setups the stiffness of the standard
non-characteristic approach can be seen as well as the key ideas and benefits of
the proposed characteristic method.

\subsection{Boundary Treatment}\label{sec:1d:bc}
\subsubsection{Continuous Problem}\label{sec:1d:bc:continuous}
On the domain $\Omega = [0, 1]$ we consider the scalar wave equation
\begin{subequations}\label{eqn:1d:wave}
  \begin{alignat}{3}
    \ddot{u} &= \partial_{1}^2 u, &&~~x \in \Omega ,&&~t \in [0, T],  \label{eqn:1d:wave:pde}
    \end{alignat}
where $\dot{u}$ and $\ddot{u}$ denote the first and second time derivatives of
the displacement $u$. To scalar wave equation~\eqref{eqn:1d:wave:pde} we add the family of boundary conditions
  \begin{alignat}{3}
    \label{eqn:1d:wave:bc}
    \tau &= -\alpha \dot{u}, &&~ ~x \in \{0, 1\},&&~t \in [0, T],
  \end{alignat}
  \end{subequations}
for boundary traction $\tau$ given by
\begin{equation}
  \label{eqn:1d:wave:tau}
  \tau = n \; \partial_{1} u,
\end{equation}
where $n$ is the outward unit normal to the domain, namely $n = -1$ at $x = 0$
and $n = 1$ at $x = 1$.
At each domain edge, the particular type of boundary condition imposed is controlled by the parameter $\alpha \geq 0$. For example $\alpha = 0$ corresponds to a Neumann boundary condition and Dirichlet boundary conditions correspond to $\alpha \rightarrow \infty$, with the latter seen by rewriting~\eqref{eqn:1d:wave:bc} as
\begin{equation}
  \dot{u} = -\frac{\tau}{\alpha}.
\end{equation}
Defining the solution energy,
\begin{equation}
  \label{eqn:1d:wave:energy}
  E = \frac{1}{2} \int_{0}^{1}
  \left(
    \dot{u}^2 + {\left(\partial_{1} u\right)}^2
  \right),
\end{equation}
leads to the following energy estimate and thus well-posedness:
\begin{lemma}\label{lemma:1d:wave:energy:estimate}
  Governing equations~\eqref{eqn:1d:wave} with energy~\eqref{eqn:1d:wave:energy}
  satisfies $\dot{E} \le 0$.
\end{lemma}
\begin{proof}
  Taking the time derivative of the energy~\eqref{eqn:1d:wave:energy}, using the
  scalar wave equation~\eqref{eqn:1d:wave:pde}, and simplifying with
  integration by parts yields
  \begin{equation}
    \label{eqn:1d:energy:pre:estimate}
    \dot{E}
    =
    \int_{0}^{1}
    \left(
      \dot{u}\left(\partial_{1}^2 u\right) + \left(\partial_{1} u\right) \left(\partial_{1} \dot{u}\right)
    \right)
    =
    {\left.\dot{u}\;\partial_{1} u\right|}_{0}^{1}
    =
    {\left.\dot{u}\;\tau \right|}_{1}
      +
    {\left.\dot{u}\;\tau \right|}_{0}.
  \end{equation}
  Applying boundary conditions~\eqref{eqn:1d:wave:bc} leads to
  \begin{equation}
    \label{eqn:1d:energy:estimate}
    \dot{E}
    =
    -{\left.\alpha\dot{u}^{2} \right|}_{1}
    -{\left.\alpha\dot{u}^{2} \right|}_{0}
    \le 0,
  \end{equation}
  since $\alpha \ge 0$.
  \qed
\end{proof}

For the numerical scheme that follows, it is useful to express the
boundary conditions~\eqref{eqn:1d:wave:bc} in characteristic form:
\begin{equation}
  \label{eqn:1d:wave:char:bc}
  \dot{u} + \tau
  =
  R \left( \dot{u} - \tau \right), ~ ~x \in \{0, 1\},~t \in [0, T].
\end{equation}
Here $R$, known as the reflection coefficient, is related to $\alpha$ by the
transformations
\begin{equation}
  \label{eqn:1d:wave:alpha:R}
  R = \frac{1 - \alpha}{1 + \alpha}
  \quad \mbox{and} \quad
  \alpha =\frac{1-R}{1+R},
\end{equation}
where if $\alpha \ge 0$ then $-1 \le R \le 1$ and vice versa.  For example, the Neumann condition ($\alpha = 0$) corresponds to $R = 1$, and the Dirichlet condition ($\alpha \to \infty$) corresponds to $R = -1$.
The terms $\dot{u} \pm \tau$ in boundary condition~\eqref{eqn:1d:wave:char:bc}
correspond to characteristic variables propagating in the $\mp n$ directions,
respectively; recall that the unit normal $n$ is included in traction
definition~\eqref{eqn:1d:wave:tau}.
This implies that the boundary condition takes the form of expressing the
incoming characteristic as a reflection of the outgoing characteristic.
Since~\eqref{eqn:1d:wave:char:bc} is
completely equivalent to~\eqref{eqn:1d:wave:bc}
the continuous energy estimate of Lemma~\ref{lemma:1d:wave:energy:estimate}
still holds.

In the semidiscrete scheme that follows it will be useful to define
\begin{equation}
  \label{eqn:1d:wave:char}
  w = \dot{u} - \tau
  \quad\mbox{and}\quad
  q = \dot{u} + \tau,
\end{equation}
with $w$ being the outward and $q$ the inward propagating characteristic
variables.
We reiterate that for the continuous problem, the boundary conditions specified
in terms of $\alpha$ and $R$ are completely equivalent.
That said, the two different perspectives on the boundary conditions will lead
to discretizations with differing numerical properties.

\subsubsection{Semidiscrete Problem}\label{sec: sem_bd}

Using the SBP finite difference operators defined in \sref{sec:SBP}, a family of
SBP-SAT semidiscretizations for governing equations~\eqref{eqn:1d:wave} is
\begin{equation}
  \label{eqn:1d:wave:disc}
  \ddot{\vv{u}} = \mm{D}_{11} \vv{u}
  +
  \sum_{\mathclap{k \in \{0, N\}}}
  \left(\mm{H}^{-1}\vv{e}_{k}\left(\tau^{*}_{k} - n_{k} \vv{b}_{k}^{T} \vv{u}\right)
  - n_{k}\mm{H}^{-1}\vv{b}_{k}\left(u^{*}_{k} - u_{k}\right)\right),
\end{equation}
where $n_{0} = -1$ and $n_{N} = 1$ are the outward pointing normals at the
boundary points.
Additionally, $\tau^{*}_{k}$ and $u^{*}_{k}$ with $k \in \{0, N\}$ are yet-to-be-defined
numerical fluxes that penalize the boundary derivative and displacement toward values that satisfy a range of desired boundary conditions.  For example a Neumann condition $-\partial_1 u = 0$ at $x = 0$ can be imposed by making the choices
\begin{subequations}
  \begin{alignat}{3}
\tau^*_0 &= 0, \\
u^*_0 & = u_0.
    \end{alignat}
\end{subequations}
We will illustrate however, that the numerical fluxes can be chosen to impose the boundary conditions using either the non-characteristic~\eqref{eqn:1d:wave:bc} or characteristic~\eqref{eqn:1d:wave:char:bc} formulations, which though equivalent in the continuous setting, do not lead to the same numerical scheme.

The semidiscrete version of solution energy~\eqref{eqn:1d:wave:energy} is
\begin{equation}
  \label{eqn:1d:wave:disc:energy}
  \begin{split}
    \mathcal{E} =&\; \frac{1}{2}
    \left(
      \dot{\vv{u}}^{T} \mm{H} \dot{\vv{u}}
      +
      \vv{u}^{T} \mm{A}_{11} \vv{u}
    \right)
       +
       \frac{1}{2\gamma}\sum_{{k \in \{0, N\}}}\left(
      \tau_{k}^{2}
      -
      {\left(\vv{b}_{k}^{T}\vv{u}\right)}^{2}
    \right),
  \end{split}
\end{equation}
with $\tau_{k}$ being the following approximation of traction:
\begin{equation}
  \label{eqn:1d:wave:disc:tau}
  \tau_{k} = n_{k} \vv{b}_{k}^{T}\vv{u} + \gamma\left(u^{*}_{k} - u_{k}\right),
\end{equation}
where $\gamma \ge 0$ is a penalty parameter.
The following theorem characterizes energy~\eqref{eqn:1d:wave:disc:energy}:
\begin{theorem}\label{thm:disc:1d:energy:seminorm}
  Energy~\eqref{eqn:1d:wave:disc:energy} is a seminorm for all $\vv{u}$ and
  $u^{*}_{k}$ if $\gamma$ is positive and sufficiently large.
\end{theorem}
\begin{proof}
  The borrowing lemma from \citet[Lemma 1]{AlmquistDunham2020} states that
  \begin{equation}
    \vv{u}^{T} \mm{A}_{11} \vv{u} \ge
    \sum_{\mathclap{k \in \{0, N\}}}\;
    \left(
      \theta {\left(\vv{d}_{k}^{T}\vv{u}\right)}^{2}
      +
      \zeta {\left(\vv{\Delta}_{k}^{T}\vv{u}\right)}^{2}
    \right).
  \end{equation}
  Here, parameters $\theta > 0$ and $\zeta > 0$ depend on the specific SBP
  operator but are independent of the grid spacing; see
  Table~\ref{tab:borrowing}.
  The operator $\vv{d}_{k}^{T} \vv{u}$ for $k \in \{0, N\}$ is an approximation
  of the boundary derivative and $\vv{\Delta}_{k}^{T} = \vv{b}_{k}^{T} -
  \vv{d}_{k}^{T}$; the specific form of $\vv{d}_{k}^{T}$ does not matter for the
  one-dimensional analysis but is critical in the multi-dimensional case%
  \footnote{
    An alternative borrowing lemma is given in \citet{VirtaMattsson2014}, and
    though this lemma yields a slightly better bound on the penalty parameter
    $\gamma$ in one-dimensional, \citet{AlmquistDunham2020} have shown that the
    bound is significantly worse in multiple dimensions.
  }.
  Expanding the difference term in energy~\eqref{eqn:1d:wave:disc:energy} and
  rewriting in terms of $\vv{\Delta}_{k}^{T}$ and $\vv{d}_{k}^{T}$ yields
  \begin{equation}
    \tau_{k}^{2}
    -
    {\left(\vv{b}_{k}^{T}\vv{u}\right)}^{2}
    =
    2\gamma n_{k} \delta_{k} \vv{d}_{k}^{T}\vv{u}
    +
    2\gamma n_{k} \delta_{k} \vv{\Delta}_{k}^{T}\vv{u}
    +
    \gamma^{2} \delta_{k}^{2},
  \end{equation}
  where we have defined $\delta_{k} = u_{k}^{*} - u_{k}$ for $k \in \{0, N\}$.
  Putting these results together and completing the square gives
  \begin{equation}
    \begin{split}
      &\vv{u}^{T} \mm{A}_{11} \vv{u}
      +
      \frac{1}{\gamma}
      \sum_{{k \in \{0, N\}}}
      \left(
        \tau_{k}^{2}
        -
        {\left(\vv{b}_{k}^{T}\vv{u}\right)}^{2}
      \right)
      \ge\\
      &~
      \sum_{\mathclap{k \in \{0, N\}}}\;\;
      \left(
        \theta{\left(n_{k} \vv{d}_{k}^{T}\vv{u} + \frac{1}{\theta}\delta_{k}\right)}^{2}
        +
        \zeta{\left(n_{k} \vv{\Delta}_{k}^{T}\vv{u} + \frac{1}{\zeta}\delta_{k}\right)}^{2}
        +
        \left(\gamma - \frac{1}{\theta} - \frac{1}{\zeta}\right)\delta_{k}^2
      \right).
    \end{split}
  \end{equation}
  Choosing $\gamma$ so that
  \begin{equation}
    \label{eqn:1d:wave:disc:gamma}
    \gamma \ge \frac{1}{\theta} + \frac{1}{\zeta},
  \end{equation}
  leads to the estimate
  \begin{equation}
    \vv{u}^{T} \mm{A}_{11} \vv{u}
    +
    \frac{1}{\gamma}
    \sum_{{k \in \{0, N\}}}
      \left(
        \tau_{k}^{2}
        -
        {\left(\vv{b}_{k}^{T}\vv{u}\right)}^{2}
      \right)
    \ge 0,
  \end{equation}
  for all $\vv{u}$ and $u^{*}$ since $\theta > 0$ and $\zeta > 0$.
  The theorem then follows because $\mm{H}$ is positive definite.
  \qed
\end{proof}
\begin{corollary}
  Semidiscretization~\eqref{eqn:1d:wave:disc} satisfies
  \begin{equation}
    \label{eqn:1d:wave:disc:energy:rate}
    \dot{\mathcal{E}} = \sum_{\mathclap{k \in \{0, N\}}}
    \left(
      \dot{u}_{k}^{*}\tau_{k}^{*}
      -
      \left(\dot{u}_{k} - \dot{u}^{*}_{k}\right)
      \left(\tau_{k} - \tau^{*}_{k}\right)
    \right),
  \end{equation}
  where $\mathcal{E}$ is defined by~\eqref{eqn:1d:wave:disc:energy}.
\end{corollary}
\begin{proof}
  Follows directly by taking the time derivative of
  energy~\eqref{eqn:1d:wave:disc:energy}, substituting
  semidiscretization~\eqref{eqn:1d:wave:disc}, and applying SBP
  property~\eqref{eqn:sbp:d2:1d}.
  \qed
\end{proof}

\subsubsection{Semidiscrete Problem: Non-Characteristic Boundary Treatment}

If the numerical fluxes in semidiscretization~\eqref{eqn:1d:wave:disc} are
defined to impose boundary conditions in the non-characteristic form~\eqref{eqn:1d:wave:bc}, namely,
\begin{subequations}
  \label{eqn:1d:wave:disc:nonchar:bc:flux}
  \begin{align}
    \tau^{*}_{k} &= -\alpha \dot{u}_{k},\\
    u^{*}_{k} &= u_{k},
  \end{align}
\end{subequations}
then the numerical scheme is equivalent to that of \citet{VirtaMattsson2014}.
In this case, the definition of $u^{*}$ implies that the terms with $\gamma$
drop out from energy~\eqref{eqn:1d:wave:disc:energy} and the penalty parameter
has no impact on the scheme.  Neumann boundary conditions can be imposed by setting $\alpha = 0$ in~\eqref{eqn:1d:wave:disc:nonchar:bc:flux}, however it becomes clear that the imposition of Dirichlet boundary conditions ($\alpha \to \infty$) is problematic with this formulation.

\begin{theorem}
  Semidiscretization~\eqref{eqn:1d:wave:disc} with numerical
  fluxes~\eqref{eqn:1d:wave:disc:nonchar:bc:flux} satisfies $\dot{\mathcal{E}}
  \le 0$ where $\mathcal{E}$ is defined by~\eqref{eqn:1d:wave:disc:energy}.
\end{theorem}
\begin{proof}
  Substituting numerical fluxes~\eqref{eqn:1d:wave:disc:nonchar:bc:flux} into
  energy rate~\eqref{eqn:1d:wave:disc:energy:rate} leads to
  \begin{equation}
    \dot{\mathcal{E}}
    =
    -\sum_{\mathclap{k \in \{0, N\}}} \alpha {\left(\dot{u}_{k}^{*}\right)}^{2} \le 0
  \end{equation}
  since $\alpha \ge 0$.
  \qed
\end{proof}

\subsubsection{Semidiscrete Problem: Characteristic Boundary Treatment}

The numerical fluxes in semidiscretization~\eqref{eqn:1d:wave:disc} can instead be chosen to impose the characteristic boundary condition~\eqref{eqn:1d:wave:char:bc}, giving rise to a different numerical discretization that defined by non-characteristic numerical flux~\eqref{eqn:1d:wave:disc:nonchar:bc:flux}. The basic idea of the characteristic method is to mimic the upwinding procedure
used in methods for first order hyperbolic equations.
Namely, we seek to modify the incoming characteristic variable while preserving
to outgoing characteristic variable.
To do this, we introduce an equation for $\dot{u}_{k}^{*}$ which tracks the time
evolution of the numerical flux.

The characteristic-based numerical fluxes are defined as
\begin{subequations}
  \label{eqn:1d:wave:disc:char:bc:flux}
  \begin{align}
    \dot{u}^{*}_{k} &= \frac{q^{*}_{k} + w^{*}_{k}}{2},\\
    \tau^{*}_{k} &= \frac{q^{*}_{k} - w^{*}_{k}}{2},
  \end{align}
\end{subequations}
where the characteristic variables are defined so that the outgoing
characteristic is preserved and the incoming characteristic satisfies the
boundary condition:
\begin{subequations}
  \label{eqn:1d:wave:disc:char:bc}
  \begin{align}
    w_{k}^{*} &= \dot{u}_{k} - \tau_{k},\\
    q_{k}^{*} &= R\;w_{k}^{*}
  \end{align}
\end{subequations}
for any $-1 \leq R \leq 1$.  Importantly the characteristic numerical flux $u^{*}_{k}$ must be tracked as an
independent variable in the solution process. To illustrate the characteristic approach, both Neumann and Dirichlet conditions can be enforced in a straightforward manner by choosing either $R = 1$ or $R = -1$, respectively, in~\eqref{eqn:1d:wave:disc:char:bc}. For $R = -1$, for example, solving~\eqref{eqn:1d:wave:disc:char:bc:flux} and~\eqref{eqn:1d:wave:disc:char:bc} for the numerical fluxes yields the specific choices
\begin{subequations}
  \label{eqn:1d:wave:disc:char:dirichlet}
  \begin{align}
    u_{k}^{*} &= 0,\\
    \tau_{k}^{*} &=  \tau_{k} - \dot{u}_{k} ,
  \end{align}
\end{subequations}
corresponding to Dirichlet boundary conditions.
\begin{theorem}
  Semidiscretization~\eqref{eqn:1d:wave:disc} with numerical
  fluxes~\eqref{eqn:1d:wave:disc:char:bc:flux} satisfies $\dot{\mathcal{E}}
  \le 0$ where $\mathcal{E}$ is defined by~\eqref{eqn:1d:wave:disc:energy}.
\end{theorem}
\begin{proof}
  By definition, $\tau^{*}_{k}$ and $\dot{u}_{k}^{*}$ satisfy the boundary
  condition:
  \begin{equation}
    \label{eqn:1d:wave:disc:char:bc:alpha}
    \tau^{*}_{k} =
    -\frac{1-R}{2}w_{k}^{*} =
    -\left(\frac{1-R}{1 + R}\right)
    \left(\frac{1+R}{2}\right) w_{k}^{*} =
    -\alpha
    \left(\frac{w_{k}^{*}+q_{k}^{*}}{2}\right) =
    -\alpha \dot{u}_{k}^{*},
  \end{equation}
  where we have used numerical flux
  definitions~\eqref{eqn:1d:wave:disc:char:bc:flux}
  and~\eqref{eqn:1d:wave:disc:char:bc} along with the
  relationship~\eqref{eqn:1d:wave:alpha:R} between $\alpha$ and $R$.
  If the grid based incoming characteristic is defined as
  \begin{equation}
    q_{k} = \dot{u}_{k} + \tau_{k},
  \end{equation}
  then it follows that
  \begin{equation}
    \label{eqn:1d:wave:disc:char:diff}
    \dot{u}_{k} - \dot{u}_{k}^{*}
    =
    \tau_{k} - \tau_{k}^{*}
    =
    \frac{q_{k} - q_{k}^{*}}{2}.
  \end{equation}
  Using~\eqref{eqn:1d:wave:disc:char:bc:alpha}
  and~\eqref{eqn:1d:wave:disc:char:diff} in energy
  rate~\eqref{eqn:1d:wave:disc:energy:rate} gives
  \begin{equation}
    \dot{\mathcal{E}} = -\sum_{\mathclap{k \in \{0, N\}}}\;
    \left(
      \alpha {\left(\dot{u}_{k}^{*}\right)}^{2}
      +
      \frac{{\left(q_{k} - q^{*}_{k}\right)}^{2}}{4}
    \right) \le 0
  \end{equation}
  since $\alpha \ge 0$.
  \qed
\end{proof}

\subsubsection{Numerical Results}

In order to integrate semidiscretization~\eqref{eqn:1d:wave:disc} in time
one can either use methods designed for second order ODEs, such as those
proposed in \citet{Duru2019}, or transform the equations into a
first order system of ODEs and use first order time integration technology such
as Runge-Kutta methods.
In this latter approach, the auxiliary variable $\vv{v} = \dot{\vv{u}}$ is
introduced and semidiscretization~\eqref{eqn:1d:wave:disc} becomes
\begin{subequations}\label{eqn:1d:wave:disc:ODE}
  \begin{alignat}{2}
    \dot{\vv{v}} &= \mm{D}_{11} \vv{u}
    +
    \sum_{\mathclap{k\in\{0,N\}}}\;
    \left(
      \mm{H}^{-1}\vv{e}_{k}\left(\tau^{*}_{k} - n_{k} \vv{b}_{k}^{T} \vv{u}\right)
      + \mm{H}^{-1}\vv{b}_{k}\left(u^{*}_{k} - u_{k}\right)
    \right),
    \\
    \dot{\vv{u}} &= \vv{v};
  \end{alignat}
\end{subequations}
in the case of the characteristic numerical
flux~\eqref{eqn:1d:wave:disc:char:bc:flux} additional ODEs are required to
track $u^{*}_{k}$ for $k \in \{0, N\}$.
System of ODEs~\eqref{eqn:1d:wave:disc:ODE} can be written more compactly as
\begin{subequations}
  \label{eqn:ODE:system}
  \begin{align}
,
  \end{align}
\end{subequations}
where linear operator $\mm{G}$ corresponds to the characteristic numerical
flux~\eqref{eqn:1d:wave:disc:char:bc:flux} and $\bar{\mm{G}}$ non-characteristic
numerical flux~\eqref{eqn:1d:wave:disc:nonchar:bc:flux}.
For each method, the time step size is controlled by the eigenvalue spectrum of
the corresponding linear operator.

\begin{figure}
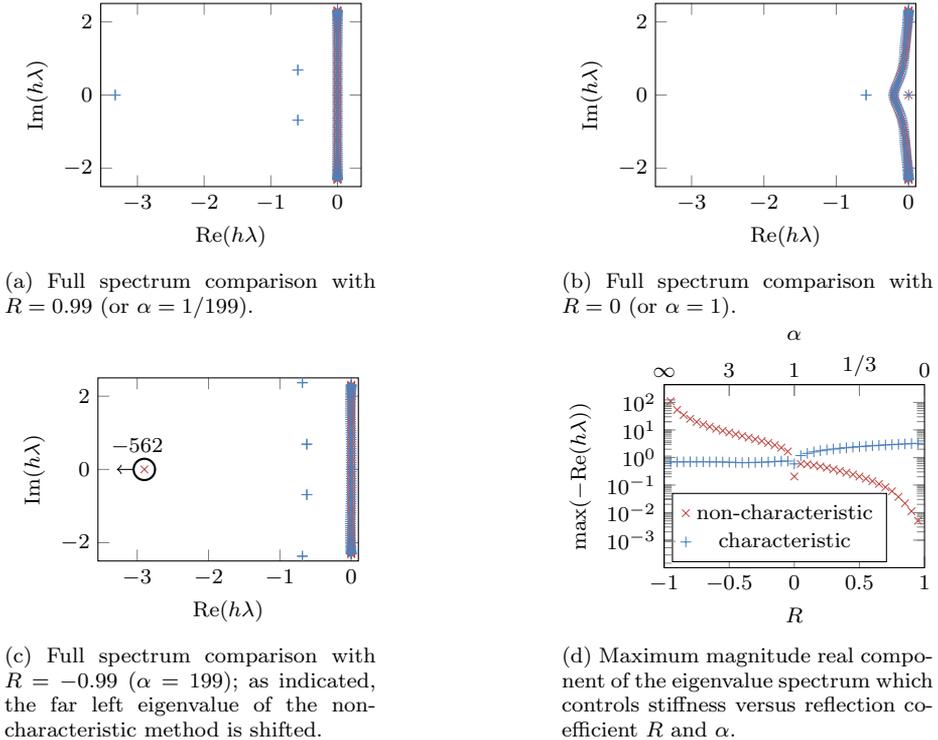

  \begin{subfigure}[t]{0.4\textwidth}
    \centering
%
     \caption{Full spectrum comparison with $R = 0.99$ (or $\alpha = 1 /
    199$).\label{fig:stiffness:1D:099}}
  \end{subfigure}
  \hfill
  \begin{subfigure}[t]{0.4\textwidth}
    \centering
%
     \caption{Full spectrum comparison with $R = 0$ (or $\alpha =
    1$).\label{fig:stiffness:1D:1}}
  \end{subfigure}
  \\
  \begin{subfigure}[t]{0.4\textwidth}
    \centering
%
     \caption{Full spectrum comparison with $R = -0.99$ ($\alpha = 199$); as
    indicated, the far left eigenvalue of the non-characteristic method is
    shifted.\label{fig:stiffness:1D:199}}
  \end{subfigure}
  \hfill
  \begin{subfigure}[t]{0.4\textwidth}
    \centering
%
     \caption{Maximum magnitude real component of the eigenvalue spectrum which
    controls stiffness versus reflection coefficient $R$ and
  $\alpha$.\label{fig:stiffness:1D:max}}
  \end{subfigure}
  \caption{Comparison of the eigenvalue spectra for the proposed characteristic
  and non-characteristic \citep{VirtaMattsson2014} treatment of boundary
  conditions for various values of reflection coefficient $R$. In all cases the
  domain is $[0, 1]$ with grid spacing $1 / 50$ and SBP interior accuracy of $2p
  = 4$. The non-characteristic method is indicated by red $\color{red} \times$ and the
  characteristic method with blue $\color{blu} +$.\label{fig:stiffness:1D}}
\end{figure}
Figures~\ref{fig:stiffness:1D:099}-\ref{fig:stiffness:1D:199} compare the
eigenvalue spectra of the operators with $R = 0.99$, $0$, and $-0.99$ (or
equivalently $\alpha = 1/99$, $1$, and $199$) using a grid with $N = 50$ and the
SBP operators with fourth order interior accuracy from \citet{Mattsson2012}.
The markers $\color{blu} +$ represent the characteristic numerical flux and
$\color{red} \times$ the non-characteristic numerical flux.
As $\alpha$ increases the spectrum associated with the non-characteristic
methods has a large magnitude negative real eigenvalue, which will severely
restrict the time step, i.e., the scheme can become arbitrarily stiff.
On the other hand, the spectrum associated with the characteristic numerical
flux is well behaved in all cases.
To further explore this, in Figure~\ref{fig:stiffness:1D:max} the maximum
magnitude real component of the eigenvalue spectra is given as a function of
$R$ (or $\alpha$).
As can be seen the maximum magnitude real eigenvalue of the non-characteristic
methods grows rapidly for $R < 0$ (or $\alpha > 1$) and is more uniform for
the characteristic method.

\begin{figure}
  \centering
  \begin{subfigure}[t]{0.4\textwidth}
    \centering
\begin{tikzpicture}[trim axis left, trim axis right]
\begin{loglogaxis}[
    ymax={1},
    ymin={1.0e-11},
    xmin={0.0009},
    xmax={0.12},
    xlabel={$h$},
    ylabel={$\|\Delta \VV{u}\|_{H}$},
    width=5cm,
    height=4cm,
    ytick={1e-1, 1e-4, 1e-7, 1e-10},
    legend pos=outer north east,
    legend cell align={left},
    legend style={draw=none},
  ]
    \addplot[color={black}, forget plot]
        table[row sep={\\}]
        {
            \\
            0.003676470588235294  0.0009338054707521556  \\
            0.001838235294117647  0.0002334513676880389  \\
            0.003676470588235294  0.0002334513676880389  \\
            0.003676470588235294  0.0009338054707521556  \\
        }
        ;
    \addplot[color={black}, forget plot]
        table[row sep={\\}]
        {
            \\
            0.003676470588235294  1.4812645010656482e-6  \\
            0.001838235294117647  9.257903131660301e-8  \\
            0.003676470588235294  9.257903131660301e-8  \\
            0.003676470588235294  1.4812645010656482e-6  \\
        }
        ;
    \addplot[color={black}, forget plot]
        table[row sep={\\}]
        {
            \\
            0.003676470588235294  2.1378816834965955e-9  \\
            0.001838235294117647  6.680880260926861e-11  \\
            0.003676470588235294  6.680880260926861e-11  \\
            0.003676470588235294  2.1378816834965955e-9  \\
        }
        ;
    \node[anchor=west] () at (axis cs:3.676471e-03, 4.669027e-04){$2$};
    \node[anchor=west] () at (axis cs:3.676471e-03, 3.703161e-07){$4$};
    \node[anchor=west] () at (axis cs:3.676471e-03, 2.672352e-10){$5$};

    \addplot[color={red}, mark={+}]
        table[row sep={\\}]
        {
            \\
            0.058823529411764705  0.22612177643808343  \\
            0.029411764705882353  0.0859543687728942  \\
            0.014705882352941176  0.025028264463180117  \\
            0.007352941176470588  0.006229654594888951  \\
            0.003676470588235294  0.001556347196596309  \\
            0.001838235294117647  0.00038903785155413304  \\
        }
        ;
        \addlegendentry{$2p = 2$: characteristic}
    \addplot[color={red}, mark={x}]
        table[row sep={\\}]
        {
            \\
            0.058823529411764705  0.23090227904292934  \\
            0.029411764705882353  0.09107555965063646  \\
            0.014705882352941176  0.024964448614618117  \\
            0.007352941176470588  0.0062289995409710724  \\
            0.003676470588235294  0.0015563424512535927  \\
            0.001838235294117647  0.0003890378168202516  \\
        }
        ;
        \addlegendentry{$2p = 2$: non-characteristic}
    \addplot[color={blu}, mark={+}]
        table[row sep={\\}]
        {
            \\
            0.058823529411764705  0.0822877307439926  \\
            0.029411764705882353  0.010554678467476425  \\
            0.014705882352941176  0.0011086930754160465  \\
            0.007352941176470588  4.778874390799815e-5  \\
            0.003676470588235294  2.5482968544872185e-6  \\
            0.001838235294117647  1.5020766207244883e-7  \\
        }
        ;
        \addlegendentry{$2p = 4$: characteristic}
    \addplot[color={blu}, mark={x}]
        table[row sep={\\}]
        {
            \\
            0.058823529411764705  0.0700386417562426  \\
            0.029411764705882353  0.006254954137733173  \\
            0.014705882352941176  0.0006354218315982737  \\
            0.007352941176470588  4.130886669612187e-5  \\
            0.003676470588235294  2.468774168442747e-6  \\
            0.001838235294117647  1.4897900560737402e-7  \\
        }
        ;
        \addlegendentry{$2p = 4$: non-characteristic}
    \addplot[color={grn}, mark={+}]
        table[row sep={\\}]
        {
            \\
            0.058823529411764705  0.18028446409985247  \\
            0.029411764705882353  0.01777050529095492  \\
            0.014705882352941176  0.0004734837362713266  \\
            0.007352941176470588  8.326214791816387e-7  \\
            0.003676470588235294  1.2135232044208862e-8  \\
            0.001838235294117647  2.2385309093063844e-10  \\
        }
        ;
        \addlegendentry{$2p = 6$: characteristic}
    \addplot[color={grn}, mark={x}]
        table[row sep={\\}]
        {
            \\
            0.058823529411764705  0.06692898953673738  \\
            0.029411764705882353  0.006162596823061943  \\
            0.014705882352941176  0.00013440534241776264  \\
            0.007352941176470588  3.644118449916648e-7  \\
            0.003676470588235294  3.5631361391609927e-9  \\
            0.001838235294117647  1.2917252971672303e-10  \\
        }
        ;
        \addlegendentry{$2p = 6$: non-characteristic}
\end{loglogaxis}
\end{tikzpicture}
     \caption{$R = 0.99$ (or $\alpha = 1 / 199$).}
  \end{subfigure}\\
  \hfill
  \begin{subfigure}[t]{0.4\textwidth}
    \centering
\begin{tikzpicture}
\begin{loglogaxis}[
    ymax={1},
    ymin={1.0e-11},
    xmin={0.0009},
    xmax={0.12},
    xlabel={$h$},
    ylabel={$\|\Delta \VV{u}\|_{H}$},
    width=5cm,
    height=4cm,
    ytick={1e-1, 1e-4, 1e-7, 1e-10},
  ]
    \addplot[color={black}]
        table[row sep={\\}]
        {
            \\
            0.003676470588235294  0.00014899605131954194  \\
            0.001838235294117647  3.7249012829885485e-5  \\
            0.003676470588235294  3.7249012829885485e-5  \\
            0.003676470588235294  0.00014899605131954194  \\
        }
        ;
    \addplot[color={black}]
        table[row sep={\\}]
        {
            \\
            0.003676470588235294  3.3177494144163815e-7  \\
            0.001838235294117647  2.0735933840102384e-8  \\
            0.003676470588235294  2.0735933840102384e-8  \\
            0.003676470588235294  3.3177494144163815e-7  \\
        }
        ;
    \addplot[color={black}]
        table[row sep={\\}]
        {
            \\
            0.003676470588235294  6.6636496086412474e-9  \\
            0.001838235294117647  2.0823905027003898e-10  \\
            0.003676470588235294  2.0823905027003898e-10  \\
            0.003676470588235294  6.6636496086412474e-9  \\
        }
        ;
    \node[anchor=west] () at (axis cs:3.676471e-03, 7.449803e-05){$2$};
    \node[anchor=west] () at (axis cs:3.676471e-03, 8.294374e-08){$4$};
    \node[anchor=west] () at (axis cs:3.676471e-03, 8.329562e-10){$5$};

    \addplot[color={red}, mark={+}]
        table[row sep={\\}]
        {
            \\
            0.058823529411764705  0.13392948879742933  \\
            0.029411764705882353  0.03270300081843555  \\
            0.014705882352941176  0.004400824918720094  \\
            0.007352941176470588  0.001015069906571923  \\
            0.003676470588235294  0.00024833662672998214  \\
            0.001838235294117647  6.173839915986293e-5  \\
        }
        ;
    \addplot[color={red}, mark={x}]
        table[row sep={\\}]
        {
            \\
            0.058823529411764705  0.1400146583195319  \\
            0.029411764705882353  0.02301081530254358  \\
            0.014705882352941176  0.004404525121558552  \\
            0.007352941176470588  0.0010152205781219284  \\
            0.003676470588235294  0.0002483267521992366  \\
            0.001838235294117647  6.17379241261392e-5  \\
        }
        ;
    \addplot[color={blu}, mark={+}]
        table[row sep={\\}]
        {
            \\
            0.058823529411764705  0.05657878404598931  \\
            0.029411764705882353  0.009838450825757763  \\
            0.014705882352941176  0.000331870298870656  \\
            0.007352941176470588  1.1656278062457311e-5  \\
            0.003676470588235294  5.841846423320025e-7  \\
            0.001838235294117647  3.3725995165327314e-8  \\
        }
        ;
    \addplot[color={blu}, mark={x}]
        table[row sep={\\}]
        {
            \\
            0.058823529411764705  0.06024308605745468  \\
            0.029411764705882353  0.004446838475684809  \\
            0.014705882352941176  0.00023370863351585732  \\
            0.007352941176470588  1.0096311580913698e-5  \\
            0.003676470588235294  5.529582357360636e-7  \\
            0.001838235294117647  3.298209570986106e-8  \\
        }
        ;
    \addplot[color={grn}, mark={+}]
        table[row sep={\\}]
        {
            \\
            0.058823529411764705  0.03174097436415166  \\
            0.029411764705882353  0.0067769968429472  \\
            0.014705882352941176  0.00031785180299593685  \\
            0.007352941176470588  6.407005747980461e-6  \\
            0.003676470588235294  1.4106648821104527e-7  \\
            0.001838235294117647  3.1324823282679175e-9  \\
        }
        ;
    \addplot[color={grn}, mark={x}]
        table[row sep={\\}]
        {
            \\
            0.058823529411764705  0.020637909292370225  \\
            0.029411764705882353  0.004958594650802197  \\
            0.014705882352941176  0.00010990796383758757  \\
            0.007352941176470588  8.453147118697981e-7  \\
            0.003676470588235294  1.1106082681068746e-8  \\
            0.001838235294117647  2.5208637680473877e-10  \\
        }
        ;
\end{loglogaxis}
\end{tikzpicture}
     \caption{$R = 0$ (or $\alpha = 1$).}
  \end{subfigure}
  \hfill
  \begin{subfigure}[t]{0.4\textwidth}
    \centering
\begin{tikzpicture}
\begin{loglogaxis}[
    ymax={1},
    ymin={1.0e-11},
    xmin={0.0009},
    xmax={0.12},
    xlabel={$h$},
    ylabel={$\|\Delta \VV{u}\|_{H}$},
    width=5cm,
    height=4cm,
    ytick={1e-1, 1e-4, 1e-7, 1e-10},
  ]
    \addplot[color={black}]
        table[row sep={\\}]
        {
            \\
            0.003676470588235294  0.000922186993184716  \\
            0.001838235294117647  0.000230546748296179  \\
            0.003676470588235294  0.000230546748296179  \\
            0.003676470588235294  0.000922186993184716  \\
        }
        ;
    \addplot[color={black}]
        table[row sep={\\}]
        {
            \\
            0.003676470588235294  1.130178075181037e-6  \\
            0.001838235294117647  7.063612969881481e-8  \\
            0.003676470588235294  7.063612969881481e-8  \\
            0.003676470588235294  1.130178075181037e-6  \\
        }
        ;
    \addplot[color={black}]
        table[row sep={\\}]
        {
            \\
            0.003676470588235294  1.6653611990816624e-8  \\
            0.001838235294117647  5.204253747130195e-10  \\
            0.003676470588235294  5.204253747130195e-10  \\
            0.003676470588235294  1.6653611990816624e-8  \\
        }
        ;
    \node[anchor=west] () at (axis cs:3.676471e-03, 4.610935e-04){$2$};
    \node[anchor=west] () at (axis cs:3.676471e-03, 2.825445e-07){$4$};
    \node[anchor=west] () at (axis cs:3.676471e-03, 2.081701e-09){$5$};

    \addplot[color={red}, mark={+}]
        table[row sep={\\}]
        {
            \\
            0.058823529411764705  0.2539814501422496  \\
            0.029411764705882353  0.10044975326074823  \\
            0.014705882352941176  0.02529182756902668  \\
            0.007352941176470588  0.006186082157190339  \\
            0.003676470588235294  0.0015382338854750114  \\
            0.001838235294117647  0.0003842771691100584  \\
        }
        ;
    \addplot[color={red}, mark={x}]
        table[row sep={\\}]
        {
            \\
            0.058823529411764705  0.3151468697171973  \\
            0.029411764705882353  0.09898277606845406  \\
            0.014705882352941176  0.02468139006585481  \\
            0.007352941176470588  0.006149674374767716  \\
            0.003676470588235294  0.0015369783219745267  \\
            0.001838235294117647  0.0003842371079368138  \\
        }
        ;
    \addplot[color={blu}, mark={+}]
        table[row sep={\\}]
        {
            \\
            0.058823529411764705  0.15602435288137395  \\
            0.029411764705882353  0.015512443396566157  \\
            0.014705882352941176  0.0006836175562566793  \\
            0.007352941176470588  3.0820828346013364e-5  \\
            0.003676470588235294  1.894485313075684e-6  \\
            0.001838235294117647  1.2122621787407574e-7  \\
        }
        ;
    \addplot[color={blu}, mark={x}]
        table[row sep={\\}]
        {
            \\
            0.058823529411764705  0.14305312497412534  \\
            0.029411764705882353  0.014421370485589561  \\
            0.014705882352941176  0.0005286586081128851  \\
            0.007352941176470588  2.9255688794290078e-5  \\
            0.003676470588235294  1.8836301253017283e-6  \\
            0.001838235294117647  1.2114909432641091e-7  \\
        }
        ;
    \addplot[color={grn}, mark={+}]
        table[row sep={\\}]
        {
            \\
            0.058823529411764705  0.18651619742638545  \\
            0.029411764705882353  0.004806496210456035  \\
            0.014705882352941176  0.0001976322705410779  \\
            0.007352941176470588  1.2100713863327298e-5  \\
            0.003676470588235294  2.9390402190866045e-7  \\
            0.001838235294117647  6.473312452570765e-9  \\
        }
        ;
    \addplot[color={grn}, mark={x}]
        table[row sep={\\}]
        {
            \\
            0.058823529411764705  0.11217626940800501  \\
            0.029411764705882353  0.005204484991306133  \\
            0.014705882352941176  7.056504910179448e-5  \\
            0.007352941176470588  1.2969897691636655e-6  \\
            0.003676470588235294  2.7756019984694373e-8  \\
            0.001838235294117647  5.277819731362915e-10  \\
        }
        ;
\end{loglogaxis}
\end{tikzpicture}
     \caption{$R = -0.99$ (or $\alpha = 199$).}
  \end{subfigure}
  \hfill~
  \caption{L$^2$-convergence comparison of the characteristic ($+$ markers) and
  non-characteristic ($\times$ markers) treatment of boundary conditions with
  various values of the reflection coefficient $R$. The red, blue, and green
  curves correspond to SBP interior orders $2$, $4$, and $6$,
  respectively.\label{fig:convergence:1D}}
\end{figure}
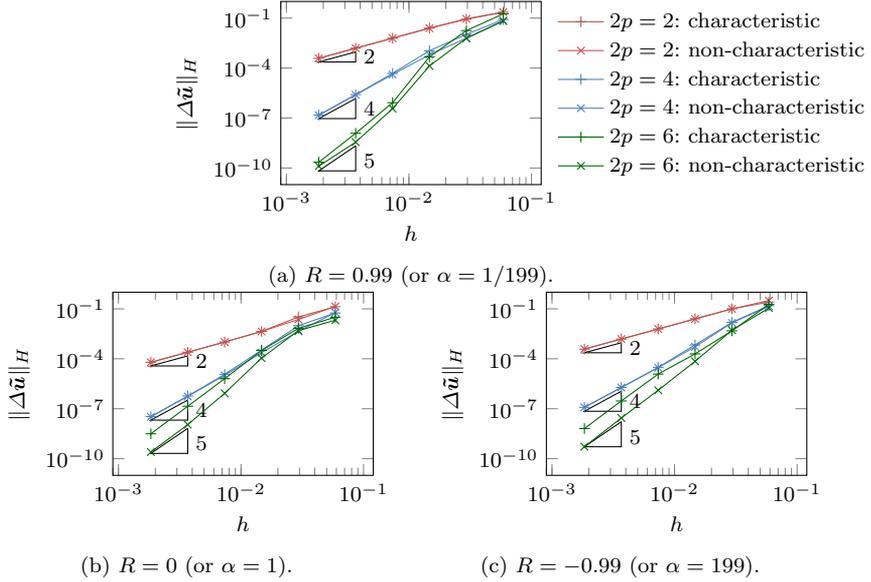
We now considering the accuracy of both the non-characteristic and
characteristic schemes.
In this test we consider the initial condition
\begin{subequations}
  \begin{align}
    u_{0}(x) &= \sin(2 \pi x)^6,\\
    \dot{u}_{0}(x) &= 0,
  \end{align}
\end{subequations}
which for times $t \in [0, 1]$ and reflection coefficients $R\in[-1, 1]$  has
the analytic solution
\begin{equation}
  \begin{split}
    u(x, t) &=
    \frac{\bar{u}_{0}(x - t) + \bar{u}_{0}(x + t) + R\left(\bar{u}_{0}(2 - x -
    t) + \bar{u}_{0}(-x + t)\right)}{2},\\
    \bar{u}_{0}(x) &=
    \begin{cases}
      u_{0}(x), & 0 \le x \le 1,\\
      0, & \text{otherwise.}
    \end{cases}
  \end{split}
\end{equation}
The L$^{2}$-convergence can be seen in Figure~\ref{fig:convergence:1D} at time
$t = 0.9$ for $R = 0.99$, $0$, and $-0.99$; the L$^2$-error in the solution is
defined as
\begin{equation}
  \|\Delta \vv{u}\|_{H} = \sqrt{\Delta \vv{u}^{T} \mm{H} \Delta \vv{u}},
\end{equation}
with $\Delta \vv{u}$ being the pointwise difference between the numerical and
exact solutions.
In the figure, the interior accuracy of the SBP scheme is denoted by $2p$ and
all the operators are from~\citet{Mattsson2012}.
The spatial resolutions used in the test are $N = 17 \times 2^r$ with $r = 0, 1,
2, 3, 4, 5$ and time integration is performed using matrix exponentiation.
As can be seen, both the non-characteristic and characteristic methods converge
at similar rates.
For the characteristic method with $2p = 6$ the overall error constant is
higher, though this can be improved by increasing the penalty parameter (not
shown) at the cost of increased stiffness.

\subsection{Nonlinear Interface}

\subsubsection{Continuous Problem}

We now consider an interface at $x = 0$ on which the traction acting on the
interface is related to the jump in particle velocity across the interface.
To do this, we consider the domain $\Omega = [-1, 1]$ and modify wave
equation~\eqref{eqn:1d:wave} to be
\begin{subequations}\label{eqn:1d:interface}
  \begin{alignat}{3}
    \label{eqn:1d:interface:pde}
    &\ddot{u} = \partial_{1}^2 u, &&~~x \in \Omega ,&&~t \in [0, T],\\
    \label{eqn:1d:interface:bc}
    &\tau = -\alpha \dot{u}, &&~ ~x \in \{-1, 1\},&&~t \in [0, T],\\
    \label{eqn:1d:interface:condition}
    &\begin{cases}
      \tau^{+} = -\tau^{-}\\
      \tau^{-} = F\left(\dot{u}^{+} - \dot{u}^{-}\right)
    \end{cases},
    &&~~x = 0 ,&&~t \in [0, T].
  \end{alignat}
\end{subequations}
In interface condition~\eqref{eqn:1d:interface:condition} the superscript $\pm$
denotes the solution on the two-sides of the interface with $\dot{u}^{-}$ being
the particle velocity as $x \rightarrow 0^{-}$ and $\dot{u}^{+}$ the particle
velocities as $x \rightarrow 0^{+}$.
The tractions $\tau^{\pm}$ are defined with normals that point out of the
respective sides of the interface, i.e., $\tau^{\pm} = \mp \partial_{1}
u^{\pm}$.
The nonlinear interface function $F$ is assumed to be odd and take the same sign
as its argument:
\begin{equation}
  V\;F(V) \ge 0
  \quad \mbox{and} \quad
  F(-V) = -F(V).
\end{equation}

Defining the solution energy,
\begin{equation}
  \label{eqn:1d:interface:energy}
  E = \frac{1}{2} \int_{-1}^{1}
  \left(
    \dot{u}^2 + {\left(\partial_{1} u\right)}^2
  \right),
\end{equation}
leads to the following energy estimate and thus well-posedness:
\begin{lemma}\label{lemma:1d:interface:energy:estimate}
  Governing equations~\eqref{eqn:1d:interface} with
  energy~\eqref{eqn:1d:interface:energy} satisfies $\dot{E} \le 0$.
\end{lemma}
\begin{proof}
  Taking the time derivative of energy~\eqref{eqn:1d:interface:energy},
  using scalar wave equation~\eqref{eqn:1d:interface:pde}, and simplifying
  with integration by parts yields
  \begin{equation}
    \dot{E}
    =
    \int_{-1}^{1}
    \left(
      \dot{u}\left(\partial_{1}^2 u\right) + \left(\partial_{1} u\right) \left(\partial_{1} \dot{u}\right)
    \right)
    =
    {\left.\dot{u}\;\partial_{1} u\right|}_{-1}^{0^{-}}
      +
      {\left.\dot{u}\;\partial_{1} u\right|}_{0^{+}}^{1}.
  \end{equation}
  Applying the definition of traction~\eqref{eqn:1d:wave:tau}, boundary
  conditions~\eqref{eqn:1d:interface:bc}, and interface
  conditions~\eqref{eqn:1d:interface:condition} leads to
  \begin{equation}
    \dot{E}
    =
    -{\left.\alpha\dot{u}^{2} \right|}_{1}
    -{\left.\alpha\dot{u}^{2} \right|}_{0}
      -\left(\dot{u}^{+} - \dot{u}^{-} \right) F\left(\dot{u}^{+} - \dot{u}^{-} \right)
    \le 0,
  \end{equation}
  since $\alpha \ge 0$ and $V\;F(V) \ge 0$.
  \qed
\end{proof}

\subsubsection{Characteristic Nonlinear Interface Condition}\label{sec:1d:interface:char}

As was done with boundary conditions in \sref{sec:1d:bc:continuous}, it will be
useful for the discretization to rewrite nonlinear interface
condition~\eqref{eqn:1d:interface:condition} in terms of the characteristic
variables propagating into and out of the interface.
Namely we let $w^{\pm}$ be the characteristic variables propagating into and
$q^{\pm}$ out of the two sides of interface:
\begin{equation}
  w^{\pm} = \dot{u}^{\pm} + \tau^{\pm}
  \quad\mbox{and}\quad
  q^{\pm} = \dot{u}^{\pm} - \tau^{\pm}.
\end{equation}
We then define the nonlinear functions $\mathcal{Q}^{\pm}(w^{-}, w^{+})$ so that
the characteristic-defined interface particle velocities and tractions,
\begin{subequations}
  \begin{align}
    \dot{u}^{\pm} &= \frac{\mathcal{Q}^{\pm} + w^{\pm}}{2},\\
    \tau^{\pm} &= \frac{\mathcal{Q}^{\pm} - w^{\pm}}{2},
  \end{align}
\end{subequations}
satisfy nonlinear interface condition~\eqref{eqn:1d:interface:condition}.
For a general $F$ there is no closed form expression of $\mathcal{Q}^{\pm}$, but
existence can be guaranteed by the implicit function theorem if $F'(V) > 0$
\citep[Proposition 1]{KozdonDunhamNordstrom2012}; \aref{app:rootfinding} show
how determining $\mathcal{Q}^{\pm}$ can be reduced to a scalar root finding
problem.

\subsubsection{Discrete Problem}
The interface problem is discretized with two blocks each with $N + 1$ grid
points.
The grid solution on $[-1, 0]$ is $\vv{u}^{-}$ and the grid solution on $[0, 1]$
is $\vv{u}^{+}$.
With this a semidiscretization in each block is then
\begin{equation}
  \label{eqn:1d:interface:disc}
  \ddot{\vv{u}}^{\pm}
  = \mm{D}_{11} \vv{u}^{\pm}
  +
  \sum_{\mathclap{k\in\{0,N\}}}\;
  \left(
    \mm{H}^{-1}\vv{e}_{k}\left(\tau^{\pm *}_{k} - n_{k} \vv{b}_{k}^{T}
    \vv{u}^{\pm}\right)
    - n_{k}\mm{H}^{-1}\vv{b}_{k}\left(u^{\pm *}_{k} - u_{k}^{\pm}\right)
  \right),
\end{equation}
where $n_{0} = -1$ and $n_{N} = 1$ are the outward pointing normals at the block
edges.
As in section \sref{sec:1d:bc}, the terms $\tau^{\pm*}_{k}$ and $u^{\pm *}_{k}$ are yet-to-be-defined
numerical fluxes which can either be set in a non-characteristic or
characteristic manner.

In each block we define the energy in the solution as
\begin{equation}
  \label{eqn:1d:interface:disc:block:energy}
  \begin{split}
    \mathcal{E}^{\pm} =&\; \frac{1}{2}
    \left(
      {\left(\dot{\vv{u}}^{\pm}\right)}^{T} \mm{H} \dot{\vv{u}}^{\pm}
      +
      {\left(\vv{u}^{\pm}\right)}^{T} \mm{A}_{11} \vv{u}^{\pm}
    \right)\\
                       &+
                       \frac{1}{2\gamma}\sum_{{k \in \{0, N\}}}\left(
      {\left(\tau_{k}^{\pm}\right)}^{2}
      -
      {\left(\vv{b}_{k}^{T}\vv{u}^{\pm}\right)}^{2}
    \right).
  \end{split}
\end{equation}
with traction $\tau_{k}^{\pm}$ being as defined~\eqref{eqn:1d:wave:disc:tau}
for a penalty parameter $\gamma$.
The energy in the whole domain is then
\begin{equation}
  \label{eqn:1d:interface:disc:energy}
  \mathcal{E} = \mathcal{E}^{-} + \mathcal{E}^{+},
\end{equation}
and the following theorem guarantees that $\mathcal{E}$ is a seminorm of the
solution.
\begin{theorem}\label{thm:interface:disc:1d:energy:seminorm}
  Energy~\eqref{eqn:1d:interface:disc:energy} is a seminorm for all
  $\vv{u}^{\pm}$ and $u^{\pm *}_{k}$ if $\gamma$ is positive and sufficiently
  large.
\end{theorem}
\begin{proof}
  Since the block energy~\eqref{eqn:1d:interface:disc:block:energy} is of the
  same form as~\eqref{eqn:1d:wave:disc:energy} the result follows directly from
  Theorem~\ref{thm:disc:1d:energy:seminorm}.
  \qed
\end{proof}
\begin{corollary}
  Semidiscretization~\eqref{eqn:1d:interface:disc} satisfies
  \begin{equation}
    \label{eqn:1d:interface:disc:energy:rate}
    \dot{\mathcal{E}} = BT + IT
  \end{equation}
  where $\mathcal{E}$ is defined by~\eqref{eqn:1d:interface:disc:energy} with the
  boundary term $BT$ and interface term $IT$ being
  \begin{alignat}{3}
    BT =&\;
    \dot{u}_{0}^{- *}\tau_{0}^{- *}
    -
    \left(\dot{u}_{0}^{-} - \dot{u}^{- *}_{0}\right)
    \left(\tau_{0}^{-} - \tau^{- *}_{0}\right)\\\notag&
    +
    \dot{u}_{N}^{+ *}\tau_{N}^{+ *}
    -
    \left(\dot{u}_{N}^{+} - \dot{u}^{+ *}_{N}\right)
    \left(\tau_{N}^{+} - \tau^{+ *}_{N}\right), \\
    \label{eqn:1d:interface:IT}
    IT =&\;
    \dot{u}_{N}^{- *}\tau_{N}^{- *}
    -
    \left(\dot{u}_{N}^{-} - \dot{u}^{- *}_{N}\right)
    \left(\tau_{N}^{-} - \tau^{- *}_{N}\right)\\\notag&
    +
    \dot{u}_{0}^{+ *}\tau_{0}^{+ *}
    -
    \left(\dot{u}_{0}^{+} - \dot{u}^{+ *}_{0}\right)
    \left(\tau_{0}^{+} - \tau^{+ *}_{0}\right).
  \end{alignat}
\end{corollary}
\begin{proof}
  Follows directly by taking the time derivative of
  energy~\eqref{eqn:1d:interface:disc:energy}, substituting
  semidiscretization~\eqref{eqn:1d:interface:disc}, and applying SBP
  property~\eqref{eqn:sbp:d2:1d}.
  \qed
\end{proof}

As shown in \sref{sec:1d:bc}, if boundary conditions are enforced with either the
non-characteristic~\eqref{eqn:1d:wave:disc:nonchar:bc:flux} or
characteristic~\eqref{eqn:1d:wave:disc:char:bc:flux} numerical fluxes
then $BT \le 0$.
All that remains to be shown is how the numerical fluxes for the interface can
be defined so that $IT \le 0$, and as with the boundaries we will show how this
can be done in a non-characteristic and characteristic manner.

\subsubsection{Non-Characteristic Interface Treatment}

In \citet{Duru2019} it was proposed to set the numerical fluxes for the
interface as
\begin{subequations}
  \label{eqn:1d:interface:nc:flux}
  \begin{align}
    \tau^{-*}_{0} &= -\tau^{+*}_{N} = F\left(\dot{u}_{0}^{+} - \dot{u}_{N}^{-}\right),\\
    u^{-*}_{0} &= u_{0}^{-},\\
    u^{+*}_{N} &= u_{N}^{+}.
  \end{align}
\end{subequations}
This leads to the following energy estimate for the interface:
\begin{theorem}
  Non-characteristic numerical flux~\eqref{eqn:1d:interface:nc:flux} leads to
  the energy dissipation $IT \le 0$ where $IT$ is defined
  by~\eqref{eqn:1d:interface:IT}.
\end{theorem}
\begin{proof}
  Substituting numerical flux~\eqref{eqn:1d:interface:nc:flux}
  into~\eqref{eqn:1d:interface:IT} gives
  \begin{equation}
    IT =
    \dot{u}_{N}^{- *}\tau_{N}^{- *}
    +
    \dot{u}_{0}^{+ *}\tau_{0}^{+ *}
    =
    -\left(\dot{u}_{0}^{+ *} - \dot{u}_{N}^{- *}\right)
    F\left(\dot{u}_{0}^{+ *} - \dot{u}_{N}^{- *}\right),
  \end{equation}
  and the result follows since $V\;F(V) \ge 0$.
  \qed
\end{proof}

\subsubsection{Characteristic Interface Treatment}

The characteristic numerical interface fluxes are defined as
\begin{subequations}
  \label{eqn:1d:interface:char:flux}
  \begin{align}
    \dot{u}_{N}^{-*} &=
    \frac{q_{N}^{-*} + w_{N}^{-*}}{2},&
    \dot{u}_{0}^{+*} &=
    \frac{q_{0}^{+*} + w_{0}^{+*}}{2},\\
    \tau_{N}^{-*} &=
    \frac{q_{N}^{-*} - w_{N}^{-*}}{2},&
    \tau_{0}^{+*} &=
    \frac{q_{0}^{+*} - w_{0}^{+*}}{2},
  \end{align}
\end{subequations}
where the characteristic variables are
\begin{subequations}
  \begin{align}
    w^{-*}_{N} &= \frac{\dot{u}_{N}^{-} - \tau_{N}^{-}}{2},&
    w^{+*}_{0} &= \frac{\dot{u}_{0}^{+} - \tau_{0}^{+}}{2},\\
    q^{-*}_{N} &= \mathcal{Q}^{-}\left(w^{-*}_{N}, w^{+*}_{0}\right),&
    q^{+*}_{0} &= \mathcal{Q}^{+}\left(w^{-*}_{N}, w^{+*}_{0}\right);
  \end{align}
\end{subequations}
see \sref{sec:1d:interface:char}.
As with the case of characteristic treatment of boundary conditions, the
numerical fluxes $u_{N}^{-*}$ and $u_{0}^{+*}$ must be tracked as variables
during the solution process.
The following theorem guarantees energy stability of the characteristic
interface treatment:
\begin{theorem}
  Characteristic numerical flux~\eqref{eqn:1d:interface:char:flux} leads to the
  energy dissipation $IT \le 0$ where $IT$ is defined
  by~\eqref{eqn:1d:interface:IT}.
\end{theorem}
\begin{proof}
  By definition, characteristic numerical
  flux~\eqref{eqn:1d:interface:char:flux} satisfies the nonlinear
  interface condition:
  \begin{equation}
    \tau^{-*}_{N} = - \tau^{+*}_{0} =
    F\left(\dot{u}^{+*}_{0} - \dot{u}^{-*}_{N}\right).
  \end{equation}
  With this, we then have that
  \begin{equation}
    \label{eqn:1d:interface:fric}
    \dot{u}_{N}^{- *}\tau_{N}^{- *} + \dot{u}_{0}^{+ *}\tau_{0}^{+ *}
    =
    -\left(\dot{u}_{0}^{+ *} - \dot{u}_{N}^{- *}\right)
    F\left(\dot{u}^{+*}_{0} - \dot{u}^{-*}_{N}\right)
    \le 0,
  \end{equation}
  since $V\;F(V) \ge 0$. Defining the grid based interface outgoing characteristic
  variables
  \begin{subequations}
    \begin{align}
      q^{-}_{N} &= \frac{\dot{u}_{N}^{-} + \tau_{N}^{-}}{2},&
      q^{+}_{0} &= \frac{\dot{u}_{0}^{+} + \tau_{0}^{+}}{2},
    \end{align}
  \end{subequations}
  it follows from characteristic numerical
  flux~\eqref{eqn:1d:interface:char:flux} that
  \begin{subequations}
    \label{eqn:1d:interface:diff}
    \begin{align}
      \dot{u}_{N}^{-} - \dot{u}_{N}^{-*} = \tau_{N}^{-} - \tau_{N}^{-*} =
      \frac{q_{N}^{-} - q_{N}^{-*}}{2},\\
      \dot{u}_{0}^{+} - \dot{u}_{0}^{+*} = \tau_{0}^{+} - \tau_{0}^{+*} =
      \frac{q_{0}^{+} - q_{0}^{+*}}{2}.
    \end{align}
  \end{subequations}
  Using relations~\eqref{eqn:1d:interface:fric}
  and~\eqref{eqn:1d:interface:diff} in interface
  term~\eqref{eqn:1d:interface:IT} gives
  \begin{equation}
    IT =
    -\left(\dot{u}_{0}^{+ *} - \dot{u}_{N}^{- *}\right)
    F\left(\dot{u}^{+*}_{0} - \dot{u}^{-*}_{N}\right)
    -\frac{{\left(q_{N}^{-} - q_{N}^{-*}\right)}^{2}}{4}
    -\frac{{\left(q_{0}^{+} - q_{0}^{+*}\right)}^{2}}{4}
    \le 0.
  \end{equation}
  \qed
\end{proof}

\subsubsection{Numerical Results}

In order to assess the stiffness and accuracy of the non-characteristic and
characteristic interface treatment, we use the nonlinear interface function
\begin{equation}
  F(V) = \beta \arcsinh\left(V\right),
\end{equation}
where the parameter $\beta > 0$ will vary in the experiment.
The outer boundary conditions are Neumann and enforced using the
non-characteristic method.
The initial condition is taken as
\begin{subequations}
  \begin{align}
    u_{0}(x) &= \exp\left(-{\left(\frac{x - \mu}{\sigma}\right)}^{2}\right),\\
    \dot{u}_{0}(x) &= \frac{2(x-\mu)}{\sigma^2} u_0(x),
  \end{align}
\end{subequations}
with $\mu = -1/2$ and $\sigma = 1/15$.
These initial conditions lead to a pulse centered at $x = \mu$ moving to the
right.
The exact solution for displacement is
\begin{equation}
  u(x, t) =
  \begin{cases}
    u_{0}(x-t) + \psi^{-}\left(t + x\right), & x < 0,\\
    u_{0}(x+t) + \psi^{+}\left(t - x\right), & x > 0.
  \end{cases}
\end{equation}
The functions $\psi^{\pm}$ are a result of the wave propagating out of the
interface, and are defined as
\begin{equation}
  \label{eqn:1d:wave:semi:int}
  \psi^{\pm}(T) =
  \frac{1}{2}
  \begin{cases}
    0, & T \le 0\\
    \int_{0}^{T}
    \mathcal{Q}^{\pm}\left(
      w^{-}_{0}\left(t\right),
      w^{+}_{0}\left(t\right)
    \right) \; dt,& T>0,
  \end{cases}
\end{equation}
with the characteristic variables propagating into the interface being
\begin{equation}
  w^{\pm}_{0}\left(t\right)
  =
  \dot{u}_{0}\left(\pm t\right)
  \pm
  u'_{0}\left(\pm t\right);
\end{equation}
the integral in~\eqref{eqn:1d:wave:semi:int} must be solved numerically since
$\mathcal{Q}^{\pm}$ does not have a closed form.

Time stepping is done by converting
semidiscretization~\eqref{eqn:1d:interface:disc} into a first order system and
using the low-storage, fourth order Runge-Kutta scheme of \citet[(5,4)
$2N$-Storage RK scheme, solution $3$]{CarpenterKennedy1994}.
The time step size is selected to be of the form
\begin{equation}
  \Delta t = \kappa\;h,
\end{equation}
where $\kappa$ is a Courant number chosen so that the simulation is stable and
accurate.
The L$^2$-error in the solution is
\begin{equation}
  \|\Delta \vv{u}\|_{H} = \sqrt{
    {\left(\Delta \vv{u}^{-}\right)}^{T} \mm{H} \Delta \vv{u}^{-} +
    {\left(\Delta \vv{u}^{+}\right)}^{T} \mm{H} \Delta \vv{u}^{+}
  },
\end{equation}
with $\Delta \vv{u}^{\pm}$ being the pointwise difference between the numerical
and exact solutions on the two sides of the interface.
To find a suitable value of $\kappa$ we start with an initial value of $\kappa =
1$ and this is successively halved until the simulation is stable and accurate.
\begin{table}
  \centering
  \begin{tabular}{rcccccc}
    \toprule
    & \multicolumn{3}{c}{characteristic} & \multicolumn{3}{c}{non-characteristic}\\
    $\beta$ & $2p = 2$ & $2p = 4$ & $2p = 6$ & $2p = 2$ & $2p = 4$ & $2p = 6$\\
    \midrule
    $32$  & $1 / 2$ & $1 / 2$ & $1 / 4$ & $1 /  32$ & $1 /  64$ & $1 /  64$\\
    $64$  & $1 / 2$ & $1 / 2$ & $1 / 4$ & $1 /  64$ & $1 / 128$ & $1 / 128$\\
    $128$ & $1 / 2$ & $1 / 2$ & $1 / 4$ & $1 / 128$ & $1 / 256$ & $1 / 256$\\
    \bottomrule
  \end{tabular}
  \caption{Stable Courant $\kappa$ for the characteristic and non-characteristic
    interface treatment for increasing values of $\beta$ and various SBP
    operator with interior accuracy $2p$.\label{tab:stiffness:1D}}
\end{table}
Table~\ref{tab:stiffness:1D} gives the Courant $\kappa$ determined for each SBP
order and $\beta$ pair.
As can be seen the Courant number for the characteristic method is purely a
function of the SBP operator and for the non-characteristic the time step size
scales inversely with $\beta$.

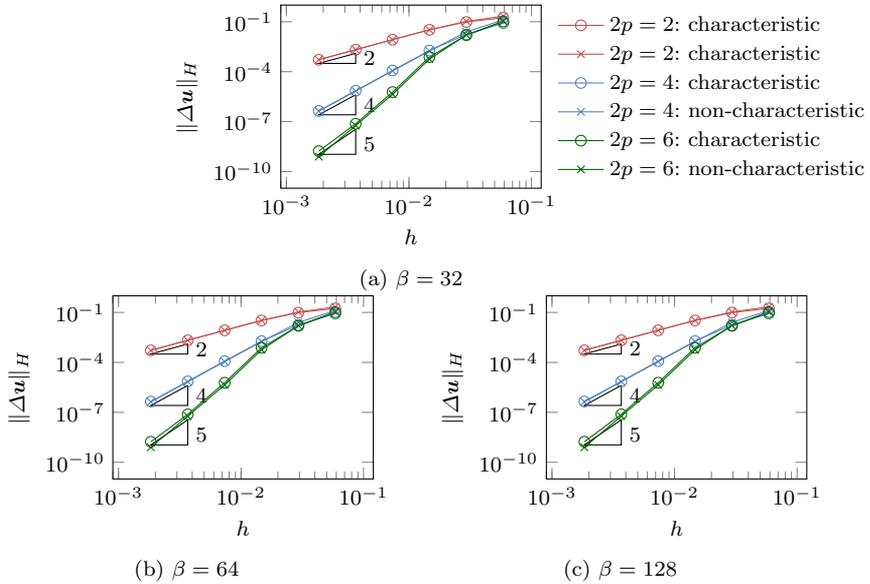
\begin{figure}
  \centering
  \begin{subfigure}[t]{0.4\textwidth}
    \centering
\begin{tikzpicture}[trim axis left, trim axis right]
\begin{loglogaxis}[
  ymax={1},
  ymin={1.0e-11},
  xmin={0.0009},
  xmax={0.12},
  xlabel={$h$},
  ylabel={$\|\Delta \vv{u}\|_{H}$},
  width=5cm,
  height=4cm,
  ytick={1e-1, 1e-4, 1e-7, 1e-10},
  legend pos=outer north east,
  legend cell align={left},
  legend style={draw=none},
  ]
    \addplot[color={black}, forget plot]
        table[row sep={\\}]
        {
            \\
            0.003676470588235294  0.0012570048041047526  \\
            0.001838235294117647  0.00031425120102618814  \\
            0.003676470588235294  0.00031425120102618814  \\
            0.003676470588235294  0.0012570048041047526  \\
        }
        ;
        \addplot[color={black}, forget plot]
        table[row sep={\\}]
        {
            \\
            0.003676470588235294  4.1036725889692965e-6  \\
            0.001838235294117647  2.5647953681058103e-7  \\
            0.003676470588235294  2.5647953681058103e-7  \\
            0.003676470588235294  4.1036725889692965e-6  \\
        }
        ;
        \addplot[color={black}, forget plot]
        table[row sep={\\}]
        {
            \\
            0.003676470588235294  3.474082877985275e-8  \\
            0.001838235294117647  1.0856508993703985e-9  \\
            0.003676470588235294  1.0856508993703985e-9  \\
            0.003676470588235294  3.474082877985275e-8  \\
        }
        ;
    \node[anchor=west] () at (axis cs:3.676471e-03, 6.285024e-04){$2$};
    \node[anchor=west] () at (axis cs:3.676471e-03, 1.025918e-06){$4$};
    \node[anchor=west] () at (axis cs:3.676471e-03, 4.342604e-09){$5$};

    \addplot[color={red}, mark={o}]
        table[row sep={\\}]
        {
            \\
            0.058823529411764705  0.17320523136869984  \\
            0.029411764705882353  0.09629908815902914  \\
            0.014705882352941176  0.03275037571675725  \\
            0.007352941176470588  0.008410111199629299  \\
            0.003676470588235294  0.002096968264744474  \\
            0.001838235294117647  0.0005238801657876518  \\
        }
        ;
        \addlegendentry{$2p = 2$: characteristic}
    \addplot[color={red}, mark={x}]
        table[row sep={\\}]
        {
            \\
            0.058823529411764705  0.21609094826621897  \\
            0.029411764705882353  0.10616654480867492  \\
            0.014705882352941176  0.03270257292610321  \\
            0.007352941176470588  0.008366729335484118  \\
            0.003676470588235294  0.0020950080068412543  \\
            0.001838235294117647  0.0005238081035225379  \\
        }
        ;
        \addlegendentry{$2p = 2$: characteristic}
    \addplot[color={blu}, mark={o}]
        table[row sep={\\}]
        {
            \\
            0.058823529411764705  0.08338313775672741  \\
            0.029411764705882353  0.017212753521918915  \\
            0.014705882352941176  0.0018667131563323328  \\
            0.007352941176470588  0.00011782443978974617  \\
            0.003676470588235294  7.130354136137796e-6  \\
            0.001838235294117647  4.407723775824383e-7  \\
        }
        ;
        \addlegendentry{$2p = 4$: characteristic}
    \addplot[color={blu}, mark={x}]
        table[row sep={\\}]
        {
            \\
            0.058823529411764705  0.12821997704526872  \\
            0.029411764705882353  0.024775139519025588  \\
            0.014705882352941176  0.0017598624806215914  \\
            0.007352941176470588  0.00011007029628041189  \\
            0.003676470588235294  6.839454314948827e-6  \\
            0.001838235294117647  4.246717759267561e-7  \\
        }
        ;
        \addlegendentry{$2p = 4$: non-characteristic}
    \addplot[color={grn}, mark={o}]
        table[row sep={\\}]
        {
            \\
            0.058823529411764705  0.10258266027018158  \\
            0.029411764705882353  0.015505043092841963  \\
            0.014705882352941176  0.000756679737431766  \\
            0.007352941176470588  5.973112820469564e-6  \\
            0.003676470588235294  7.456693566331302e-8  \\
            0.001838235294117647  1.6831186892759116e-9  \\
        }
        ;
        \addlegendentry{$2p = 6$: characteristic}
    \addplot[color={grn}, mark={x}]
        table[row sep={\\}]
        {
            \\
            0.058823529411764705  0.10648899006483742  \\
            0.029411764705882353  0.01591622530155028  \\
            0.014705882352941176  0.0005689943817070845  \\
            0.007352941176470588  4.398137245611646e-6  \\
            0.003676470588235294  5.7901381299754586e-8  \\
            0.001838235294117647  7.844671146686644e-10  \\
        }
        ;
        \addlegendentry{$2p = 6$: non-characteristic}
\end{loglogaxis}
\end{tikzpicture}
     \caption{$\beta = 32$}
  \end{subfigure}\\
  \hfill
  \begin{subfigure}[t]{0.4\textwidth}
    \centering
\begin{tikzpicture}
\begin{loglogaxis}[
  ymax={1},
  ymin={1.0e-11},
  xmin={0.0009},
  xmax={0.12},
  xlabel={$h$},
  ylabel={$\|\Delta \vv{u}\|_{H}$},
  width=5cm,
  height=4cm,
  ytick={1e-1, 1e-4, 1e-7, 1e-10},
  ]
    \addplot[color={black}]
        table[row sep={\\}]
        {
            \\
            0.003676470588235294  0.0012669739902238884  \\
            0.001838235294117647  0.0003167434975559721  \\
            0.003676470588235294  0.0003167434975559721  \\
            0.003676470588235294  0.0012669739902238884  \\
        }
        ;
    \addplot[color={black}]
        table[row sep={\\}]
        {
            \\
            0.003676470588235294  4.132528775097858e-6  \\
            0.001838235294117647  2.5828304844361614e-7  \\
            0.003676470588235294  2.5828304844361614e-7  \\
            0.003676470588235294  4.132528775097858e-6  \\
        }
        ;
    \addplot[color={black}]
        table[row sep={\\}]
        {
            \\
            0.003676470588235294  3.525459415910581e-8  \\
            0.001838235294117647  1.1017060674720566e-9  \\
            0.003676470588235294  1.1017060674720566e-9  \\
            0.003676470588235294  3.525459415910581e-8  \\
        }
        ;
    \node[anchor=west] () at (axis cs:3.676471e-03, 6.334870e-04){$2$};
    \node[anchor=west] () at (axis cs:3.676471e-03, 1.033132e-06){$4$};
    \node[anchor=west] () at (axis cs:3.676471e-03, 4.406824e-09){$5$};

    \addplot[color={red}, mark={o}]
        table[row sep={\\}]
        {
            \\
            0.058823529411764705  0.174153018704673  \\
            0.029411764705882353  0.09693123982588336  \\
            0.014705882352941176  0.033001316040755395  \\
            0.007352941176470588  0.008477006982827566  \\
            0.003676470588235294  0.00211362781169718  \\
            0.001838235294117647  0.000528042247463518  \\
        }
        ;
    \addplot[color={red}, mark={x}]
        table[row sep={\\}]
        {
            \\
            0.058823529411764705  0.2175361964934934  \\
            0.029411764705882353  0.10694946418345153  \\
            0.014705882352941176  0.03295545764856463  \\
            0.007352941176470588  0.008432698392086773  \\
            0.003676470588235294  0.002111623317039814  \\
            0.001838235294117647  0.0005279686349854534  \\
        }
        ;
    \addplot[color={blu}, mark={o}]
        table[row sep={\\}]
        {
            \\
            0.058823529411764705  0.08382702402982761  \\
            0.029411764705882353  0.017427879233179962  \\
            0.014705882352941176  0.001886559169601905  \\
            0.007352941176470588  0.00011862662302394476  \\
            0.003676470588235294  7.180243199963409e-6  \\
            0.001838235294117647  4.440504840899452e-7  \\
        }
        ;
    \addplot[color={blu}, mark={x}]
        table[row sep={\\}]
        {
            \\
            0.058823529411764705  0.12929060299761305  \\
            0.029411764705882353  0.02509590346271453  \\
            0.014705882352941176  0.0017774342902201104  \\
            0.007352941176470588  0.00011083617956641886  \\
            0.003676470588235294  6.887547958496431e-6  \\
            0.001838235294117647  4.278272889667239e-7  \\
        }
        ;
    \addplot[color={grn}, mark={o}]
        table[row sep={\\}]
        {
            \\
            0.058823529411764705  0.10333650348711322  \\
            0.029411764705882353  0.015556865214212944  \\
            0.014705882352941176  0.0007576730960345756  \\
            0.007352941176470588  5.998560329966968e-6  \\
            0.003676470588235294  7.56287298349577e-8  \\
            0.001838235294117647  1.702222130590974e-9  \\
        }
        ;
    \addplot[color={grn}, mark={x}]
        table[row sep={\\}]
        {
            \\
            0.058823529411764705  0.10753032961288433  \\
            0.029411764705882353  0.016137688485768888  \\
            0.014705882352941176  0.0005700479109706865  \\
            0.007352941176470588  4.40786052748971e-6  \\
            0.003676470588235294  5.875765693184302e-8  \\
            0.001838235294117647  7.954482339561575e-10  \\
        }
        ;
\end{loglogaxis}
\end{tikzpicture}
     \caption{$\beta = 64$}
  \end{subfigure}
  \hfill
  \begin{subfigure}[t]{0.4\textwidth}
    \centering
\begin{tikzpicture}
\begin{loglogaxis}[
  ymax={1},
  ymin={1.0e-11},
  xmin={0.0009},
  xmax={0.12},
  xlabel={$h$},
  ylabel={$\|\Delta \vv{u}\|_{H}$},
  width=5cm,
  height=4cm,
  ytick={1e-1, 1e-4, 1e-7, 1e-10},
  ]
    \addplot[color={black}]
        table[row sep={\\}]
        {
            \\
            0.003676470588235294  0.0012719020209180506  \\
            0.001838235294117647  0.00031797550522951265  \\
            0.003676470588235294  0.00031797550522951265  \\
            0.003676470588235294  0.0012719020209180506  \\
        }
        ;
    \addplot[color={black}]
        table[row sep={\\}]
        {
            \\
            0.003676470588235294  4.146840077404588e-6  \\
            0.001838235294117647  2.5917750483778676e-7  \\
            0.003676470588235294  2.5917750483778676e-7  \\
            0.003676470588235294  4.146840077404588e-6  \\
        }
        ;
    \addplot[color={black}]
        table[row sep={\\}]
        {
            \\
            0.003676470588235294  3.551712274635052e-8  \\
            0.001838235294117647  1.1099100858234538e-9  \\
            0.003676470588235294  1.1099100858234538e-9  \\
            0.003676470588235294  3.551712274635052e-8  \\
        }
        ;
    \node[anchor=west] () at (axis cs:3.676471e-03, 6.359510e-04){$2$};
    \node[anchor=west] () at (axis cs:3.676471e-03, 1.036710e-06){$4$};
    \node[anchor=west] () at (axis cs:3.676471e-03, 4.439640e-09){$5$};

    \addplot[color={red}, mark={o}]
        table[row sep={\\}]
        {
            \\
            0.058823529411764705  0.17462240328513043  \\
            0.029411764705882353  0.0972441961654238  \\
            0.014705882352941176  0.03312526659095378  \\
            0.007352941176470588  0.008510068041040876  \\
            0.003676470588235294  0.002121863025526836  \\
            0.001838235294117647  0.0005300996932284901  \\
        }
        ;
    \addplot[color={red}, mark={x}]
        table[row sep={\\}]
        {
            \\
            0.058823529411764705  0.2182563120758799  \\
            0.029411764705882353  0.10733717009552084  \\
            0.014705882352941176  0.0330803583604627  \\
            0.007352941176470588  0.008465305217108022  \\
            0.003676470588235294  0.0021198367015300844  \\
            0.001838235294117647  0.0005300253145839631  \\
        }
        ;
    \addplot[color={blu}, mark={o}]
        table[row sep={\\}]
        {
            \\
            0.058823529411764705  0.08404580854370347  \\
            0.029411764705882353  0.017533934390842595  \\
            0.014705882352941176  0.0018963739625958659  \\
            0.007352941176470588  0.00011902571054330884  \\
            0.003676470588235294  7.204988179620357e-6  \\
            0.001838235294117647  4.4567575841005704e-7  \\
        }
        ;
    \addplot[color={blu}, mark={x}]
        table[row sep={\\}]
        {
            \\
            0.058823529411764705  0.12982208606952358  \\
            0.029411764705882353  0.02525321849702212  \\
            0.014705882352941176  0.0017861428907794503  \\
            0.007352941176470588  0.00011121668412472139  \\
            0.003676470588235294  6.911400129007647e-6  \\
            0.001838235294117647  4.293918069271166e-7  \\
        }
        ;
    \addplot[color={grn}, mark={o}]
        table[row sep={\\}]
        {
            \\
            0.058823529411764705  0.10370708759889652  \\
            0.029411764705882353  0.015582315794785017  \\
            0.014705882352941176  0.0007581829484968254  \\
            0.007352941176470588  6.011938335886943e-6  \\
            0.003676470588235294  7.617122625276657e-8  \\
            0.001838235294117647  1.7119541751452498e-9  \\
        }
        ;
    \addplot[color={grn}, mark={x}]
        table[row sep={\\}]
        {
            \\
            0.058823529411764705  0.10804846051077131  \\
            0.029411764705882353  0.01624794371839921  \\
            0.014705882352941176  0.0005706199047665837  \\
            0.007352941176470588  4.413309785133423e-6  \\
            0.003676470588235294  5.919520457725087e-8  \\
            0.001838235294117647  8.012417120819705e-10  \\
        }
        ;
\end{loglogaxis}
\end{tikzpicture}
     \caption{$\beta = 128$}
  \end{subfigure}
  \hfill~
  \caption{L$^2$-convergence comparison of the characteristic ($+$ markers) and
    non-characteristic ($\times$ markers) treatment of the nonlinear interface
    condition for various values of the strength parameter $\beta$.  The red,
    blue, and green curves correspond to SBP interior orders $2$, $4$, and $6$,
  respectively.\label{fig:convergence:interface:1D}}
\end{figure}
Figure~\ref{fig:convergence:interface:1D} shows the convergence results for the
solution at time $t = 1$ for the non-characteristic and characteristic
methods.
The spatial resolutions used in the test are $N = 17 \times 2^r$ with $r = 0, 1,
2, 3, 4, 5$.
In the figure SBP interior orders $2p = 2$, $4$, and $6$ are used and the
nonlinear strength parameters as $\beta = 32$, $64$, and $128$.
As can be seen both methods converge similarly at the expected rates.

\section{Multi-Dimensional Model Problem}\label{sec:model:problem}

\subsection{Continuous Problem}
Let $\Omega \subset \Re^{d}$ be a bounded domain with boundary $\partial
\Omega$.
The boundary is split into two distinct parts: a Dirichlet boundary $\partial
\Omega_{D}$ and a characteristic boundary $\partial \Omega_{C}$.
Additionally, let $\Gamma_{I} \subset \Re^{d-1}$ be a set of interfaces in the
domain.
Unless otherwise noted summation over repeated subscripts is implied, e.g.,
$u_{i} v_{i} = \sum_{i=1}^{d} u_{i} v_{i}$, $u_{ii} = \sum_{i=1}^{d} u_{ii}$,
and $u_{i}C_{ij}v_{j} = \sum_{i=1}^{d}\sum_{j=1}^{d} u_{i}C_{ij}v_{j}$.

As a model problem we consider the scalar second-order, anisotropic wave equation
for displacement $u$:
\begin{subequations}\label{eqn:wave}
  \begin{alignat}{2}
    \label{eqn:wave:pde}
    &\rho \ddot{u} = \partial_{i} C_{ij} \partial_{j} u, &&~
    ~\vec{x} \in \Omega ,~t \in [0, T],\\
    \label{eqn:wave:Dirichlet}
    &u = g_D, &&~ ~\vec{x} \in \partial\Omega_D,~t \in [0, T],\\
    \label{eqn:wave:char:bc:orig}
    &Z \dot{u} + \tau
    = R (Z \dot{u} - \tau) + g_C, &&~ ~\vec{x} \in \partial\Omega_C,~t \in [0,
    T],\\
    \label{eqn:wave:interface:orig}
    &\begin{cases}
       \tau^- = -\tau^+,\\
      \tau^\pm = F(V^{\pm}),
       \end{cases}
    &&~ ~\vec{x} \in \Gamma_{I},~t \in [0, T].
  \end{alignat}
\end{subequations}
Here, the density $\rho > 0$ and the components of the stiffness matrix $C_{ij}$
are taken to be spatially varying. The stiffness matrix is assumed
to be symmetric positive definite: $C_{ij} = C_{ji}$ and $v_{i} C_{ij} v_{j} \ge
0$ with equality only when $v_{i} = 0$ for all $i$.
At interfaces and boundaries the traction $\tau$ is
\begin{equation}
  \label{eqn:def:traction}
  \tau = n_{i} C_{ij} \partial_{j}u,
\end{equation}
where the vector $n_{i}$ is the outward unit normal to boundary.
On $\partial\Omega_C$ the reflection coefficient satisfies $-1 \le R \le 1$
where the shear impedance $Z > 0$ is defined by $Z^2 = \rho n_{i}C_{ij}n_{j}$.
On the interface $\Gamma_{I}$, relationship~\eqref{eqn:wave:interface:orig}
specifies force balance and a nonlinear conditions, respectively.
The normal vector is defined so that $n_{i}^{-}$ points away from the minus side
and $n_{i}^{+}$ points away from the plus side with $n_{i}^{+} = -n_{i}^{-}$.
The superscripts on the material parameters denote which side of the interface
the material parameters are evaluated on.
We define the jump in $\dot{u}$ across the interface by
\begin{equation}
  V^{\pm} = \dot{u}^{\mp} - \dot{u}^{\pm}.
\end{equation}
Physically, The nonlinear function $F(V)$ is the frictional strength of the
interface and is assumed to satisfy $VF(V) \ge 0$.
Force balance and $V^{+} = -V^{-}$ imply that $F(V^{+}) = -F(V^{-})$.

\subsection{Domain Decomposition}\label{sec:domain:decomposition}
Let $\BB$ be a partitioning of $\Omega \subset \Re^{d}$ into $N_{b}$
non-overlapping, curvilinear blocks (quadrilaterals when $d=2$ and hexahedrons
when $d=3$).
For each $B \in \BB$ there is a diffeomorphic mapping $\vec{x}^{B}$ between $B$
and the reference block $\hat{B} = [0, 1]^d$ such that $\vec{x}^{B}(\vec{\xi})
\in B$ for all $\vec{\xi} \in \hat{B}$; the reference block is the same as the
reference domain discussed in \sref{sec:SBP:MD} and the same face numbering is
used.
We use the notation $\hat{\partial}_{i}$ to denote the partial derivative with
respect to $\xi_{i}$.
The Jacobian determinant is denoted as $J^{B}$. For example with $d=2$
\begin{equation}
  J^{B} =
  \left(\hat{\partial}_{1} x_{1}^{B}\right)
  \left(\hat{\partial}_{2} x_{2}^{B}\right)
  -
  \left(\hat{\partial}_{1} x_{2}^{B}\right)
  \left(\hat{\partial}_{2} x_{1}^{B}\right).
\end{equation}
Note that typically the metric terms are computed by first computing
$\hat{\partial}_{l} x_{i}^{B}$ and then metric identities are employed to
calculate $\partial_{j} \xi_{m}^{B}$; see, for example, \citet{Kopriva2006jsc}.

Each block $B \in \BB $ has $2d$ faces, and we let $\partial B^{f}$ for $f = 1,
2, \dots, 2d$ be the faces in physical space and $\partial \hat{B}^{f}$ be the
faces in the reference space.
We assume that each face $B^{f}$ corresponds to either a Dirichlet boundary,
characteristic boundary, nonlinear interface, or a purely computational
interface (i.e., an artificial interface introduced in the partitioning of
$\Omega$).
We let $n_{i}^{B^f}$ denote the outward pointing normal to face $f$ of block $B$
in physical space and $\hat{n}_{i}^{B^f} \equiv \hat{n}_{i}^{f}$ denote the same
outward pointing normal in the reference space; see also~\eqref{eqn:hat:n:kron}.
The relationship between $n_{i}$ and $\hat{n}_{i}$ is
\begin{equation}
  S_{J}^{B^{f}} n^{B^{f}}_{i} = J^{B} \left(\partial_{i} \xi_{k}^{B}\right)
  \hat{n}_{k}^{f},
\end{equation}
where the surface Jacobian $S_{J}^{B^{f}}$ is the normalization factor so that
$n_{i}^{B^{f}}$ is a unit vector.
Given the face numbering convention and properties of the reference unit norm
$\hat{n}_{i}^{f}$, the surface Jacobian with $d = 2$ is
\begin{equation}
  {\left(S_{J}^{B^{f}}\right)}^{2} = {\left(J^{B}\right)}^{2}
  \left(\partial_{i} \xi_{\ceil*{\frac{f}{2}}}^{B}\right)
  \left(\partial_{i} \xi_{\ceil*{\frac{f}{2}}}^{B}\right)~~\text{(no
  summation over $f$)}.
\end{equation}

Before writing down the transformed governing equations, it is useful to define
a few metric term scaled quantities.
For each $B \in \BB$ we define the transformed density and stiffness matrix as
\begin{subequations}
  \begin{align}
    \hat{\rho} &= J \rho,\\
    \hat{C}_{ij} &=
    J
    \left(\partial_{l} \xi_{i}\right)
    C_{lm}
    \left(\partial_{m} \xi_{j}\right);
  \end{align}
\end{subequations}
in this equation, and those that follow, unless needed the superscript $B$
denoting the block number is suppressed.
Similarly, on face $\partial B^{f}$ the shear impedance and traction are defined
as
\begin{subequations}
  \begin{align}
    {\left(\hat{Z}^{f}\right)}^{2}
    &=
    \hat{\rho} \hat{n}^{f}_{i} \hat{C}_{ij} \hat{n}^{f}_{j}
    =
    {\left(S_{J}^{f} Z^{f}\right)}^{2},\\
    \label{eqn:transformed:traction}
    \hat{\tau}^{f}
    &=
    \hat{n}^{f}_{i} \hat{C}_{ij} \hat{\partial}_{j} u
    =
    S_{J}^{f} \tau^{f};
  \end{align}
\end{subequations}
unless needed for clarity, the superscript $B^{f}$ is reduced to $f$.
Finally, the scaled boundary data and interface function are
\begin{subequations}
  \begin{align}
    \hat{g}_{C}^{f} &= S_{J}^{f} g_{C},\\
    \hat{F}(V) &= S_{J}^{f} F(V).
  \end{align}
\end{subequations}

With these, for each $B \in \BB$ governing equations~\eqref{eqn:wave} become
\begin{subequations}\label{eqn:wave:transformed}
  \begin{equation}
    \label{eqn:wave:transformed:pde}
    \hat{\rho} \ddot{u} = \hat{\partial}_{i} \hat{C}_{ij} \hat{\partial}_{j} u, ~
    ~\vec{\xi} \in [0, 1]^{d} ,~t \in [0, T].
  \end{equation}
  For each face $\partial B^{f}$ the boundary and interface condition are
  \begin{alignat}{2}
    \label{eqn:wave:transformed:Dirichlet}
    &u = g_D^{f},
    &&~ \text{if } \partial B^{f} \cap \partial \Omega_{D} \ne \emptyset,\\
    \label{eqn:wave:transformed:char:bc:orig}
    &\hat{Z}^{f} \dot{u}^{f} + \hat{\tau}^{f}
    = R^{f} \left(\hat{Z}^{f} \dot{u}^{f} - \hat{\tau}^{f}\right) +
    \hat{g}_C^{f},
    &&~ \text{if } \partial B^{f} \cap \partial \Omega_{C} \ne \emptyset,\\
    \label{eqn:wave:transformed:interface:orig}
    &\begin{cases}
    	\hat{\tau}^{f^{-}} = -\hat{\tau}^{f^{+}}\\
      \hat{\tau}^{f^{\pm}} = \hat{F}\left(V^{f^{\pm}}\right),\\
    \end{cases}
    &&~ \text{if } \partial B^{f} \cap \Gamma_{I} \ne \emptyset,\\
    \label{eqn:wave:transformed:computational:orig}
    &\begin{cases}
      \hat{\tau}^{f^{-}} =
      -
      \hat{\tau}^{f^{+}},\\
      \dot{u}^{f^{-}} = \dot{u}^{f^{+}},
    \end{cases}
    &&~ \text{otherwise}.
  \end{alignat}
\end{subequations}
where $V^{f^{\pm}} = \dot{u}^{f^{\mp}} - \dot{u}^{f^{\pm}}$.
Here the notation $f^{\pm}$ denotes the two sides of the interface with $f^{-}$
denoting the interior value and $f^{+}$ denoting the exterior (neighboring
block) value.
Namely, let face $\partial B^{f}$ of block $B \in \BB$ be connected to block $C
\in \BB$ along face $\partial C^{f'}$, then $\partial B^{f^{-}} = \partial
B^{f}$ and $\partial B^{f^{+}} = \partial C^{f'}$.
By definition $S_{J}^{f^{+}} = S_{J}^{f^{-}}$ and $\hat{n}^{f^{+}}_{i} = -
\hat{n}^{f^{-}}_{i}$.
Interface conditions~\eqref{eqn:wave:transformed:computational:orig} are not
present in the original governing equations~\eqref{eqn:wave}, and are added to
account for continuity of the solution across computational block interfaces.

\subsection{Characteristic Variables}

As in the one-dimensional case, it is useful to introduce the characteristic
variables
\begin{subequations}
  \begin{align}
    \hat{q}^{f} &= \hat{Z}^{f} \dot{u}^{f} + \hat{\tau}^{f},\\
    \hat{w}^{f} &= \hat{Z}^{f} \dot{u}^{f} - \hat{\tau}^{f},
  \end{align}
\end{subequations}
where $\hat{q}^{f}$ and $\hat{w}^{f}$ propagate out of and into the block
face $f$, respectively; recall that $\hat{\tau}^{f}$ includes the outward
normal.
As before, the displacement and traction can be recovered from the
characteristic variables:
\begin{subequations}
  \label{eqn:wave:transformed:char:phys}
  \begin{align}
    \dot{u}^{f} &= \frac{\hat{q}^{f} + \hat{w}^{f}}{2\hat{Z}^{f}},\\
    \hat{\tau}^{f} &= \frac{\hat{q}^{f} - \hat{w}^{f}}{2}.
  \end{align}
\end{subequations}
With this, the characteristic boundary
condition~\eqref{eqn:wave:transformed:char:bc:orig} can be written as
\begin{equation}
  \label{eqn:wave:transformed:char:bc}
  \hat{q}^{f} = R^{f} \hat{w}^{f} + \hat{g}_C^{f}.
\end{equation}

Interface conditions~\eqref{eqn:wave:transformed:interface:orig}
and~\eqref{eqn:wave:transformed:computational:orig} can both be rewritten in
terms
of the characteristic variables:
\begin{equation}
  \label{eqn:wave:transformed:interface}
  \hat{q}^{f^{\pm}} = \hat{\mathcal{Q}}^{f^{\pm}}\left(\hat{w}^{f^{\pm}}, \hat{w}^{f^{\mp}}\right),
\end{equation}
where the superscript $f^{\pm}$ denote the variable on either side of the
interface.
For computational interface~\eqref{eqn:wave:transformed:computational:orig},
$\hat{\mathcal{Q}}^{f^{\pm}}$ is a linear function:
\begin{equation}
  \label{eqn:computational:Q}
  \hat{\mathcal{Q}}^{f^{\pm}}\left(\hat{w}^{f^{\pm}}, \hat{w}^{f^{\mp}}\right)
  =
  \frac{
    2\hat{Z}^{f^{\pm}} \hat{w}^{f^{\mp}}
    +
    (\hat{Z}^{f^{\pm}} - \hat{Z}^{f^{\mp}}) \hat{w}^{f^{\pm}}
  }{\hat{Z}^{f^{+}} + \hat{Z}^{f^{-}}};
\end{equation}
when $\hat{Z}^{f^{+}} = \hat{Z}^{f^{-}}$ this reduces to transmission of the
characteristic variable across the interface:
$\hat{\mathcal{Q}}^{f^{\pm}}\left(\hat{w}^{f^{\pm}}, \hat{w}^{f^{\mp}}\right) =
\hat{w}^{f^{\mp}}$.

As discussed in one-dimension, for nonlinear interface
condition~\eqref{eqn:wave:transformed:interface:orig}, in general there is no
closed form expression for $\hat{\mathcal{Q}}^{\pm}$.
As shown in \aref{app:rootfinding}, for a given $\dot{u}^{f^{\pm}}$ and
$\hat{\tau}^{f^{\pm}}$ the function $\hat{\mathcal{Q}}^{\pm}$ can be found to be
consistent with the interface condition by solving the nonlinear system
\begin{subequations}
  \label{eqn:Q}
  \begin{align}
    \label{eqn:dotu:Q}
    \dot{u}^{f^{\pm}} &=
    \frac{\hat{\mathcal{Q}}^{f^{\pm}} + \hat{w}^{f^{\pm}}}{2\hat{Z}^{f^{\pm}}},\\
    \label{eqn:tau:Q}
    \hat{\tau}^{f^{\pm}} &=
    \frac{\hat{\mathcal{Q}}^{f^{\pm}} - \hat{w}^{f^{\pm}}}{2} =
    \hat{F}\left(V^{f^{\pm}}\right).
  \end{align}
\end{subequations}

\subsection{Energy Analysis}
To guide the development of the numerical scheme, we now develop an energy
estimate for governing equation~\eqref{eqn:wave:transformed:pde}.
We define a seminorm $E(u)$ and then show that $\dot{E}(u(\cdot,
t)) \le 0$ when $g^{f}_{D} = \hat{g}_{C}^{f} = 0$ for all $t > 0$; with non-zero
boundary data energy growth due to the boundary conditions must be allowed.

For the transformed system~\eqref{eqn:wave:transformed}, the energy in block $B
\in \BB$ is
\begin{equation}
  \label{eqn:wave:transformed:energy:block}
  E^{B} = \frac{1}{2} \int_{\hat{B}} \left(\hat{\rho} \dot{u}^2 +
  \left(\hat{\partial}_{i} u\right)
  \hat{C}_{ij}
  \left(\hat{\partial}_{j} u\right)
  \right);
\end{equation}
this is valid seminorm of $u$, namely $E^{B} \ge 0$ for all $u$, because the
stiffness matrix is symmetric positive definite.
The total energy in the domain is then
\begin{equation}
  \label{eqn:wave:transformed:energy}
  E = \sum_{B \in \BB} E^{B}.
\end{equation}
\begin{lemma}\label{lemma:wave:energy:estimate}
  Governing equations~\eqref{eqn:wave:transformed} with
  energy~\eqref{eqn:wave:transformed:energy} satisfy $\dot{E} \le 0$ if
  $g^{f}_{D} = \hat{g}^{f}_{C} = 0$.
\end{lemma}
\begin{proof}
  Taking the time derivative of block
  energy~\eqref{eqn:wave:transformed:energy:block}, substituting in governing
  equations~\eqref{eqn:wave:transformed:pde}, and applying the divergence
  theorem gives
  \begin{equation}
    \label{eqn:con:faceenergy}
    \dot{E}^{B} =
    \sum_{f=1}^{2d}\int_{\partial \hat{B}^{f}} \dot{u}^{f}\hat{\tau}^{f}.
  \end{equation}
  If face $f$ is a Dirichlet boundary then applying boundary
  condition~\eqref{eqn:wave:transformed:Dirichlet} with $g^{f}_{D} = 0$ gives
  \begin{equation}
    \int_{\partial \hat{B}^{f}} \dot{u}^{f}\hat{\tau}^{f}
    =
    \int_{\partial \hat{B}^{f}} g^{f}_{D} \hat{\tau}^{f}
    = 0.
  \end{equation}
  Similarly, if face $f$ is a characteristic boundary, applying physical to
  characteristic variable transformation~\eqref{eqn:wave:transformed:char:phys}
  and using characteristic boundary
  condition~\eqref{eqn:wave:transformed:char:bc} gives
  \begin{equation}
    \begin{split}
      \int_{\partial \hat{B}^{f}} \dot{u}^{f}\hat{\tau}^{f}
    &=
    -\int_{\partial \hat{B}^{f}}
    \frac{1}{4\hat{Z}}
    \left(1+R^{f}\right)
    \left(1-R^{f}\right)
    \hat{w}^2\\
    &=
    -\int_{\partial \hat{B}^{f}}
    \frac{1}{4\hat{Z}}
    \left(1-{\left(R^{f}\right)}^2\right)
    {\left(\hat{w}^{f}\right)}^2
    \le 0,
    \end{split}
  \end{equation}
  since $-1 \le R^f \le 1$.
  If face $f$ is an interface, adding both contributions gives
  \begin{equation}
    \begin{split}
      \int_{\partial \hat{B}^{f}} \dot{u}^{f^{-}}\hat{\tau}^{f^{-}}
      +
      \int_{\partial \hat{B}^{f}} \dot{u}^{f^{+}}\hat{\tau}^{f^{+}}
      &=
      -\int_{\partial \hat{B}^{f}}
      \left(\dot{u}^{f^{+}} - \dot{u}^{f^{-}}\right)\hat{\tau}^{f^{-}}
      =
      -\int_{\partial \hat{B}^{f}}
      V^{f^{-}}\hat{\tau}^{f^{-}},
      \end{split}
  \end{equation}
  where we have used that $\hat{\tau}^{f^{+}} = - \hat{\tau}^{f^{-}}$
  by~\eqref{eqn:wave:transformed:computational:orig}
  and~\eqref{eqn:wave:transformed:interface:orig}.
  For the case of computational interface
  condition~\eqref{eqn:wave:transformed:computational:orig} the particle
  velocity is continuous, thus $V^{f^{-}} = 0$ and the interface leads to a zero
  rate of change in energy.
  When the interface is governed by nonlinear interface
  condition~\eqref{eqn:wave:transformed:interface:orig} we have that
  \begin{equation}
    \int_{\partial \hat{B}^{f}} \dot{u}^{f^{-}}\hat{\tau}^{f^{-}}
    +
    \int_{\partial \hat{B}^{f}} \dot{u}^{f^{+}}\hat{\tau}^{f^{+}}
    =
    -
    \int_{\partial \hat{B}^{f}}
    V^{f^{-}}\hat{F}\left(V^{f^{-}}\right) \le 0,
  \end{equation}
  by the condition that $V\hat{F}(V) \ge 0$.
  \qed
\end{proof}

\section{Multi-Block Semidiscretization}\label{sec:discretization}
A single block semidiscretization of~\eqref{eqn:wave:transformed} with weak
enforcement of boundary conditions is
\begin{equation}
  \begin{split}
  \label{eqn:disc:wave:orig}
    \MM{\rho} \ddot{\VV{u}}
    =&\;
    \MM{D}_{ij}^{(\hat{C}_{ij})} \VV{u}
    +
    \sum_{f=1}^{2d}
    \MM{H}^{-1}
    {\left(\Mm{L}^{f}\right)}^{T}
    \mm{H}^{f}
    \left(
    \vv{\hat{\tau}}^{*f}
    - n^{f}_{i}\mm{\hat{C}}_{ij}^{f}\Mm{B}^{f}_{j}\VV{u}
    \right)
    \\&
    -
    \sum_{f=1}^{2d}
    \MM{H}^{-1}
    {\left(\Mm{B}^{f}_{j}\right)}^{T}
    n^{f}_{i}\mm{\hat{C}}_{ij}^{f}
    \mm{H}^{f}
    \left(
    \vv{u}^{*f}
    - \Mm{L}^{f}\VV{u}
    \right).
  \end{split}
\end{equation}
which after multiplying by $\VV{H}$ and applying the multidimensional SBP
property~\eqref{eqn:multi:sbp} gives a form that is more convenient for
analysis:
\begin{equation}
  \label{eqn:disc:wave}
  \begin{split}
    \MM{\rho} \MM{H} \ddot{\VV{u}}
    =&\;
    -\MM{A}_{ij}^{(\hat{C}_{ij})} \VV{u}
    +
    \sum_{f=1}^{2d}
    {\left(\Mm{L}^{f}\right)}^{T}
    \mm{H}^{f}
    \vv{\hat{\tau}}^{*f}
    \\&
    -
    \sum_{f=1}^{2d}
    {\left(\Mm{B}^{f}_{j}\right)}^{T}
    n^{f}_{i}\mm{\hat{C}}_{ij}^{f}
    \mm{H}^{f}
    \left(
    \vv{u}^{*f}
    - \vv{u}^{f}
    \right).
  \end{split}
\end{equation}
Here we have defined
\begin{equation}
    \vv{u}^{f} = \Mm{L}^{f}\VV{u},
\end{equation}
and $\MM{\rho}$ is a diagonal matrix of density $\rho$ evaluated at the grid
points.  The numerical flux vectors $\vv{\hat{\tau}}^{*f}$ and $\vv{u}^{*f}$,
which depend on the specific boundary or interface condition, are discussed in
detail below.

We define the energy in the domain as
\begin{equation}
  \label{eqn:disc:energy}
  \mathcal{E} = \sum_{B \in \BB} \mathcal{E}^{B},
\end{equation}
where the energy in block $B$ is
\begin{equation}
  \label{eqn:disc:block:energy}
  \begin{split}
    \mathcal{E}^{B} =&\;
    \frac{1}{2} \dot{\VV{u}}^{T} \MM{H}\MM{\rho}\dot{\VV{u}}
    +
    \frac{1}{2} \VV{u}^{T} \MM{A}_{ij}^{(\hat{C}_{ij})}\VV{u}\\
    +&
    \frac{1}{2}\sum_{f=1}^{2d}
    \left(
    {\left(\vv{\hat{\tau}}^{f}\right)}^{T}
    {\left(\mm{X}^{f}\right)}^{-1}
    \mm{H}^{f}
    {\left(\vv{\hat{\tau}}^{f}\right)}
    -
    {\left(\vv{\hat{T}}^{f}\right)}^{T}
    {\left(\mm{X}^{f}\right)}^{-1}
    \mm{H}^{f}
    {\left(\vv{\hat{T}}^{f}\right)}
    \right).
  \end{split}
\end{equation}
Here we have defined the diagonal matrix
\begin{equation}
  \label{eqn:Xf}
  \mm{X}^{f}
  =
  n^{f}_{i}n^{f}_{j}\mm{\hat{C}}_{ij}^{f}
  \mm{\Gamma}^{f},
\end{equation}
where $\mm{\Gamma}^{f}$ is a diagonal penalty parameter matrix the entries of
which must be sufficient large; a lower bound for $\mm{\Gamma}^{f}$ is given
by~\eqref{eqn:penalty:param}.
Additionally we define the block face tractions
\begin{subequations}
  \label{eqn:T:tau:hat:f}
  \begin{align}
    \vv{\hat{T}}^{f} &= n^{f}_{i}\mm{\hat{C}}_{ij}^{f}\Mm{B}^{f}_{j}\VV{u},\\
    \label{eqn:tau:hat:f}
    \vv{\hat{\tau}}^{f} &= \vv{\hat{T}}^{f}
    +
    \mm{X}^{f}
    \left(
    \vv{u}^{*f} - \Mm{L}^{f}\VV{u}
    \right).
  \end{align}
\end{subequations}
Essentially, discrete energy~\eqref{eqn:disc:block:energy} is a
direct discretization of continuous
energy~\eqref{eqn:wave:transformed:energy:block} with an additional penalty on
the faces for the mismatch between two alternative approximations of the
traction $\hat{\tau}^{f}$~\eqref{eqn:transformed:traction}.

The discrete energy satisfies the following theorem; proof in
\aref{app:disc:energy:seminorm}.
\begin{theorem}\label{thm:disc:energy:seminorm}
  Energy~\eqref{eqn:disc:block:energy} is a seminorm of the solution if
  $\mm{\Gamma}^{f}$ is positive and sufficiently large.
\end{theorem}

\begin{corollary}
  For each block $B \in \BB$, the energy rate of change is
  \begin{equation}
    \dot{\mathcal{E}}^{B} = \sum_{f=1}^{2d} \dot{\mathcal{E}}^{f},
  \end{equation}
  with the energy face rate of change being
  \begin{equation}
    \label{eqn:disc:faceenergy}
    \dot{\mathcal{E}}^{f}
    =
    {\left(\vv{\hat{\tau}}^{*f}\right)}^{T}
    \mm{H}^{f}
    \dot{\vv{u}}^{f}
    +
    {\left(\vv{\hat{\tau}}^{f}\right)}^{T}
    \mm{H}^{f}
    {\left(
        \dot{\vv{u}}^{*f}
        - \dot{\vv{u}}^{f}
    \right)}.
  \end{equation}
\end{corollary}
\begin{proof}
  For a single block, taking the time derivative of block
  energy~\eqref{eqn:disc:block:energy} and using
  discretization~\eqref{eqn:disc:wave} gives
  \begin{equation}
    \begin{split}
      \dot{\mathcal{E}}^{f} =&\;
      {\left(\vv{\hat{\tau}}^{f}\right)}^{T}
      {\left(\mm{X}^{f}\right)}^{-1}
      \mm{H}^{f}
      {\left(\dot{\vv{\hat{\tau}}}^{f}\right)}
      -
      {\left(\vv{\hat{T}}\right)}^{T}
      {\left(\mm{X}^{f}\right)}^{-1}
      \mm{H}^{f}
      {\left(\dot{\vv{\hat{T}}}\right)}
      \\
                             &+
                             {\left(\vv{\hat{\tau}}^{*f}\right)}^{T}
                             \mm{H}^{f}
                             \dot{\vv{u}}^{f}
                             -
                             {\left(
                                 \vv{u}^{*f}
                                 - \vv{u}^{f}
                             \right)}^{T}
                             \mm{H}^{f}
                             \dot{\vv{\hat{T}}}.
    \end{split}
  \end{equation}
  Using the definition of $\vv{\hat{\tau}}^{f}$~\eqref{eqn:tau:hat:f} the
  rate of change of face energy simplifies to face
  rate~\eqref{eqn:disc:faceenergy}.
  \qed
\end{proof}

Discrete face energy rate~\eqref{eqn:disc:faceenergy} is of the same form as
the continuous counterpart~\eqref{eqn:con:faceenergy}, namely a boundary
integral of the particle velocity times the traction at the boundary.
Stability is now reduced to showing that if a face $f$ is on a physical
boundary that $\dot{\mathcal{E}}^{f} \le 0$ and if on a block interface that
$\dot{\mathcal{E}}^{f^{-}} + \dot{\mathcal{E}}^{f^{+}} \le 0$.

In the remainder of this section numerical fluxes are given for characteristic
boundary conditions as well as the characteristic treatment of computational and
nonlinear interfaces.
In \aref{app:standard:approach} the typical SBP-SAT numerical fluxes for
Dirichlet, Neumann, and Characteristic boundary conditions as well as
computational interfaces are given, e.g., those from \citet{VirtaMattsson2014}
with the improved Dirichlet penalty parameter of \citet{AlmquistDunham2020}.
The standard approaches specify $u^{*}$ directly, i.e., they do not require an
additional block face variable be integrated in time.

\subsection{Characteristic Boundary Conditions}\label{subsec:characteristic:bcs}

When block face $f$ corresponds to a characteristic
boundary~\eqref{eqn:wave:transformed:char:bc} we choose values of
$\vv{\hat{\tau}}^{*f}$ and $\dot{\vv{u}}^{*f}$ which preserves the outgoing
characteristic variable while also satisfying the boundary condition:
\begin{subequations}
  \begin{align}
    \mm{\hat{Z}}^{f} \dot{\vv{u}}^{*f} - \vv{\hat{\tau}}^{*f}
    &=
    \mm{\hat{Z}}^{f} \dot{\vv{u}}^{f} - \vv{\hat{\tau}}^{f},\\
    \nonumber
    \mm{\hat{Z}}^{f} \dot{\vv{u}}^{*f} + \vv{\hat{\tau}}^{*f}
    &=
    \mm{R}^{f} \left(\mm{\hat{Z}}^{f} \dot{\vv{u}}^{*f} -
    \vv{\hat{\tau}}^{*f}\right) + \vv{\hat{g}}_{C}^{f};
  \end{align}
\end{subequations}
Solving these equations for the numerical fluxes gives
\begin{subequations}
  \label{eqn:disc:bcs:Characteristic}
  \begin{align}
    \dot{\vv{u}}^{*f} &=
    \frac{\mm{I}+\mm{R}^{f}}{2}
    \left(
    \dot{\vv{u}}^{f}
    -
    {\left(\mm{\hat{Z}}^{f}\right)}^{-1}
    \vv{\hat{\tau}}^{f}
    \right)
    +
    \frac{1}{2}{\left(\mm{\hat{Z}}^{f}\right)}^{-1} \vv{\hat{g}}_{C}^{f},
    \\
    \vv{\hat{\tau}}^{*f} &=
    -\frac{\mm{I}-\mm{R}^{f}}{2}
    \left(
    \mm{\hat{Z}}^{f}
    \dot{\vv{u}}^{f}
    -
    \vv{\hat{\tau}}^{f}
    \right)
    +
    \frac{1}{2} \vv{\hat{g}}_{C}^{f}.
  \end{align}
\end{subequations}
As in the one-dimensional formulation, numerical flux $\vv{u}^{*f}$ must be
stored along the face and integrated in time.

Using the characteristic boundary treatment~\eqref{eqn:disc:bcs:Characteristic}
with $g_{c} = 0$ in face energy rate of change~\eqref{eqn:disc:faceenergy} gives
\begin{equation}
  \begin{split}
    \dot{\mathcal{E}}^{f} =&\;
    -
    {\left(
    \dot{\vv{u}}^{f}
    \right)}^{T}
    \frac{\mm{I}-\mm{R}^{f}}{2}
    \mm{\hat{Z}}^{f}
    \mm{H}^{f}
    \dot{\vv{u}}^{f}
    -
    {\left(\vv{\hat{\tau}}^{f}\right)}^{T}
    \frac{\mm{I}+\mm{R}^{f}}{2}
    {\left(\mm{\hat{Z}}^{f}\right)}^{-1}
    \mm{H}^{f}
    \vv{\hat{\tau}}^{f}.
  \end{split}
\end{equation}
Since the reflection coefficient satisfies $-1 \le R \le 1$, the boundary
treatment is energy stable: $\dot{\mathcal{E}}_{f} \le 0$.

\subsection{Characteristic Interface}\label{subsec:characteristic:interface}
For characteristic interfaces, computational or nonlinear, the aim is to define
the numerical fluxes to satisfy the interface condition in a way that preserves
the characteristic variables propagating into the interface.
As noted in \sref{sec:domain:decomposition}, the nonlinear and computational
interface conditions can be enforced using the functions
$\hat{\mathcal{Q}}^{f^{\pm}}$~\eqref{eqn:wave:transformed:interface}.
Thus we define $\vv{\hat{\tau}}^{*f^\pm}$ and $\dot{\vv{u}}^{*f^\pm}$ so that
they satisfy
\begin{subequations}
  \begin{align}
    \vv{\hat{w}}^{*f^{\pm}} &=
    \mm{\hat{Z}}^{f^\pm}\dot{\vv{u}}^{f^{\pm}} - \vv{\hat{\tau}}^{f^{\pm}}
    =
    \mm{\hat{Z}}^{f^\pm}\dot{\vv{u}}^{*f^{\pm}} - \vv{\hat{\tau}}^{*f^{\pm}},\\
    \vv{\hat{q}}^{*f^{\pm}}
    &=
    \hat{\mathcal{Q}}^{f^{\pm}}\left( \vv{\hat{w}}^{*f^{\pm}}, \vv{\hat{w}}^{*f^{\mp}}\right)
    =
    \mm{\hat{Z}}^{f^\pm}\dot{\vv{u}}^{*f^{\pm}} + \vv{\hat{\tau}}^{*f^{\pm}}.
  \end{align}
\end{subequations}
Solving for the numerical fluxes then gives
\begin{subequations}
  \label{eqn:disc:interface}
  \begin{align}
    \dot{\vv{u}}^{*f^{\pm}} &=
    \frac{1}{2}
    {\left(\mm{\hat{Z}}^{f^{\pm}}\right)}^{-1}
    \left(\vv{\hat{q}}^{*f^{\pm}} + \vv{\hat{w}}^{*f^{\pm}}\right),\\
    \vv{\hat{\tau}}^{*f^{\pm}} &=
    \frac{1}{2}
    \left(\vv{\hat{q}}^{*f^{\pm}} - \vv{\hat{w}}^{*f^{\pm}}\right).
  \end{align}
\end{subequations}
Since $\vv{\hat{\tau}}^{*f^{\pm}}$ and $\dot{\vv{u}}^{*f^{\pm}}$ satisfy the
interface conditions, it follows that for a computational interface:
\begin{subequations}
  \begin{align}
    \vv{\hat{\tau}}^{*f^{-}} = - \vv{\hat{\tau}}^{*f^{+}},\\
    \dot{\vv{u}}^{*f^{-}} = \dot{\vv{u}}^{*f^{+}},
  \end{align}
\end{subequations}
and for the nonlinear interface:
\begin{subequations}
  \begin{align}
    \vv{\hat{\tau}}^{*f^{-}} = - \vv{\hat{\tau}}^{*f^{+}},\\
    \vv{\hat{\tau}}^{*f^{\pm}} = \hat{F}\left(\vv{V}^{*f^{\pm}}\right),
  \end{align}
\end{subequations}
where $\vv{V}^{*f^{\pm}} = \dot{\vv{u}}^{*f^{\mp}} - \dot{\vv{u}}^{*f^{\pm}}$.
Since it is required that $V \hat{F}(V) \ge 0$, for both the computational and
nonlinear interface treatment
\begin{equation}
  \label{eqn:disc:tau:V}
  {\left(\vv{\hat{\tau}}^{*f^{\pm}}\right)}^{T} \vv{V}^{*f^{\pm}} \ge 0;
\end{equation}
in the computational interface  $\vv{V}^{*f^{\pm}} = \vv{0}$.

In order to analyze the interface treatment, it is useful to define the grid
based characteristic variables
\begin{equation}
  \begin{split}
    \vv{\hat{q}}^{f^{\pm}} = \mm{\hat{Z}}^{f^{\pm}} \dot{\vv{u}}^{f^{\pm}} +
    \vv{\hat{\tau}}^{f^{\pm}},
  \end{split}
\end{equation}
so that we can write
\begin{subequations}
  \begin{align}
    \dot{\vv{u}}^{f^{\pm}} &=
    \frac{1}{2}{\left(\mm{\hat{Z}}^{f^{\pm}}\right)}^{-1} \left(\vv{\hat{q}}^{f^{\pm}} +
    \vv{w}^{*f^{\pm}}\right),\\
    \vv{\hat{\tau}}^{f^{\pm}} &=
    \frac{1}{2}\left(\vv{\hat{q}}^{f^{\pm}} - \vv{w}^{*f^{\pm}}\right);
  \end{align}
\end{subequations}
identical expressions can be written for the numerical fluxes
$\vv{\hat{q}}^{*f^{\pm}}$, $\dot{\vv{u}}^{*f^{\pm}}$, and
$\vv{\hat{\tau}}^{*f^{\pm}}$.
Using these in the face energy rate of change~\eqref{eqn:disc:faceenergy} gives
\begin{equation}
  \label{eqn:char:inter:disp}
  \begin{split}
    \dot{\mathcal{E}}^{f^{\pm}} =&\;
    {\left(\vv{\hat{\tau}}^{*f^{\pm}}\right)}^{T}
    \mm{H}^{f}
    \dot{\vv{u}}^{*f^{\pm}}
    -
    {\left(\vv{\hat{\tau}}^{*f^{\pm}} - \vv{\hat{\tau}}^{f^{\pm}}\right)}^{T}
    \mm{H}^{f}
    \left(
    \dot{\vv{u}}^{*f^{\pm}}
    -
    \dot{\vv{u}}^{f^{\pm}}
    \right)\\
    =&\;
    {\left(\vv{\hat{\tau}}^{*f^{\pm}}\right)}^{T}
    \mm{H}^{f}
    \dot{\vv{u}}^{*f^{\pm}}
    -
    \frac{1}{4}
    {\left(
    \vv{\hat{q}}^{*f^{\pm}} - \vv{\hat{q}}^{f^{\pm}}
    \right)}^{T}
    {\left(\mm{\hat{Z}}^{f^{\pm}}\right)}^{-1}
    \mm{H}^{f}
    \left(
    \vv{\hat{q}}^{*f^{\pm}} - \vv{\hat{q}}^{f^{\pm}}
    \right).
  \end{split}
\end{equation}
Adding the two sides of an interface together yields
\begin{equation}
  \begin{split}
    \dot{\mathcal{E}}^{f^{-}} + \dot{\mathcal{E}}^{f^{+}} =&\;
    -{\left(\vv{\hat{\tau}}^{*f^{-}}\right)}^{T}
    \mm{H}^{f}
    \vv{V}^{*f^{-}}\\
    &
    -
    \frac{1}{4}
    {\left(
    \vv{\hat{q}}^{*f^{-}} - \vv{\hat{q}}^{f^{-}}
    \right)}^{T}
    {\left(\mm{\hat{Z}}^{f^{-}}\right)}^{-1}
    \mm{H}^{f}
    \left(
    \vv{\hat{q}}^{*f^{-}} - \vv{\hat{q}}^{f^{-}}
    \right)\\
    &
    -
    \frac{1}{4}
    {\left(
    \vv{\hat{q}}^{*f^{+}} - \vv{\hat{q}}^{f^{+}}
    \right)}^{T}
    {\left(\mm{\hat{Z}}^{f^{+}}\right)}^{-1}
    \mm{H}^{f}
    \left(
    \vv{\hat{q}}^{*f^{+}} - \vv{\hat{q}}^{f^{+}}
    \right).
  \end{split}
\end{equation}
Here we have used that $\vv{\hat{\tau}}^{*f^{+}} = -\vv{\hat{\tau}}^{*f^{-}}$.
Energy stability results since this face energy rate of change is non-positive
due to the positivity result of~\eqref{eqn:disc:tau:V} and the fact that the
second and third terms are in quadratic form.

\section{Two-Dimensional Numerical Experiments}\label{sec:numerical:experiments}
To test semidiscrete scheme~\eqref{eqn:disc:wave} we introduce a velocity
variable $\VV{v} = \dot{\VV{u}}$ and write the method as a first order
system of equations:
\begin{subequations}
  \begin{align}
    \dot{\VV{u}} =&\; \VV{v},\\
    \dot{\VV{v}}
    =&\;
    \MM{\rho}^{-1}
    \MM{D}_{ij}^{(\hat{C}_{ij})} \VV{u}
    +
    \sum_{f=1}^{2d}
    \MM{\rho}^{-1}
    \MM{H}^{-1}
    {\left(\Mm{L}^{f}\right)}^{T}
    \mm{H}^{f}
    \left(
    \vv{\hat{\tau}}^{*f}
    - n^{f}_{i}\mm{\hat{C}}_{ij}^{f}\Mm{B}^{f}_{j}\VV{u}
    \right)
    \\&\notag
    -
    \sum_{f=1}^{2d}
    \MM{\rho}^{-1}
    \MM{H}^{-1}
    {\left(\Mm{B}^{f}_{j}\right)}^{T}
    n^{f}_{i}\mm{\hat{C}}_{ij}^{f}
    \mm{H}^{f}
    \left(
    \vv{u}^{*f}
    - \Mm{L}^{f}\VV{u}
    \right);
  \end{align}
\end{subequations}
as needed, additional variables are also introduced to track the numerical
fluxes.
The error is measured using the discrete L$^{2}$-norm
\begin{equation}
    \|\Delta\VV{u}\|_{H} =
    \sqrt{\sum_{b=1}^{N_{b}}
    {\left(\Delta\VV{u}^B\right)}^{T} \MM{J}^B \MM{H}^B \Delta\VV{u}^{B}},
\end{equation}
where $\Delta\VV{u}$ is the difference between the numerical and analytic
solution evaluated at the grid points.
In all cases the penalty parameter is chosen to be at the stability limit, i.e.,
the equality condition of~\eqref{eqn:penalty:param}.

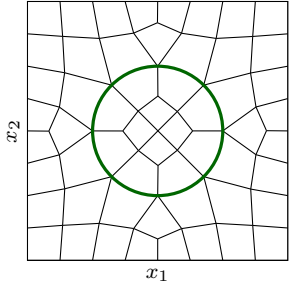
\begin{figure}
  \begin{center}
\begin{tikzpicture}[]
  \begin{axis}[ylabel = {$x_{2}$}, xmin = {-2}, xmax = {2}, ymax = {2}, xlabel =
    {$x_{1}$}, width=5cm, height=5cm, ticks=none, ymin = {-2}]\addplot+ [no marks, solid, black]coordinates {
(1.0, 0.0)
(0.5372134, 0.0002884168)
};
\addplot+ [no marks, solid, black]coordinates {
(0.7071067811865476, 0.7071067811865475)
(0.3017369, 0.3000469)
};
\addplot+ [no marks, solid, black]coordinates {
(0.5372134, 0.0002884168)
(0.3017369, 0.3000469)
};
\addplot+ [no marks, solid, black]coordinates {
(6.123233995736766e-17, 1.0)
(0.00333461, 0.5383348)
};
\addplot+ [no marks, solid, black]coordinates {
(0.3017369, 0.3000469)
(0.00333461, 0.5383348)
};
\addplot+ [no marks, solid, black]coordinates {
(-0.7071067811865475, 0.7071067811865476)
(-0.3003842, 0.3066312)
};
\addplot+ [no marks, solid, black]coordinates {
(0.00333461, 0.5383348)
(-0.3003842, 0.3066312)
};
\addplot+ [no marks, solid, black]coordinates {
(-1.0, 1.2246467991473532e-16)
(-0.5360446, -0.0008435451)
};
\addplot+ [no marks, solid, black]coordinates {
(-0.3003842, 0.3066312)
(-0.5360446, -0.0008435451)
};
\addplot+ [no marks, solid, black]coordinates {
(-0.7071067811865475, -0.7071067811865476)
(-0.2957995, -0.3036193)
};
\addplot+ [no marks, solid, black]coordinates {
(-0.5360446, -0.0008435451)
(-0.2957995, -0.3036193)
};
\addplot+ [no marks, solid, black]coordinates {
(-3.8285686989269494e-16, -1.0)
(0.005527262, -0.5353951)
};
\addplot+ [no marks, solid, black]coordinates {
(-0.2957995, -0.3036193)
(0.005527262, -0.5353951)
};
\addplot+ [no marks, solid, black]coordinates {
(0.7071067811865476, -0.7071067811865475)
(0.3022778, -0.2983624)
};
\addplot+ [no marks, solid, black]coordinates {
(0.005527262, -0.5353951)
(0.3022778, -0.2983624)
};
\addplot+ [no marks, solid, black]coordinates {
(0.3022778, -0.2983624)
(0.5372134, 0.0002884168)
};
\addplot+ [no marks, solid, black]coordinates {
(0.3022778, -0.2983624)
(0.003680641, 0.0008290927)
};
\addplot+ [no marks, solid, black]coordinates {
(-0.2957995, -0.3036193)
(0.003680641, 0.0008290927)
};
\addplot+ [no marks, solid, black]coordinates {
(0.003680641, 0.0008290927)
(0.3017369, 0.3000469)
};
\addplot+ [no marks, solid, black]coordinates {
(0.003680641, 0.0008290927)
(-0.3003842, 0.3066312)
};
\addplot+ [no marks, solid, black]coordinates {
(2.0, 2.0)
(2.0, 1.5)
};
\addplot+ [no marks, solid, black]coordinates {
(1.5, 2.0)
(1.459547, 1.459424)
};
\addplot+ [no marks, solid, black]coordinates {
(2.0, 2.0)
(1.5, 2.0)
};
\addplot+ [no marks, solid, black]coordinates {
(2.0, 1.5)
(1.459547, 1.459424)
};
\addplot+ [no marks, solid, black]coordinates {
(1.0, 2.0)
(0.8982224, 1.429057)
};
\addplot+ [no marks, solid, black]coordinates {
(1.5, 2.0)
(1.0, 2.0)
};
\addplot+ [no marks, solid, black]coordinates {
(1.459547, 1.459424)
(0.8982224, 1.429057)
};
\addplot+ [no marks, solid, black]coordinates {
(0.5, 2.0)
(0.3667017, 1.509209)
};
\addplot+ [no marks, solid, black]coordinates {
(1.0, 2.0)
(0.5, 2.0)
};
\addplot+ [no marks, solid, black]coordinates {
(0.8982224, 1.429057)
(0.3667017, 1.509209)
};
\addplot+ [no marks, solid, black]coordinates {
(0.0, 2.0)
(-0.003366247, 1.671689)
};
\addplot+ [no marks, solid, black]coordinates {
(0.5, 2.0)
(0.0, 2.0)
};
\addplot+ [no marks, solid, black]coordinates {
(0.3667017, 1.509209)
(-0.003366247, 1.671689)
};
\addplot+ [no marks, solid, black]coordinates {
(-0.5, 2.0)
(-0.3693048, 1.511127)
};
\addplot+ [no marks, solid, black]coordinates {
(0.0, 2.0)
(-0.5, 2.0)
};
\addplot+ [no marks, solid, black]coordinates {
(-0.003366247, 1.671689)
(-0.3693048, 1.511127)
};
\addplot+ [no marks, solid, black]coordinates {
(-1.0, 2.0)
(-0.8962192, 1.428906)
};
\addplot+ [no marks, solid, black]coordinates {
(-0.5, 2.0)
(-1.0, 2.0)
};
\addplot+ [no marks, solid, black]coordinates {
(-0.3693048, 1.511127)
(-0.8962192, 1.428906)
};
\addplot+ [no marks, solid, black]coordinates {
(-1.5, 2.0)
(-1.460039, 1.460656)
};
\addplot+ [no marks, solid, black]coordinates {
(-1.0, 2.0)
(-1.5, 2.0)
};
\addplot+ [no marks, solid, black]coordinates {
(-0.8962192, 1.428906)
(-1.460039, 1.460656)
};
\addplot+ [no marks, solid, black]coordinates {
(-2.0, 2.0)
(-2.0, 1.5)
};
\addplot+ [no marks, solid, black]coordinates {
(-1.5, 2.0)
(-2.0, 2.0)
};
\addplot+ [no marks, solid, black]coordinates {
(-1.460039, 1.460656)
(-2.0, 1.5)
};
\addplot+ [no marks, solid, black]coordinates {
(-2.0, -2.0)
(-2.0, -1.5)
};
\addplot+ [no marks, solid, black]coordinates {
(-1.5, -2.0)
(-1.460145, -1.459283)
};
\addplot+ [no marks, solid, black]coordinates {
(-2.0, -2.0)
(-1.5, -2.0)
};
\addplot+ [no marks, solid, black]coordinates {
(-2.0, -1.5)
(-1.460145, -1.459283)
};
\addplot+ [no marks, solid, black]coordinates {
(-1.0, -2.0)
(-0.8988851, -1.428805)
};
\addplot+ [no marks, solid, black]coordinates {
(-1.5, -2.0)
(-1.0, -2.0)
};
\addplot+ [no marks, solid, black]coordinates {
(-1.460145, -1.459283)
(-0.8988851, -1.428805)
};
\addplot+ [no marks, solid, black]coordinates {
(-0.5, -2.0)
(-0.3680867, -1.508369)
};
\addplot+ [no marks, solid, black]coordinates {
(-1.0, -2.0)
(-0.5, -2.0)
};
\addplot+ [no marks, solid, black]coordinates {
(-0.8988851, -1.428805)
(-0.3680867, -1.508369)
};
\addplot+ [no marks, solid, black]coordinates {
(0.0, -2.0)
(0.000342336, -1.669245)
};
\addplot+ [no marks, solid, black]coordinates {
(-0.5, -2.0)
(0.0, -2.0)
};
\addplot+ [no marks, solid, black]coordinates {
(-0.3680867, -1.508369)
(0.000342336, -1.669245)
};
\addplot+ [no marks, solid, black]coordinates {
(0.5, -2.0)
(0.3679869, -1.50956)
};
\addplot+ [no marks, solid, black]coordinates {
(0.0, -2.0)
(0.5, -2.0)
};
\addplot+ [no marks, solid, black]coordinates {
(0.000342336, -1.669245)
(0.3679869, -1.50956)
};
\addplot+ [no marks, solid, black]coordinates {
(1.0, -2.0)
(0.8958907, -1.428344)
};
\addplot+ [no marks, solid, black]coordinates {
(0.5, -2.0)
(1.0, -2.0)
};
\addplot+ [no marks, solid, black]coordinates {
(0.3679869, -1.50956)
(0.8958907, -1.428344)
};
\addplot+ [no marks, solid, black]coordinates {
(1.5, -2.0)
(1.459981, -1.46053)
};
\addplot+ [no marks, solid, black]coordinates {
(1.0, -2.0)
(1.5, -2.0)
};
\addplot+ [no marks, solid, black]coordinates {
(0.8958907, -1.428344)
(1.459981, -1.46053)
};
\addplot+ [no marks, solid, black]coordinates {
(2.0, -2.0)
(2.0, -1.5)
};
\addplot+ [no marks, solid, black]coordinates {
(1.5, -2.0)
(2.0, -2.0)
};
\addplot+ [no marks, solid, black]coordinates {
(1.459981, -1.46053)
(2.0, -1.5)
};
\addplot+ [no marks, solid, black]coordinates {
(-2.0, 1.0)
(-1.429238, 0.8987769)
};
\addplot+ [no marks, solid, black]coordinates {
(-2.0, 1.5)
(-2.0, 1.0)
};
\addplot+ [no marks, solid, black]coordinates {
(-1.460039, 1.460656)
(-1.429238, 0.8987769)
};
\addplot+ [no marks, solid, black]coordinates {
(-2.0, 0.5)
(-1.509204, 0.3672318)
};
\addplot+ [no marks, solid, black]coordinates {
(-2.0, 1.0)
(-2.0, 0.5)
};
\addplot+ [no marks, solid, black]coordinates {
(-1.429238, 0.8987769)
(-1.509204, 0.3672318)
};
\addplot+ [no marks, solid, black]coordinates {
(-2.0, 0.0)
(-1.671448, -0.003451758)
};
\addplot+ [no marks, solid, black]coordinates {
(-2.0, 0.5)
(-2.0, 0.0)
};
\addplot+ [no marks, solid, black]coordinates {
(-1.509204, 0.3672318)
(-1.671448, -0.003451758)
};
\addplot+ [no marks, solid, black]coordinates {
(-2.0, -0.5)
(-1.510406, -0.3695582)
};
\addplot+ [no marks, solid, black]coordinates {
(-2.0, 0.0)
(-2.0, -0.5)
};
\addplot+ [no marks, solid, black]coordinates {
(-1.671448, -0.003451758)
(-1.510406, -0.3695582)
};
\addplot+ [no marks, solid, black]coordinates {
(-2.0, -1.0)
(-1.427553, -0.8944508)
};
\addplot+ [no marks, solid, black]coordinates {
(-2.0, -0.5)
(-2.0, -1.0)
};
\addplot+ [no marks, solid, black]coordinates {
(-1.510406, -0.3695582)
(-1.427553, -0.8944508)
};
\addplot+ [no marks, solid, black]coordinates {
(-2.0, -1.0)
(-2.0, -1.5)
};
\addplot+ [no marks, solid, black]coordinates {
(-1.427553, -0.8944508)
(-1.460145, -1.459283)
};
\addplot+ [no marks, solid, black]coordinates {
(2.0, -1.0)
(1.429215, -0.8987309)
};
\addplot+ [no marks, solid, black]coordinates {
(2.0, -1.5)
(2.0, -1.0)
};
\addplot+ [no marks, solid, black]coordinates {
(1.459981, -1.46053)
(1.429215, -0.8987309)
};
\addplot+ [no marks, solid, black]coordinates {
(2.0, -0.5)
(1.509184, -0.3672392)
};
\addplot+ [no marks, solid, black]coordinates {
(2.0, -1.0)
(2.0, -0.5)
};
\addplot+ [no marks, solid, black]coordinates {
(1.429215, -0.8987309)
(1.509184, -0.3672392)
};
\addplot+ [no marks, solid, black]coordinates {
(2.0, 0.0)
(1.671369, 0.003443455)
};
\addplot+ [no marks, solid, black]coordinates {
(2.0, -0.5)
(2.0, 0.0)
};
\addplot+ [no marks, solid, black]coordinates {
(1.509184, -0.3672392)
(1.671369, 0.003443455)
};
\addplot+ [no marks, solid, black]coordinates {
(2.0, 0.5)
(1.510172, 0.3695378)
};
\addplot+ [no marks, solid, black]coordinates {
(2.0, 0.0)
(2.0, 0.5)
};
\addplot+ [no marks, solid, black]coordinates {
(1.671369, 0.003443455)
(1.510172, 0.3695378)
};
\addplot+ [no marks, solid, black]coordinates {
(2.0, 1.0)
(1.426582, 0.8944412)
};
\addplot+ [no marks, solid, black]coordinates {
(2.0, 0.5)
(2.0, 1.0)
};
\addplot+ [no marks, solid, black]coordinates {
(1.510172, 0.3695378)
(1.426582, 0.8944412)
};
\addplot+ [no marks, solid, black]coordinates {
(2.0, 1.0)
(2.0, 1.5)
};
\addplot+ [no marks, solid, black]coordinates {
(1.426582, 0.8944412)
(1.459547, 1.459424)
};
\addplot+ [no marks, solid, black]coordinates {
(-3.8285686989269494e-16, -1.0)
(-0.3680867, -1.508369)
};
\addplot+ [no marks, solid, black]coordinates {
(-0.7071067811865475, -0.7071067811865476)
(-0.8988851, -1.428805)
};
\addplot+ [no marks, solid, black]coordinates {
(-0.7071067811865475, -0.7071067811865476)
(-1.427553, -0.8944508)
};
\addplot+ [no marks, solid, black]coordinates {
(-1.0, 1.2246467991473532e-16)
(-1.510406, -0.3695582)
};
\addplot+ [no marks, solid, black]coordinates {
(-1.0, 1.2246467991473532e-16)
(-1.509204, 0.3672318)
};
\addplot+ [no marks, solid, black]coordinates {
(-0.7071067811865475, 0.7071067811865476)
(-1.429238, 0.8987769)
};
\addplot+ [no marks, solid, black]coordinates {
(-0.7071067811865475, 0.7071067811865476)
(-0.8962192, 1.428906)
};
\addplot+ [no marks, solid, black]coordinates {
(6.123233995736766e-17, 1.0)
(-0.3693048, 1.511127)
};
\addplot+ [no marks, solid, black]coordinates {
(6.123233995736766e-17, 1.0)
(0.3667017, 1.509209)
};
\addplot+ [no marks, solid, black]coordinates {
(0.7071067811865476, 0.7071067811865475)
(0.8982224, 1.429057)
};
\addplot+ [no marks, solid, black]coordinates {
(0.7071067811865476, 0.7071067811865475)
(1.426582, 0.8944412)
};
\addplot+ [no marks, solid, black]coordinates {
(1.0, 0.0)
(1.510172, 0.3695378)
};
\addplot+ [no marks, solid, black]coordinates {
(1.0, 0.0)
(1.509184, -0.3672392)
};
\addplot+ [no marks, solid, black]coordinates {
(0.7071067811865476, -0.7071067811865475)
(1.429215, -0.8987309)
};
\addplot+ [no marks, solid, black]coordinates {
(0.7071067811865476, -0.7071067811865475)
(0.8958907, -1.428344)
};
\addplot+ [no marks, solid, black]coordinates {
(0.3679869, -1.50956)
(-3.8285686989269494e-16, -1.0)
};
\draw[very thick, grn] (axis cs:0, 0) circle[radius=1];
\end{axis}

\end{tikzpicture}
     \caption{Two-dimensional domain used for numerical results in
    \sref{sec:numerical:experiments}. The thick green line is the interface
    between the two subdomains $\Omega_{1}$ and $\Omega_{2}$. The thin black
    lines show the finite difference block interfaces.\label{fig:square:circle}}
  \end{center}
\end{figure}
We consider the two-dimensional square domain $\Omega = [-2, 2]^2$.
Inside of $\Omega$ we define the unit circle $\Gamma_{I} = \{(x_{1}, x_{2}) |
x_{1}^2+x_{2}^2=1 \}$ to partition the domain into a closed unit disk
$\Omega_{1} = \{(x_{1}, x_{2}) | x_{1}^2+x_{2}^2\le1 \}$ and the remainder
$\Omega_2 = \text{cl}(\Omega \setminus \Omega_{1})$.
The interface $\Gamma_{I}$ is governed by the nonlinear condition
\begin{align}
  \tau^{\pm} = \beta \arcsinh\left(V^{\pm}\right) + g_{\tau}^{\pm},
\end{align}
where $\beta > 0$ and $g_{\tau}^{\pm}$ is a time and space dependent forcing
function; around $V = 0$ with $g_{\tau}^{\pm} = 0$ the linearization of the
interface condition is $\tau^{\pm} = \beta V$.
The right and left boundaries of $\Omega$ are taken to be Dirichlet, the top
and bottom boundaries Neumann; the Dirichlet and Neumann boundary conditions are
enforced using the standard approach described in the
\aref{app:standard:approach}.
As shown in Figure~\ref{fig:square:circle}, the domain is decomposed into 56
finite difference blocks and the all interfaces, nonlinear and computational,
are enforced using the characteristic approach described in
\sref{subsec:characteristic:interface} through the introduction of the auxiliary
variable $\vv{u}^{*f}$ on each interface.
Given the unstructured connectivity of the blocks it is necessary to use the
same $(N+1) \times (N+1)$ grid of points in each block; we refer to $N$ as the
block size.
Time stepping is performed using the low-storage, fourth order Runge-Kutta
scheme of \citet[(5,4) $2N$-Storage RK scheme, solution
$3$]{CarpenterKennedy1994}.

In order to assess the stiffness and accuracy of the scheme in two spatial
dimensions we use the method of manufactured solutions (MMS) \citep{Roache}.
In particular, we assume an analytic solution and compute the necessary boundary,
interface, and volume data.
The manufactured solution is taken to be
\begin{equation}\label{eqn:mms}
  u(x_{1}, x_{2}, t) =
  \begin{cases}
    \sin(t) \frac{e}{1+e}\left( 2 - e^{-r^2}\right) r \sin(\theta), &\quad
    (x_{1}, x_{2}) \in \Omega_1,\\
    \sin(t) \left({(r- 1)}^2 \cos(\theta) + (r - 1)\sin(\theta)\right), &\quad
    (x_{1}, x_{2}) \in \Omega_2,
  \end{cases}
\end{equation}
where $r^2 = x_{1}^2 + x_{2}^2$ and $\theta = \atan2(x_{2}, x_{1})$.
The boundary, interface, and forcing data are found by using assumed
solution~\eqref{eqn:mms} in governing equations~\eqref{eqn:wave}.
In order to avoid order reduction with time dependent data, we found it
necessary to define the Dirichlet boundary data by integrating $\dot{g}_{D}$
using the Runge-Kutta method.
Solution~\eqref{eqn:mms} satisfies force balance along $\Gamma_{I}$, i.e.,
continuity of traction $\tau$, and the interface data $g_{\tau}^{\pm}$ is used
to enforce the assumed solution.
In the MMS test the material properties are $\rho = 1$ and $C_{ij} =
\delta_{ij}$, with $\delta_{ij}$ being the kroneckor delta; after the mesh
warping the  effective material parameters $\hat{C}_{ij}$ are spatially varying.

\begin{table}
  \centering
  \sisetup{
    round-mode = places,
    round-precision = 2
  }%
  \begin{tabular}{rcccccc}
    \toprule
    & \multicolumn{3}{c}{characteristic} & \multicolumn{3}{c}{non-characteristic}\\
    \cmidrule(lr){2-4}
    \cmidrule(lr){5-7}
    $\beta$ &
    $\kappa$ & $\|\Delta\VV{u}\|_{H}$ with $\kappa$ & (with $2\kappa$) &
    $\kappa$ & $\|\Delta\VV{u}\|_{H}$ with $\kappa$ & (with $2\kappa$)
    \\
    \midrule
      $1$ & $1 / 2$ & \num{1.3932094993896921e-9} & (\num{1.6372160477702888e10}) & $1 / 2\phantom{^1}$ & \num{1.232508736634718e-9} & (\num{1.6372160477702888e10}) \\
      $4$ & $1 / 2$ & \num{1.3883082787068281e-9} & (\num{1.6372160477702888e10}) & $1 / 2\phantom{^1}$ & \num{1.242230995069987e-9} & (\num{1.6372160477702888e10}) \\
     $16$ & $1 / 2$ & \num{1.3870582193742948e-9} & (\num{1.6372160477702888e10}) & $1 / 2^3$           & \num{1.281850798312621e-9} & (\num{1.9012182973780944e-2}) \\
     $64$ & $1 / 2$ & \num{1.3882634681810410e-9} & (\num{1.6372160477702888e10}) & $1 / 2^5$           & \num{1.339934930623836e-9} & (\num{2.5132868682904300e-2}) \\
    $128$ & $1 / 2$ & \num{1.3886536855236286e-9} & (\num{1.6372160477702888e10}) & $1 / 2^6$           & \num{1.357500106402760e-9} & (\num{2.6216859792538010e-2}) \\
    \bottomrule
  \end{tabular}
  \caption{Stable Courant $\kappa$ for the characteristic and non-characteristic
  methods for increasing values of $\beta$ using the SBP operator with interior
  accuracy $2p = 6$. Shown also are the L$^{2}$-errors for the stable Courant
  number $\kappa$ and the unstable Courant number $2\kappa$. As seen, the
  characteristic method time step is independent of $\beta$ and the
  non-characteristic method requires a Courant that scales inversely with
$\beta$.\label{tab:stiffness:2D}}
\end{table}
To compare the stiffness of the standard and characteristic nonlinear interface
treatment we vary the nonlinear interface parameter $\beta$ and decrease the
time step size until the simulation is stable for a fixed block size $N = 48$.
For a non-stiff method, the time step size should be on the order of the
effective grid spacing for all $\beta > 0$.
In particular, we define the time step size to be
\begin{equation}
  \Delta t = \kappa \bar{h},
\end{equation}
where $\kappa$ is the Courant number and a non-stiff scheme should have $\kappa
\sim 1$; since the material properties are taken to be unity the wave speed in
this problem is $1$.
The effective grid spacing $\bar{h}$ is defined as
\begin{equation}
  \bar{h} = \min(\bar{h}_1, \bar{h}_2),
  ~
  \bar{h}_{r} = \frac{1}{N}
  \sqrt{
    {\left(\hat{\partial}_{r} x_{1}\right)}^{2}
    +
    {\left(\hat{\partial}_{r} x_{2}\right)}^{2}
  }.
\end{equation}
Table~\ref{tab:stiffness:2D} gives the Courant number $\kappa$ required for
stability of the two methods with various values of $\beta$ using SBP interior
order $2p = 6$.
Here the value of $\kappa$ was repeatedly halved until the error in the
simulation at time $t = 0.1$ no longer decreased dramatically.
To demonstrate that the stable Courant $\kappa$ is close to its maximum value,
Table~\ref{tab:stiffness:2D} also reports the L$^{2}$-error with a time step
defined by $\kappa$ and $2\kappa$, and as can be seen the former time step leads
to an accurate simulation and the latter an inaccurate one.
As can be seen the characteristic method requires a similar time step for all
values of the parameter $\beta$ whereas the non-characteristic method requires
a significantly reduced time step as $\beta$ increases.
Though not shown, results with SBP interior orders $2p = 2$ and $2p = 4$ are
similar; for $2p = 2$ the characteristic method can use a Courant of $\kappa =
1$ for all values of $\beta$ as can the non-characteristic method with $\beta
= 1$.

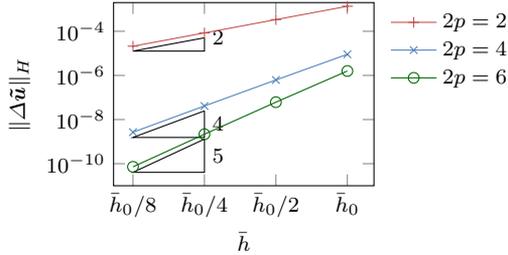
\begin{figure}
  \centering
\begin{tikzpicture}
\begin{loglogaxis}[
    ymax={0.002},
    ymin={1.0e-11},
    xmin={0.002},
    xmax={0.025},
    xlabel={$\bar{h}$},
    ylabel={$\|\Delta \VV{u}\|_{H}$},
    width=5cm,
    height=4cm,
    ytick={1e-2, 1e-4, 1e-6, 1e-8, 1e-10},
    xtick={0.01931342688192358, 0.00965671344096179, 0.004828356720480895, 0.0024141783602404476},
    xticklabels={
      $\bar{h}_{0}$,
      $\bar{h}_{0} / 2$,
      $\bar{h}_{0} / 4$,
      $\bar{h}_{0} / 8$
    },
    legend pos=outer north east,
    legend style={draw=none},
  ]
    \addplot[color={black}, forget plot]
        table[row sep={\\}]
        {
            \\
            0.00482836  5.07274e-5 \\
            0.00241418  1.26819e-5 \\
            0.00482836  1.26819e-5 \\
            0.00482836  5.07274e-5 \\
        }
        ;
    \addplot[color={black}, forget plot]
        table[row sep={\\}]
        {
          \\
          0.00482836  2.4726e-8  \\
          0.00241418  1.54538e-9 \\
          0.00482836  1.54538e-9 \\
          0.00482836  2.4726e-8  \\
        }
        ;
        \addplot[color={black}, forget plot]
        table[row sep={\\}]
        {
          \\
          0.00482836  1.29807e-9  \\
          0.00241418  4.05645e-11 \\
          0.00482836  4.05645e-11 \\
          0.00482836  1.29807e-9  \\
        }
        ;
        \node[anchor=west] () at (axis cs:0.00482836, 5.1695413859064964e-5){$2$};
        \node[anchor=west] () at (axis cs:0.00482836, 6.181509999991914e-9){$4$};
        \node[anchor=west] () at (axis cs:0.00482836, 2.29467994533007e-10){$5$};

        \addplot[color={red}, mark={+}]
        table[row sep={\\}]
        {
          \\
          0.01931342688192358    0.0013636640526359927  \\
          0.00965671344096179    0.00033890400487821703 \\
          0.004828356720480895   8.454569807108786e-5   \\
          0.0024141783602404476  2.1120460836678772e-5  \\
        }
        ;
        \addlegendentry{$2p = 2$}
        \addplot[color={blu}, mark={x}]
        table[row sep={\\}]
        {
          \\
          0.01931342688192358    9.081431805704101e-6  \\
          0.00965671344096179    6.221082179393138e-7  \\
          0.004828356720480895   4.121000453744455e-8  \\
          0.0024141783602404476  2.6732076054102486e-9 \\
        }
        ;
        \addlegendentry{$2p = 4$}
        \addplot[color={grn}, mark={o}]
        table[row sep={\\}]
        {
          \\
          0.01931342688192358    1.5858108178976352e-6 \\
          0.00965671344096179    6.134934113328974e-8  \\
          0.004828356720480895   2.1634424243046745e-9 \\
          0.0024141783602404476  7.200537658686824e-11 \\
        }
        ;
        \addlegendentry{$2p = 6$}
  \end{loglogaxis}
\end{tikzpicture}

   \caption{Convergence results for MMS solution~\eqref{eqn:mms} with $\beta =
    128$ using SBP interior orders $2p = 2, 4, 6$ with the characteristic
    nonlinear interface treatment. The value of $\bar{h}_{0} \approx 0.019$
    corresponds to block size $N = 17$.\label{fig:accuracy:2D}}
\end{figure}
To investigate the convergence of the two-dimensional, characteristic method we
now run the same MMS solution~\eqref{eqn:mms} to time $t_f = 1$ using $\beta
= 128$ with different levels of refinement and a fixed Courant number $\kappa
= 1/2$.
Figure~\ref{fig:accuracy:2D} shows the convergence of the scheme using mesh
levels $N = 17 \times 2^r$ where $r = 0, 1, 2, 3$.
As can be seen the convergence order is similar to the one-dimensional case.

\begin{figure}
\centering
\begin{subfigure}[t]{0.3\textwidth}
  \begin{tikzpicture}
    \node[inner sep=0pt, anchor=south west] (wave1) at (0,0)
    {\includegraphics[width=3.4cm, trim=241 241 241 241, clip]{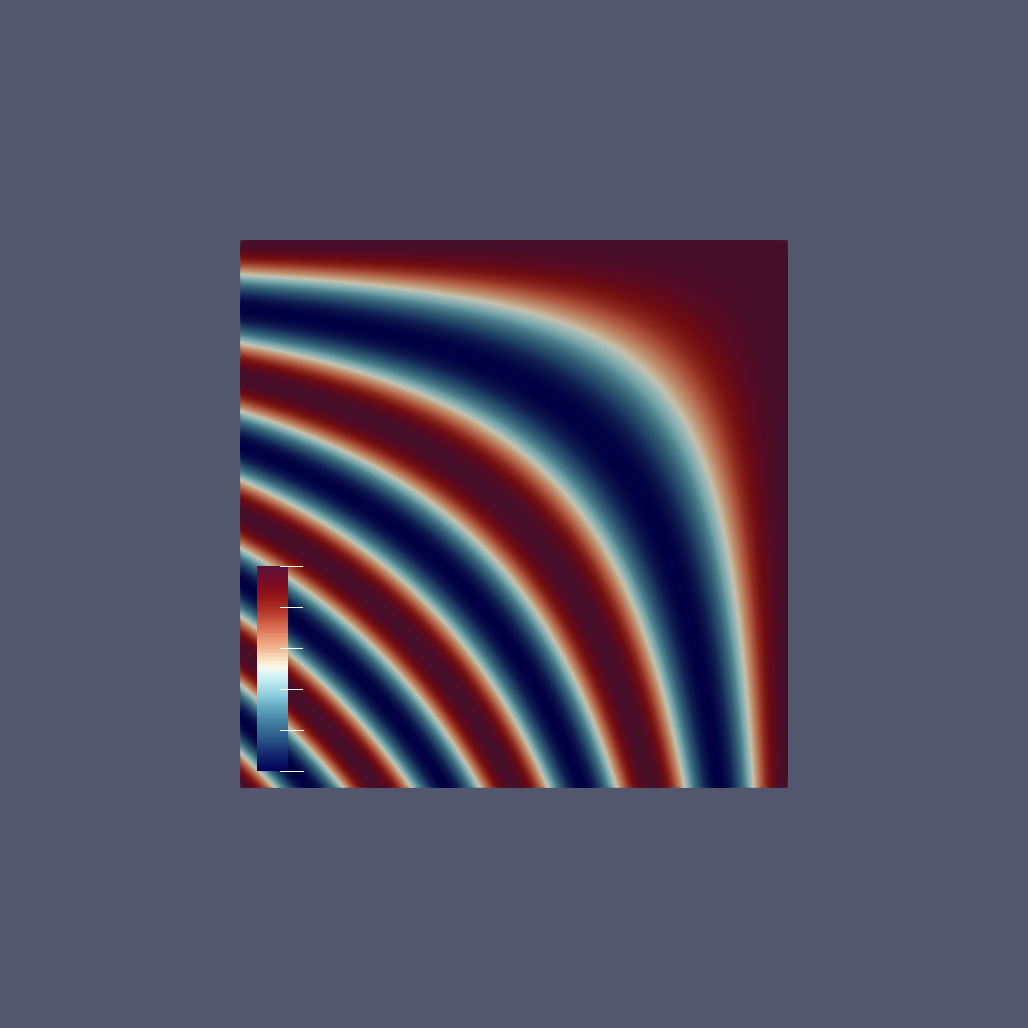}};
    \begin{axis}[
        xmin = {-2},
        xmax = {2},
        ymax = {2},
        width=5cm,
        height=5cm,
        ticks=none,
        ymin = {-2}
      ]
      \draw[very thin, grn] (axis cs:0, 0) circle[radius=1];
      \node[anchor=west] () at (axis cs:-1.5, -0.37){\tiny \color{white} $0.5$};
      \node[anchor=west] () at (axis cs:-1.5, -0.67){\tiny \color{white} $0.6$};
      \node[anchor=west] () at (axis cs:-1.5, -0.97){\tiny \color{white} $0.7$};
      \node[anchor=west] () at (axis cs:-1.5, -1.27){\tiny \color{white} $0.8$};
      \node[anchor=west] () at (axis cs:-1.5, -1.57){\tiny \color{white} $0.9$};
      \node[anchor=west] () at (axis cs:-1.5, -1.87){\tiny \color{white} $1.0$};
    \end{axis}
  \end{tikzpicture}
  \caption{material parameter $c_{11}$}
\end{subfigure}
\hfill
\begin{subfigure}[t]{0.3\textwidth}
  \begin{tikzpicture}
    \node[inner sep=0pt, anchor=south west] (wave1) at (0,0)
    {\includegraphics[width=3.4cm, trim=241 241 241 241, clip]{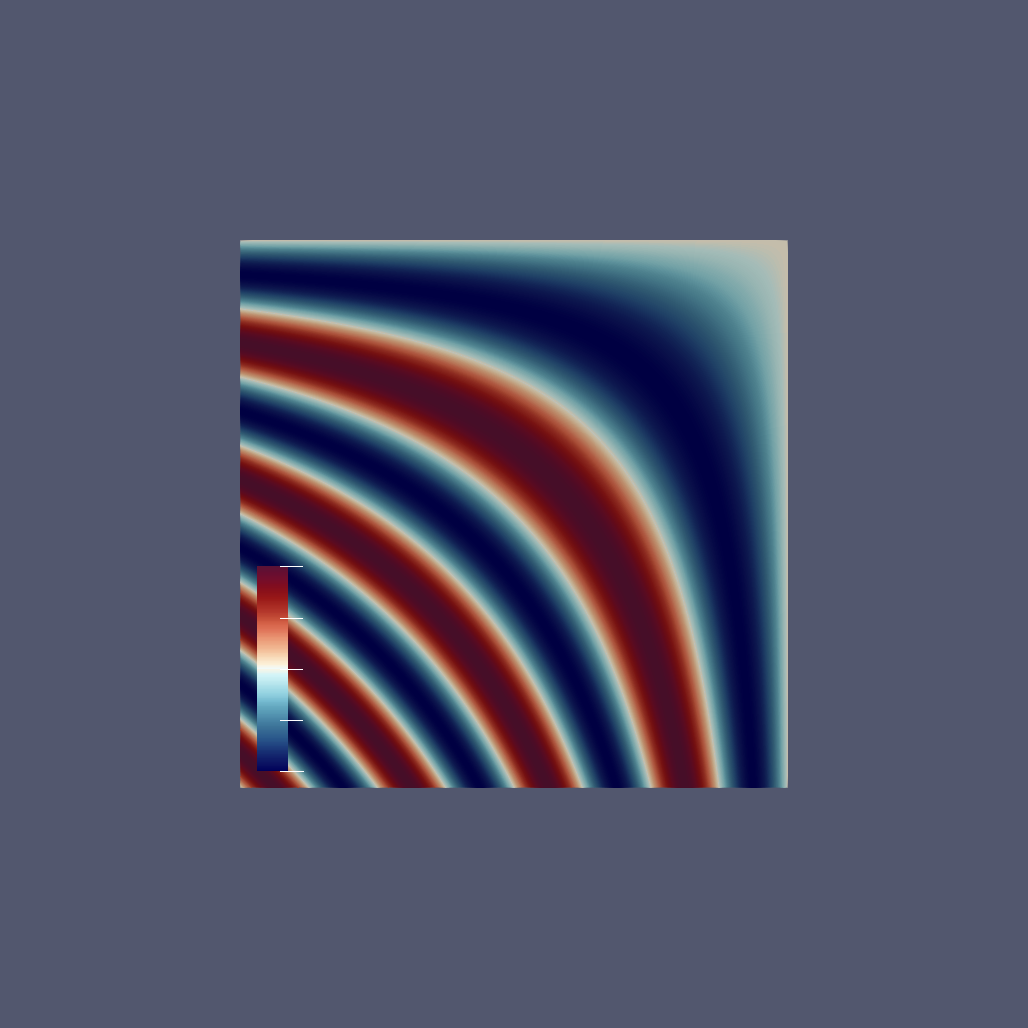}};
    \begin{axis}[
        xmin = {-2},
        xmax = {2},
        ymax = {2},
        width=5cm,
        height=5cm,
        ticks=none,
        ymin = {-2}
      ]
      \draw[very thin, grn] (axis cs:0, 0) circle[radius=1];
      \node[anchor=west] () at (axis cs:-1.5, -0.370){\tiny \color{white} $-0.25$};
      \node[anchor=west] () at (axis cs:-1.5, -0.745){\tiny \color{white} $-0.125$};
      \node[anchor=west] () at (axis cs:-1.5, -1.120){\tiny \color{white} $\phantom{-}0.0$};
      \node[anchor=west] () at (axis cs:-1.5, -1.495){\tiny \color{white} $\phantom{-}0.125$};
      \node[anchor=west] () at (axis cs:-1.5, -1.870){\tiny \color{white} $\phantom{-}0.25$};
    \end{axis}
  \end{tikzpicture}
  \caption{material parameter $c_{12}$}
\end{subfigure}
\hfill
\begin{subfigure}[t]{0.3\textwidth}
  \begin{tikzpicture}
    \node[inner sep=0pt, anchor=south west] (wave1) at (0,0)
    {\includegraphics[width=3.4cm, trim=241 241 241 241, clip]{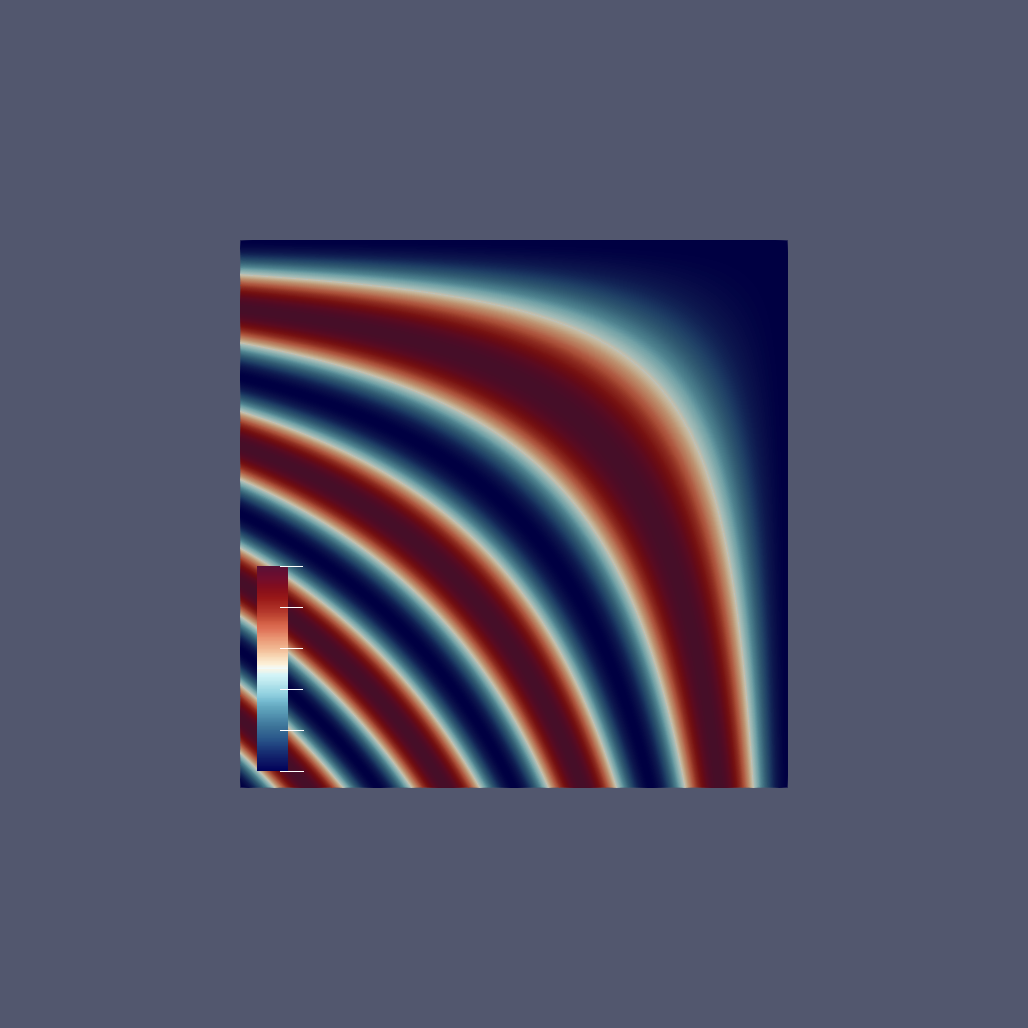}};
    \begin{axis}[
        xmin = {-2},
        xmax = {2},
        ymax = {2},
        width=5cm,
        height=5cm,
        ticks=none,
        ymin = {-2}
      ]
      \draw[very thin, grn] (axis cs:0, 0) circle[radius=1];
      \node[anchor=west] () at (axis cs:-1.5, -0.37){\tiny \color{white} $0.5$};
      \node[anchor=west] () at (axis cs:-1.5, -0.67){\tiny \color{white} $0.6$};
      \node[anchor=west] () at (axis cs:-1.5, -0.97){\tiny \color{white} $0.7$};
      \node[anchor=west] () at (axis cs:-1.5, -1.27){\tiny \color{white} $0.8$};
      \node[anchor=west] () at (axis cs:-1.5, -1.57){\tiny \color{white} $0.9$};
      \node[anchor=west] () at (axis cs:-1.5, -1.87){\tiny \color{white} $1.0$};
    \end{axis}
  \end{tikzpicture}
  \caption{material parameter $c_{22}$}
\end{subfigure}

\begin{subfigure}[t]{0.3\textwidth}
  \begin{tikzpicture}
    \node[inner sep=0pt, anchor=south west] (wave1) at (0,0)
    {\includegraphics[width=3.4cm, trim=241 241 241 241, clip]{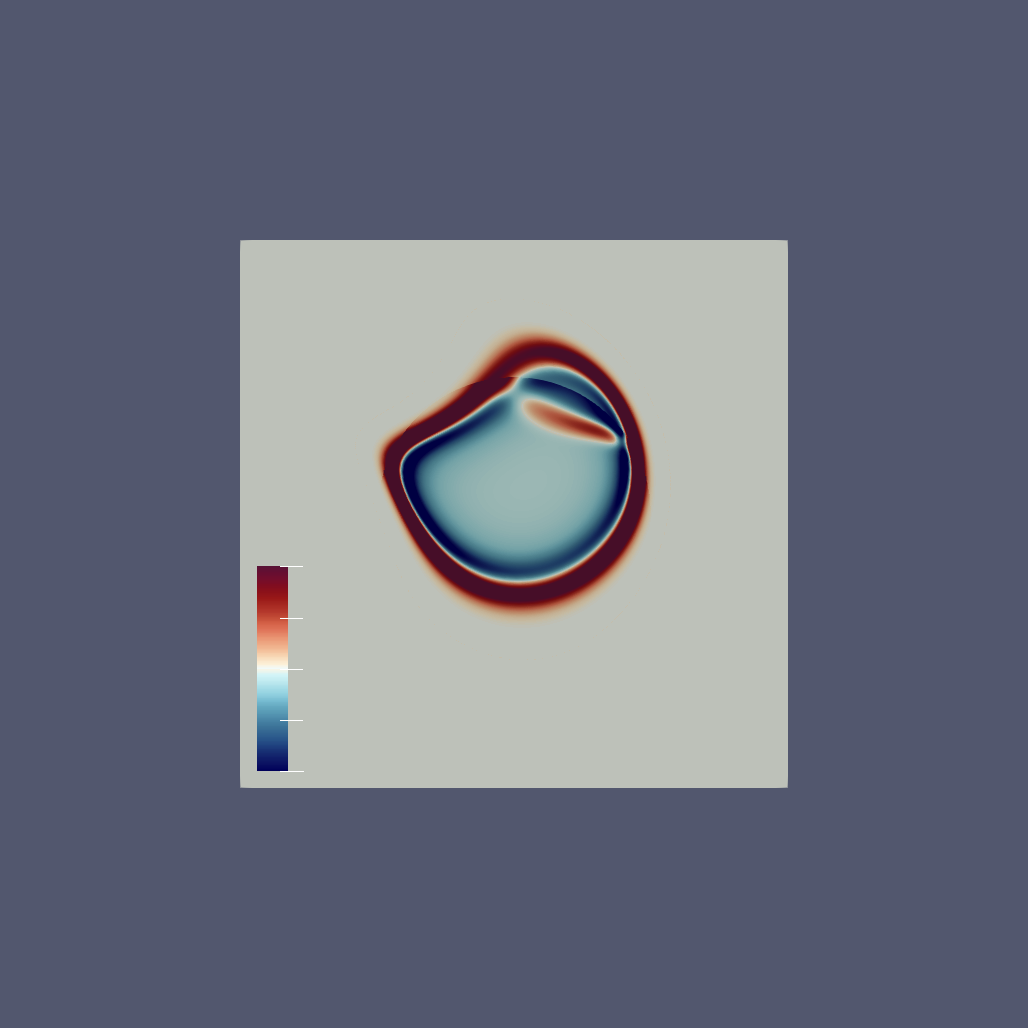}};
    \begin{axis}[
        xmin = {-2},
        xmax = {2},
        ymax = {2},
        width=5cm,
        height=5cm,
        ticks=none,
        ymin = {-2}
      ]
      \draw[very thin, grn] (axis cs:0, 0) circle[radius=1];
      \node[anchor=west] () at (axis cs:-1.5, -0.370){\tiny \color{white} $-0.05$};
      \node[anchor=west] () at (axis cs:-1.5, -0.745){\tiny \color{white} $-0.025$};
      \node[anchor=west] () at (axis cs:-1.5, -1.120){\tiny \color{white} $\phantom{-}0.0$};
      \node[anchor=west] () at (axis cs:-1.5, -1.495){\tiny \color{white} $\phantom{-}0.025$};
      \node[anchor=west] () at (axis cs:-1.5, -1.870){\tiny \color{white} $\phantom{-}0.05$};
    \end{axis}
  \end{tikzpicture}
  \caption{displacement $u$ at $t = 1$}
\end{subfigure}
\hfill
\begin{subfigure}[t]{0.3\textwidth}
  \begin{tikzpicture}
    \node[inner sep=0pt, anchor=south west] (wave1) at (0,0)
    {\includegraphics[width=3.4cm, trim=241 241 241 241, clip]{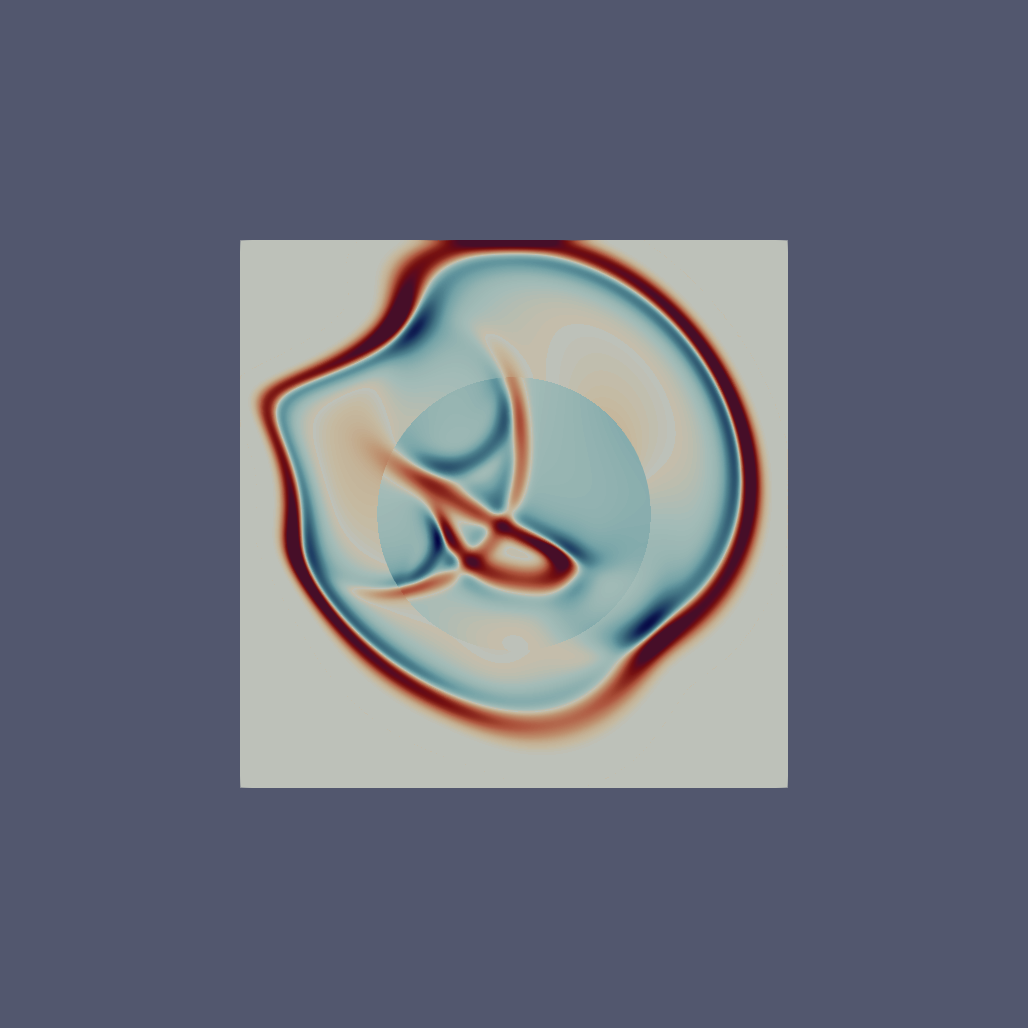}};
    \begin{axis}[
        xmin = {-2},
        xmax = {2},
        ymax = {2},
        width=5cm,
        height=5cm,
        ticks=none,
        ymin = {-2}
      ]
      \draw[very thin, grn] (axis cs:0, 0) circle[radius=1];
    \end{axis}
  \end{tikzpicture}
  \caption{displacement $u$ at $t = 2$}
\end{subfigure}
\hfill
\begin{subfigure}[t]{0.3\textwidth}
  \begin{tikzpicture}
    \node[inner sep=0pt, anchor=south west] (wave1) at (0,0)
    {\includegraphics[width=3.4cm, trim=241 241 241 241, clip]{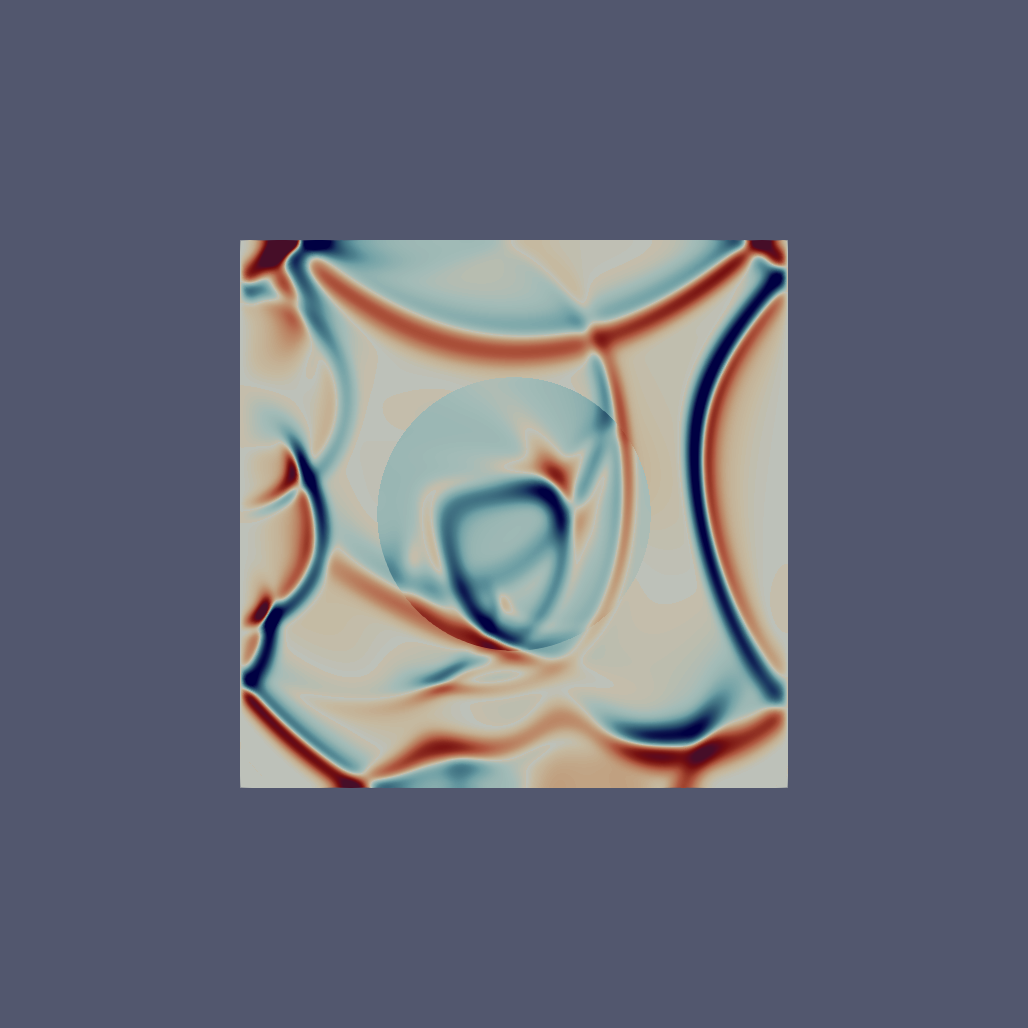}};
    \begin{axis}[
        xmin = {-2},
        xmax = {2},
        ymax = {2},
        width=5cm,
        height=5cm,
        ticks=none,
        ymin = {-2}
      ]
      \draw[very thin, grn] (axis cs:0, 0) circle[radius=1];
    \end{axis}
  \end{tikzpicture}
  \caption{displacement $u$ at $t = 3$}
\end{subfigure}
   \caption{Variable material parameters $C_{ij}$ and displacement field $u$ snap
  shots with block size $N = 136$ with the mesh show in
  Figure~\ref{fig:square:circle}. The colormap for the displacement field is
  saturated to show features at later times and the green curve indicates the
  location of the nonlinear interface.\label{fig:wavefield}}
\end{figure}
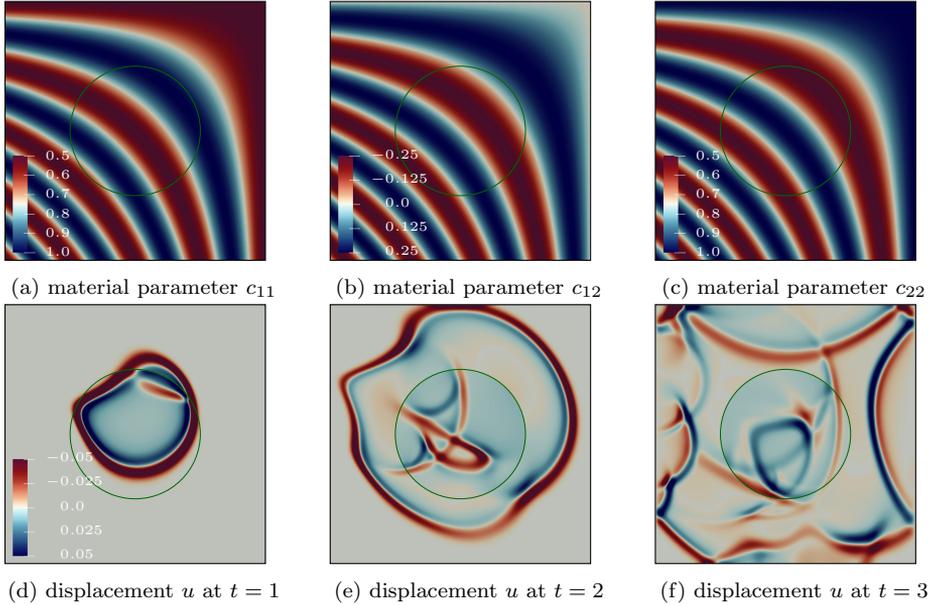
As a final test, we explore the self-convergence and energy dissipation
properties of the characteristic method with variable material properties and no
body or boundary data.
The same two-dimensional spatial domain is used, but now the material parameters
are taken to be
\begin{subequations}
  \begin{align}
    \rho &= 1,\\
    C_{11} &= \cos(\theta)^2 + \frac{1}{2} \sin(\theta)^2,\\
    C_{12} &= -\frac{1}{2} \cos(\theta) \sin(\theta),\\
    C_{22} &= \sin(\theta)^2 + \frac{1}{2} \cos(\theta)^2,
  \end{align}
\end{subequations}
where the angle $\theta = \frac{\pi}{4} \left(2 - x_{1}\right)\left(2 -
x_{2}\right)$; colormaps of the material parameters are shown in
Figure~\ref{fig:wavefield}.
The Courant number $\kappa = 1/2$ is used for all the simulations and the
material parameters lead to a maximum wave speed of $1$, i.e., maximum
eigenvalue of the matrix defined by $C_{ij} / \rho$.
The initial displacement is taken to be the product of two off-center Gaussians
\begin{equation}
  u_{0} =
  \exp\left(
  -\frac{{\left(x_{1} - \mu_{1}\right)}^{2}}{2\sigma_{1}}
  -\frac{{\left(x_{2} - \mu_{2}\right)}^{2}}{2\sigma_{2}}
  \right),
\end{equation}
where $\mu_1 = 0.1$,  $\mu_2 =  0.2$,  $\sigma_1 = 0.0025$, and $\sigma_2 =
0.005$, and the initial velocity is $\dot{u}_{0} = 0$.
A nonlinear parameter $\beta = 1$ is used in order to highlight the effect of the
nonlinear interface condition; larger values of $\beta$ lead to a more
continuous solution across the interface since the sliding velocity $V$ will be
lower.
Snapshots of the displacement field at various times are shown in
Figure~\ref{fig:wavefield} for the block size $N = 136 = 17 \times 8$ and SBP
interior order $2p = 6$.
As can be seen in the figure, there is a discontinuity in the displacement
across the interface as well as reflected waves.

For the self-convergence study we run the simulation until time $t = 3$ using
$N_{r} = 17 \times 2^r$ with $r = 1, 2, 3$.
The error is estimated by taking the difference between neighboring resolutions,
and the rate is estimated by
\begin{equation}
  \text{rate} =
  \log_{2}(\|\Delta_{1}\|_{H_{1}})
  -
  \log_{2}(\|\Delta_{2}\|_{H_{2}}),
\end{equation}
where $\Delta_{r}$ is the difference between the solutions using $N_{r}$ and
$N_{r+1}$ and $H_{r}$ indicates that the norm is taken with respect to the
metrics defined by $N_{r}$.
With this, we get an estimate convergence rate for this problem of $4.4$ using
the SBP operators with interior accuracy $2p = 6$.

\begin{figure}
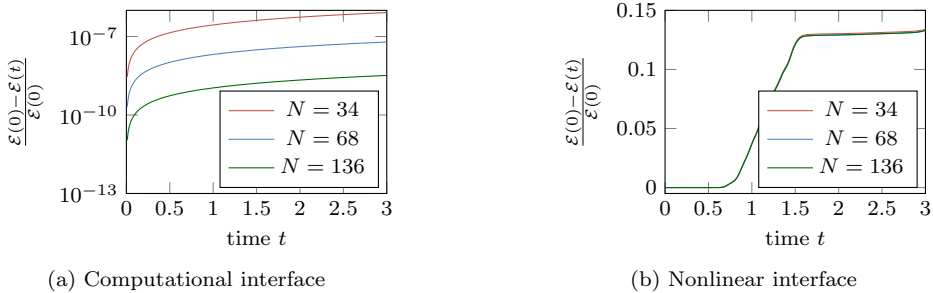

  \begin{subfigure}[t]{0.4\textwidth}
    \centering

     \caption{Computational interface}
  \end{subfigure}
  \hfill
  \begin{subfigure}[t]{0.4\textwidth}
    \centering
%
     \caption{Nonlinear interface}
  \end{subfigure}
  \caption{Normalized dissipated energy for a computational and nonlinear
    interface $\Gamma_{I}$ with energy is computed discretely
    using~\eqref{eqn:disc:energy} and positive values indicating
  dissipation.\label{fig:time:energy}}
\end{figure}

Using same material properties and initial condition,
Figure~\ref{fig:time:energy} show the dissipated energy when $\Gamma_I$ is taken
to be a computational interface and a nonlinear interface with $\beta = 1$;
energy is measured using the discrete energy norm~\eqref{eqn:disc:energy}.
In both cases the energy decreases in time as the theory predicts.
In the case of the computational interface the dissipation is purely numerical,
and as the results show the dissipation decreases as the resolution increases.
In the case of the nonlinear interface the amount of energy dissipated is
larger since the continuous formulation supports energy dissipation on interface
$\Gamma_{I}$.

\section{Concluding Remarks}\label{sec:conclusions}

We have developed a characteristic based method for handling boundary and
interface conditions with SBP finite difference methods for the second order,
scalar wave equation.
The key idea of the method is the introduction of an additional unknown on the
block boundaries which evolves in time and acts as local Dirichlet data for the
block.
The rate of change of the boundary unknown is defined in an upwind fashion that
modifies the incoming characteristic variable, which is similar to the technique
previously used to remove stiffness for the wave equation in first order form
with nonlinear interfaces \citep{KozdonDunhamNordstrom2012}.

The main benefit of the scheme is that, when compared with the standard
approach \citep{VirtaMattsson2014, Duru2019}, the scheme is non-stiff for all
characteristic boundary conditions and a class of nonlinear interface
conditions that can be written in characteristic form;
we note that at the continuous level the equation we consider are the same as
\citet{Duru2019} and that our schemes only differ at the discrete level.
One benefit of this approach is that it enables the use of a wider class of time
stepping methods for earthquake rupture problems with nonlinear interfaces.

The energy method was used to show that the proposed scheme is stable.
Numerical experiments showed that the proposed scheme is non-stiff, confirmed
the stability results, and also demonstrated the accuracy of the scheme. The
analysis presented is dimension independent, thus the results equally apply to
three dimensions. That said, the penalty parameter scales with $d$, and thus
there may be more restrictive time step in three dimensions;
see~\eqref{eqn:penalty:param}.

One area for future work includes more general wave equations, such as linear
elasticity. In this case, there will be multiple displacement and characteristic
variables. One important question to consider is whether auxiliary interface
variables are required for all components or only a subset.
The work \citet{AlmquistDunham2021} will be relevant to these exstensions,
particularly if one is willing to be restricted to fully-compatible SBP
operators.

One of the disadvantages of the traditional SBP finite difference formulations
is that the mesh must be conforming across block interfaces. Computationally,
this means that some regions of the domain will have finer mesh resolution than
the physics dictates which increases the computational cost through both
increased overhead per time step and reduced time steps size. These limitations
have lead to a recent interest in SBP-SAT methods for non-conforming interfaces
\citep{mattsson2010stable, wang2016high, kozdon2016stable}. Though we see no
obvious reason that the methods presented here would not extend to the
non-conforming interface case, it remains for this to be rigorously
demonstrated.

\section*{Acknowledgements}
We thank the two anonymous reviewers whose insightful feedback substantially
improved the paper.

\appendix\label{appendix}

\section{Definition of Two-Dimensional SBP Operators}\label{app:2D:sbp}

As an example of how to construct multidimensional SBP operators, we consider
the two dimensional SBP finite difference operators.
We describe the operators on the reference block $\hat{B} = [0,1] \times [0,1]$,
where faces $1$ and $2$ are the right and left faces with faces $3$ and $4$
being the top and bottom faces, respectively.
For simplicity we let the domain $\hat{B}$ be discretized with an $(N+1) \times
(N+1)$ grid points with the grid nodes located at $\dof{\vv{\xi}}{kl} = (kh,
lh)$ for $0 \leq k,l \leq N$ with $h = 1/N$.
The projection of $u$ onto the grid is denoted $\VV{u}$, where $\dof{\VV{u}}{kl}
\approx u(kh, lh)$ and is stored as a vector with with $\xi_1$ being the fastest
index; see~\eqref{eqn:dof:ordering}.
With this, the volume norm matrix can be written as
\begin{equation}
  \MM{H} = \mm{H} \otimes \mm{H}.
\end{equation}
We define the face restriction operators as
\begin{equation}
  \Mm{L}^{1} = \mm{I} \otimes \vv{e}_{0}^T, \qquad
  \Mm{L}^{2} = \mm{I} \otimes \vv{e}_{N}^T, \qquad
  \Mm{L}^{3} = \vv{e}_{0}^T \otimes \mm{I}, \qquad
  \Mm{L}^{4} = \vv{e}_{0}^T \otimes \mm{I},
\end{equation}
where the $\mm{I}$ is the  $(N+1) \times (N+1)$ identity matrix.
More generally the restriction to a single grid line in the $\xi_1$ and $\xi_2$
directions, respectively, are
\begin{equation}
  \Mm{L}_{l:} = \vv{e}_{l}^T \otimes \mm{I}, \qquad
  \Mm{L}_{:l} = \mm{I} \otimes \vv{e}_{l}^T.
\end{equation}
In order to construct $\MM{A}_{ii}^{(C)}$, no summation over $i$, we
construct individual one-dimensional second derivative matrices for each grid
line with varying coefficients $C$ and place them in the correct block;
expanding a single second derivative matrix with the tensor product and the
identity matrix only works in the constant coefficient case.
To do this it is useful to define $\VV{C}$ as the projection of $C$ onto the
grid with the coefficients along the individual grid lines being
\begin{equation}
  \mm{C}_{:l} = \text{diag}(C_{0l},\hdots,C_{Nl}), \qquad
  \mm{C}_{k:} = \text{diag}(C_{k0},\hdots,C_{kN}).
\end{equation}
The second derivative operators are the sum of the operators along each grid
line
\begin{subequations}
  \begin{align}
    \MM{A}_{11}^{(C)} &= (\mm{H} \otimes \mm{I}) \left[\sum_{l=0}^N
    \Mm{L}_{:l}^T \mm{A}_{11}^{\left(C_{:l}\right)} \Mm{L}_{:l} \right], &
    \MM{A}_{22}^{(C)} &= (\mm{I} \otimes \mm{H}) \left[\sum_{k=0}^N
    \Mm{L}_{k:}^T \mm{A}_{11}^{\left(C_{k:}\right)} \Mm{L}_{k:} \right],
  \end{align}
  and a tensor product is used for the mixed derivative operators
  \begin{align}
    \MM{A}_{12}^{(C)} &= (\mm{I} \otimes \mm{Q}^T) \MM{C}(\mm{Q} \otimes \mm{I}), &
    \MM{A}_{21}^{(C)} &= (\mm{Q}^T \otimes \mm{I}) \MM{C}(\mm{I} \otimes \mm{Q}).
  \end{align}
\end{subequations}

The boundary derivatives parallel to a face are given by the first derivative
operator $\mm{D}_{1}$ and those perpendicular with the boundary derivative
operators from the SBP definition~\ref{def:sbp2}:
\begin{subequations}
  \begin{align}
    \Mm{B}^{1}_{1} &= \mm{I} \otimes \vv{b}_0^T, &
    \Mm{B}^{1}_{2} &= \mm{D}_{1} \otimes \vv{e}_{0}^{T},\\
    \Mm{B}^{2}_{1} &= \mm{I} \otimes \vv{b}_N^T,&
    \Mm{B}^{2}_{2} &= \mm{D}_{1} \otimes \vv{e}_{N}^{T},\\
    \MM{B}^{3}_{1} &= \vv{e}_{0}^{T} \otimes \mm{D}_{1},&
    \MM{B}^{3}_{2} &= \vv{b}_0^T \otimes \mm{I},\\
    \MM{B}^{4}_{1} &= \vv{e}_{N}^{T} \otimes \mm{D}_{1}, &
    \MM{B}^{4}_{2} &= \vv{b}_N^T \otimes \mm{I}.
  \end{align}
\end{subequations}

\section{Proof of Theorem~\ref{thm:disc:energy:seminorm}}\label{app:disc:energy:seminorm}
To show that energy~\eqref{eqn:disc:block:energy} is positive we need the following definition from
\citet[Definition 2.4]{Mattsson2012}:
\begin{equation}
  \label{eqn:Mattsson:Def24}
  \MM{A}_{ij}^{(c)} = \MM{D}_{i}^{T} \MM{C} \MM{H} \MM{D}_{j} +
  \MM{R}_{ij}^{(c)}.
\end{equation}
The remainder matrix $\MM{R}_{ij}^{(c)}$ is symmetric positive
semidefinite if the coefficient $c$ is always positive; the remainder matrix is
zero when $i \ne j$.  The remainder matrix can be further decomposed using the
borrowing lemma from \citet[Lemma 1]{AlmquistDunham2020}:
\begin{equation}
  \label{eqn:AlmquistDunham:Lem1}
  \MM{R}_{ii}^{(c)}
  =
  \MM{S}_{ii}^{(c)} + \sum_{f=2i-1}^{2i}
  \zeta^{f}
  {\left(\Mm{\Delta}^{f}_{i}\right)}^{T}
  \mm{H}^{f} \mm{C}^{f,\min}
  \Mm{\Delta}^{f}_{i}
  \quad\text{(no summation over $i$)}.
\end{equation}
Here the matrix $\MM{S}_{ii}^{(c)}$ (no summation over $i$) is a positive
semidefinite and the matrix $\Mm{\Delta}_{i}^{f} = \Mm{B}_{i}^{f} -
\Mm{D}_{i}^{f}$ is the difference between the boundary derivative matrix from
$\MM{D}_{ii}$ (no summation over $i$) and the first derivative matrix
$\MM{D}_{i}$ at the boundary.
Each element of the diagonal matrix $\mm{C}^{f,min}$ is the minimum value of $c$
in the $m_{b}$ points orthogonal to the boundary where $m_{b}$ depends on the
order of accuracy of the SBP operator.
The positive constant $\zeta^f = h^f_{\bot} \bar{\zeta}$ where $h^{f}_{\bot}$ is
the grid spacing orthogonal to the face and $\bar{\zeta}$ is a constant which
depends on the SBP operator.
The $(m_{b}, \bar{\zeta})$ values used for the operators in this paper are given
in Table~\ref{tab:borrowing}; see \citet[Table 1]{AlmquistDunham2020}.
\begin{table}
  \centering
  \begin{tabular}{cccc}
    \toprule
    SBP interior order $2p$ & $\bar{\theta}$ & $\bar{\zeta}$ & $m_{b}$\\
    \midrule
    $2$ & $1 / 2$ & $1.0$ & $2$\\
    $4$ & $17 / 48$ & $0.5776$ & $4$\\
    $6$ & $13649 / 43200$ & $0.3697$ & $7$\\
    \bottomrule
  \end{tabular}
  \caption{Borrowing parameters and SBP norm $\mm{H}$ matrix corner value for
  used operators \citep[Table 1]{AlmquistDunham2020}.\label{tab:borrowing}}
\end{table}

Since $\MM{H}$ is diagonal and positive, it is clear that for any face $f$
\begin{equation}
  \VV{v}^{T} \MM{H} \VV{v}
  \ge
  \theta^{f}
  {\left(\Mm{L}^{f}\VV{v}_{i}\right)}^{T}
  \mm{H}^{f}
  \Mm{L}^{f}\VV{v}_{j}
  =
  \theta^{f}
  {\left(\VV{v}_{i}^{f}\right)}^{T}
  \mm{H}^{f}
  \VV{v}_{j}^{f},
\end{equation}
where $\theta^{f}$ is the value of the $\dof{\vv{H}}{00}$ where $\vv{H}$ is the
norm matrix orthogonal to the face.
It then follows that
\begin{equation}
  \begin{split}
    \VV{v}^{T} \MM{H} \VV{v}
    \ge
    \sum_{f = 1}^{2d}
    \frac{\theta^{f}}{d}
    {\left(\VV{v}^{f}\right)}^{T}
    \mm{H}^{f}
    \VV{v}^{f},
  \end{split}
\end{equation}
where the factor of $1/d$ is needed to avoid over counting corners and, when $d =
3$, edges.
Since the coefficient matrix $C_{ij}$ is positive definite, this can be extended
to included the variable coefficients:
\begin{equation}
  \label{eqn:vol:to:face}
  \begin{split}
    \VV{v}_{i}^{T} \MM{\hat{C}}_{ij} \MM{H} \VV{v}_{j}
    \ge
    \sum_{f = 1}^{2d}
    \frac{\theta^{f}}{d}
    {\left(\VV{v}_{i}^{f}\right)}^{T}
    \mm{\hat{C}}_{ij}^{f}
    \mm{H}^{f}
    \VV{v}_{j}^{f}.
  \end{split}
\end{equation}

We now turn to considering the discrete block
energy~\eqref{eqn:disc:block:energy}. The first term satisfies
\begin{equation}
  \label{eqn:app:rho}
  \frac{1}{2} \dot{\VV{u}}^{T} \MM{H}\MM{\rho}\dot{\VV{u}} \ge 0,
\end{equation}
because it is in quadratic form and $\MM{H}$ and $\MM{\rho}$ are diagonal,
positive matrices. The remaining terms will be shown to combine in a manner
that is also positive semidefinite.

Combing relations~\eqref{eqn:Mattsson:Def24},~\eqref{eqn:AlmquistDunham:Lem1},
and~\eqref{eqn:vol:to:face} we have that
\begin{equation}
  \label{eqn:app:A}
  \begin{split}
    \VV{u}^{T} \MM{A}_{ij}^{(\hat{C}_{ij})}\VV{u}
    \ge\;&
    \sum_{f=1}^{2d}\frac{\theta^{f}}{d}{\left(\Mm{D}_{i}^{f}\VV{u}\right)}^{T}
    \mm{\hat{C}}_{ij}^f \mm{H}^f
    \Mm{D}_{j}^{f} \VV{u}\\
    &
    +
    \sum_{k=1}^{d}
    \left(
    \sum_{f=2k-1}^{2k}
    \zeta^{f}
    {\left(\Mm{\Delta}^{f}_{k} \VV{u}\right)}^{T}
    \mm{H}^{f} \mm{\hat{C}}_{kk}^{f,\min}
    \Mm{\Delta}^{f}_{k}
    \VV{u}
    \right).
  \end{split}
\end{equation}
We now considering the face term of the discrete block
energy~\eqref{eqn:disc:block:energy}. Defining $\vv{\delta}^{f}_{u} =
\vv{u}^{*f} - \vv{u}^{f}$ and using the definition of $\vv{\hat{\tau}}^{f}$ and
$\vv{\hat{T}}^{f}$ in~\eqref{eqn:T:tau:hat:f} gives
\begin{equation}
  \begin{split}
    {\left(\vv{\hat{\tau}}^{f}\right)}^{T}
    {\left(\mm{X}^{f}\right)}^{-1}
    \mm{H}^{f}
    {\left(\vv{\hat{\tau}}^{f}\right)}
    -
    {\left(\vv{\hat{T}}\right)}^{T}
    {\left(\mm{X}^{f}\right)}^{-1}
    \mm{H}^{f}
    {\left(\vv{\hat{T}}\right)}
    & =
    2
    {\left(\vv{\hat{T}}\right)}^{T}
    \mm{H}^{f}
    \vv{\delta}^{f}_{u}
    +
    {\left(
    \vv{\delta}^{f}_{u}
    \right)}^{T}
    \mm{X}^{f}
    \mm{H}^{f}
    \vv{\delta}^{f}_{u}.
  \end{split}
\end{equation}
It is useful to note that $\vv{\hat{T}}$ can be rewritten using
$\Mm{\Delta}^{f}_{k}$ as
\begin{equation}
  \vv{\hat{T}}^{f} =
  \hat{n}^{f}_{i}\mm{\hat{C}}_{ij}^{f}\Mm{B}^{f}_{j}\VV{u}
  =
  \hat{n}^{f}_{i}\mm{\hat{C}}_{ij}^{f}\Mm{D}^{f}_{j}\VV{u}
  +
  \hat{n}^{f}_{k}\mm{\hat{C}}_{kk}^{f}\Mm{\Delta}^{f}_{k}\VV{u},
  \quad k = \ceil*{\frac{f}{2}};
\end{equation}
this follows because only when $f \in (2j, 2j-1)$ is $\Mm{B}^{f}_{j} \ne
\Mm{D}^{f}_{j}$. Using this along with the definition of $\mm{X}^{f}$
in~\eqref{eqn:Xf} leads to,
\begin{equation}
  \label{eqn:app:tau}
  \begin{split}
    &{\left(\vv{\hat{\tau}}^{f}\right)}^{T}
    {\left(\mm{X}^{f}\right)}^{-1}
    \mm{H}^{f}
    {\left(\vv{\hat{\tau}}^{f}\right)}
    -
    {\left(\vv{\hat{T}}\right)}^{T}
    {\left(\mm{X}^{f}\right)}^{-1}
    \mm{H}^{f}
    {\left(\vv{\hat{T}}\right)}\\
    &\quad=
    2\hat{n}_{i}^{f}
    {\left(\Mm{D}^{f}_{j}\VV{u}\right)}^{T}
    \mm{H}^{f}\mm{\hat{C}}_{ij}^{f}
    \vv{\delta}^{f}_{u}
    +
    \hat{n}^{f}_{i}
    \hat{n}^{f}_{j}
    {\left(\vv{\delta}^{f}_{u}\right)}^{T}
    \mm{\hat{C}}_{ij}^{f}
    \mm{\Gamma}^{f}
    \mm{H}^{f}
    \vv{\delta}^{f}_{u}
    +
    2
    \hat{n}_{k}^{f}
    {\left(\Mm{\Delta}^{f}_{k}\VV{u}\right)}^{T}
    \mm{H}^{f}
    \mm{\hat{C}}_{kk}^{f}
    \vv{\delta}^{f}_{u},
  \end{split}
\end{equation}
where $k = \ceil*{\frac{f}{2}}$.

Returning to the remaining terms of block energy~\eqref{eqn:disc:block:energy},
we use~\eqref{eqn:app:A} and~\eqref{eqn:app:tau} to write
\begin{equation}
  \label{eqn:app:A:tau}
  \begin{split}
    &\VV{u}^{T} \MM{A}_{ij}^{(C_{ij})}\VV{u}
    +
    \sum_{f=1}^{2d}
    \left(
    {\left(\vv{\hat{\tau}}^{f}\right)}^{T}
    {\left(\mm{X}^{f}\right)}^{-1}
    \mm{H}^{f}
    {\left(\vv{\hat{\tau}}^{f}\right)}
    -
    {\left(\vv{\hat{T}}\right)}^{T}
    {\left(\mm{X}^{f}\right)}^{-1}
    \mm{H}^{f}
    {\left(\vv{\hat{T}}\right)}
    \right)
    \\
    &\quad \ge
    \sum_{f=1}^{2d}
    \left(
      \frac{\theta^{f}}{d}{\left(\Mm{D}_{i}^{f} \VV{u}\right)}^{T} \mm{\hat{C}}_{ij}^f \mm{H}^f \Mm{D}_{j}^{f} \VV{u}
    +
    2 {\left(\Mm{D}^{f}_{j}\VV{u}\right)}^{T}
    \mm{\hat{C}}_{ij}^{f} \mm{H}^{f}
    \hat{n}_{i}^{f} \vv{\delta}_{u}^{f}
    \right)
    \\
    &\quad\quad +
    \sum_{k=1}^{d}
    \sum_{f=2k-1}^{2k}
    \left(
    \zeta^{f} {\left(\Mm{\Delta}^{f}_{k} \VV{u}\right)}^{T}
    \mm{\hat{C}}_{kk}^{f,\min} \mm{H}^{f}
    \Mm{\Delta}^{f}_{k} \VV{u}
    +
    2 {\left(\Mm{\Delta}^{f}_{k}\VV{u}\right)}^{T}
    \mm{\hat{C}}_{kk}^{f} \mm{H}^{f}
    \hat{n}_{k}^{f} \vv{\delta}^{f}_{u}
    \right)
    \\
    &\quad\quad +
    \sum_{f=1}^{2d}
    {\left( \vv{\delta}^{f}_{u}\right)}^{T}
    \hat{n}^{f}_{i}
    \mm{\hat{C}}_{ij}^{f} \mm{\Gamma}^{f} \mm{H}^{f}
    \hat{n}^{f}_{j} \vv{\delta}^{f}_{u}.
  \end{split}
\end{equation}
If we choose
\begin{subequations}
  \begin{align}
    \label{eqn:penalty:param}
    \mm{\Gamma}^{f} &\ge \frac{d}{\theta^{f}} \mm{I} + \frac{1}{\zeta^{f}}
    \mm{P}^{f},\\
    \mm{P}^{f} &= \mm{\hat{C}}_{kk}^{f} {\left(\mm{\hat{C}}_{kk}^{f,\min}\right)}^{-1},
  ~ k = \ceil*{\frac{f}{2}},
  \end{align}
\end{subequations}
then we have that
\begin{equation}
  \label{eqn:app:delta}
  \begin{split}
    &
    \sum_{f=1}^{2d}
    {\left( \vv{\delta}^{f}_{u}\right)}^{T}
    \hat{n}^{f}_{i}
    \mm{\hat{C}}_{ij}^{f} \mm{\Gamma}^{f} \mm{H}^{f}
    \hat{n}^{f}_{j} \vv{\delta}^{f}_{u}
    \\&\quad
    \ge
    \sum_{f=1}^{2d}
    \frac{d}{\theta^{f}}
    {\left( \vv{\delta}^{f}_{u}\right)}^{T}
    \hat{n}^{f}_{i}
    \mm{\hat{C}}_{ij}^{f} \mm{H}^{f}
    \hat{n}^{f}_{j} \vv{\delta}^{f}_{u}
    +
    \sum_{f=1}^{2d}
    \frac{1}{\zeta^{f}}
    {\left( \vv{\delta}^{f}_{u}\right)}^{T}
    \hat{n}^{f}_{i}
    \mm{\hat{C}}_{ij}^{f} \mm{P}^{f} \mm{H}^{f}
    \hat{n}^{f}_{j} \vv{\delta}^{f}_{u}
    \\&\quad
    =
    \sum_{f=1}^{2d}
    \frac{d}{\theta^{f}}
    {\left( \vv{\delta}^{f}_{u}\right)}^{T}
    \hat{n}^{f}_{i}
    \mm{\hat{C}}_{ij}^{f} \mm{H}^{f}
    \hat{n}^{f}_{j} \vv{\delta}^{f}_{u}
    +
    \sum_{k=1}^{d}
    \sum_{f = 2k-1}^{2k}
    \frac{1}{\zeta^{f}}
    {\left( \vv{\delta}^{f}_{u}\right)}^{T}
    \mm{\hat{C}}_{kk}^{f} \mm{P}^{f} \mm{H}^{f}
    \vv{\delta}^{f}_{u}
    \\&\quad
    =
    \sum_{f=1}^{2d}
    \frac{d}{\theta^{f}}
    {\left( \vv{\delta}^{f}_{u}\right)}^{T}
    \hat{n}^{f}_{i}
    \mm{\hat{C}}_{ij}^{f} \mm{H}^{f}
    \hat{n}^{f}_{j} \vv{\delta}^{f}_{u}
    +
    \sum_{k=1}^{d}
    \sum_{f = 2k-1}^{2k}
    \frac{1}{\zeta^{f}}
    {\left(\mm{P}^{f} \vv{\delta}^{f}_{u}\right)}^{T}
    \mm{\hat{C}}_{kk}^{f,\min}
    \mm{H}^{f}\mm{P}^{f}
    \vv{\delta}^{f}_{u},
  \end{split}
\end{equation}
where we have used that $\hat{n}^{f}_{i} \mm{\hat{C}}_{ij}^{f}
\hat{n}^{f}_{j} = \mm{\hat{C}}_{kk}^{f}$ with $k = \ceil*{\frac{f}{2}}$ (no
summation over $k$). Though a similar transformation could be used on the first
summation it is not needed and complicates the analysis that follows.
Returning to~\eqref{eqn:app:A:tau} then gives with~\eqref{eqn:app:delta}
\begin{equation}
  \begin{split}
    &\VV{u}^{T} \MM{A}_{ij}^{(C_{ij})}\VV{u}
    +
    \sum_{f=1}^{2d}
    \left(
    {\left(\vv{\hat{\tau}}^{f}\right)}^{T}
    {\left(\mm{X}^{f}\right)}^{-1}
    \mm{H}^{f}
    {\left(\vv{\hat{\tau}}^{f}\right)}
    -
    {\left(\vv{\hat{T}}\right)}^{T}
    {\left(\mm{X}^{f}\right)}^{-1}
    \mm{H}^{f}
    {\left(\vv{\hat{T}}\right)}
    \right)
    \\
    &\quad \ge
    \sum_{f=1}^{2d}
    \left(
      \frac{\theta^{f}}{d}{\left(\Mm{D}_{i}^{f} \VV{u}\right)}^{T} \mm{\hat{C}}_{ij}^f \mm{H}^f \Mm{D}_{j}^{f} \VV{u}
    +
    2 {\left(\Mm{D}^{f}_{j}\VV{u}\right)}^{T}
    \mm{\hat{C}}_{ij}^{f} \mm{H}^{f}
    \hat{n}_{i}^{f} \vv{\delta}_{u}^{f}
    \right)
    \\
    &\qquad +
    \sum_{f=1}^{2d}
    \frac{d}{\theta^{f}}
    {\left( \vv{\delta}^{f}_{u}\right)}^{T}
    \hat{n}^{f}_{i}
    \mm{\hat{C}}_{ij}^{f} \mm{H}^{f}
    \hat{n}^{f}_{j} \vv{\delta}^{f}_{u}
    \\
    &\quad\quad +
    \sum_{k=1}^{d}
    \sum_{f=2k-1}^{2k}
    \left(
    \zeta^{f} {\left(\Mm{\Delta}^{f}_{k} \VV{u}\right)}^{T}
    \mm{\hat{C}}_{kk}^{f,\min} \mm{H}^{f}
    \Mm{\Delta}^{f}_{k} \VV{u}
    +
    2 {\left(\Mm{\Delta}^{f}_{k}\VV{u}\right)}^{T}
    \mm{\hat{C}}_{kk}^{f} \mm{H}^{f}
    \hat{n}_{k}^{f} \vv{\delta}^{f}_{u}
    \right)
    \\
    &\quad\quad +
    \sum_{k=1}^{d}
    \sum_{f = 2k-1}^{2k}
    \frac{1}{\zeta^{f}}
    {\left(\mm{P}^{f} \vv{\delta}^{f}_{u}\right)}^{T}
    \mm{\hat{C}}_{kk}^{f,\min}
    \mm{H}^{f}\mm{P}^{f}
    \vv{\delta}^{f}_{u}
    \\
    &\quad =
    \sum_{f=1}^{2d}
    \frac{\theta^{f}}{d}
    {\left(\Mm{D}_{i}^{f} \VV{u} + \frac{d}{\theta^{f}} \hat{n}_{i}^{f}
    \vv{\delta}_{u}^{f} \right)}^{T}
    \mm{\hat{C}}_{ij}^f \mm{H}^f
    {\left(\Mm{D}_{j}^{f} \VV{u} + \frac{d}{\theta^{f}} \hat{n}_{j}^{f}
    \vv{\delta}_{u}^{f} \right)}
    \\
    &\quad\quad +
    \sum_{k=1}^{d}
    \sum_{f=2k-1}^{2k}
    \zeta^{f}
    {\left(\Mm{\Delta}^{f}_{k} \VV{u} +
    \frac{1}{\zeta^{f}}\hat{n}^{f}_{k}\mm{P}^{f}\vv{\delta}_{u}^{f}\right)}^{T}
    \mm{\hat{C}}_{kk}^{f,\min} \mm{H}^{f}
    {\left(\Mm{\Delta}^{f}_{k} \VV{u} +
    \frac{1}{\zeta^{f}}\hat{n}^{f}_{k}\mm{P}^{f}\vv{\delta}_{u}^{f}\right)},
  \end{split}
\end{equation}
where we have used that $\mm{\hat{C}}_{kk}^{f} = \hat{n}^{f}_{k} \mm{\hat{C}}_{kk}^{f}
\hat{n}_{k}^{t}$ (no summation over $k$) and $\mm{\hat{C}}_{kk}^{f,\min} \mm{P}^{f}
= \mm{C}_{kk}^{f}$ (no summation over $k$). Since this expression is in quadratic form, it is
non-negative and the when combine with~\eqref{eqn:app:rho} shows that the block
energy~\eqref{eqn:disc:block:energy} is non-negative.

\section{Non-Characteristic Boundary and Interface Treatment}
\label{app:standard:approach}

The standard approach for SBP-SAT for Dirichlet~\eqref{eqn:wave:Dirichlet},
and characteristic boundaries~\eqref{eqn:wave:char:bc:orig} as well as
computational and nonlinear interfaces from \citet{VirtaMattsson2014} and
\citet{Duru2019} are presented in the notation of this paper; Neumann boundary
treatment is the same as the characteristic boundary treat with $R = 1$.

\subsection{Dirichlet Boundary Conditions}
When block face $f$ is on a Dirichlet
boundary~\eqref{eqn:wave:transformed:Dirichlet} then the numerical fluxes are
chosen to be
\begin{subequations}
  \label{eqn:disc:bcs:Dirichlet}
  \begin{align}
    \vv{u}^{*f} &= \vv{g}_{D},\\
    \vv{\hat{\tau}}^{*f} &= \vv{\hat{\tau}}^{f};
  \end{align}
\end{subequations}
Using these numerical fluxes, the face energy rate of
change~\eqref{eqn:disc:faceenergy} is
\begin{equation}
  \begin{split}
    \dot{\mathcal{E}}^{f} =&\;
    {\left(\vv{\hat{\tau}}^{f}\right)}^{T}
    \mm{H}^{f}
    \dot{\vv{u}}^{f}
    +
    {\left(\vv{\hat{\tau}}^{f}\right)}^{T}
    \mm{H}^{f}
    {\left(
    \dot{\vv{g}}^{f}_{D}
    - \dot{\vv{u}}^{f}
    \right)}
    =
    {\left(\vv{\hat{\tau}}^{f}\right)}^{T}
    \mm{H}^{f}
    \dot{\vv{g}}^{f}_{D},
  \end{split}
\end{equation}
which with $g_{D} = 0$ gives $\dot{\mathcal{E}}^{f} = 0$ and does not lead to energy
growth.

\subsection{Characteristic (and Neumann) Boundary Conditions}
In order to define the standard treatment of characteristic boundary
conditions~\eqref{eqn:wave:char:bc:orig}, it is useful to
solve~\eqref{eqn:wave:char:bc:orig} for $\tau$:
\begin{equation}
  \tau = -\alpha \dot{u} + \nu g_C,
\end{equation}
with $\alpha = -Z (1 - R) / (R + 1) \le 0$ and $\nu = 1 / (R + 1)$.
We note again that the Neumann boundary condition is attained when $R = 1$ in
which case $\alpha = 0$ and $\nu = 1$.
With this, if block face $f$ is on a characteristic boundary then the numerical
fluxes are chosen to be
\begin{subequations}
  \label{eqn:disc:bcs:characteristic}
  \begin{align}
    \vv{u}^{*f} &= \vv{u}^{f},\\
    \vv{\hat{\tau}}^{*f} &= -\mm{\hat{\alpha}}^{f} \dot{\vv{u}}^{f} +
    \mm{\hat{\nu}}^{f} \vv{g}_{C},
  \end{align}
\end{subequations}
where the parameters $\mm{\hat{\alpha}}$ and $\mm{\hat{\nu}}$ are diagonal
matrices of $S_{J}^{f} \alpha$ and $S_{J}^{f} \nu$ evaluated at each point on
face $f$.
Using these numerical fluxes in~\eqref{eqn:disc:faceenergy} give
\begin{equation}
  \begin{split}
    \dot{\mathcal{E}}^{f} =&\;
    -{\left(\dot{\vv{u}}^{f}\right)}^{T}
    \mm{\hat{\alpha}}
    \mm{H}^{f}
    \dot{\vv{u}}^{f}
    +
    {\left(\vv{g}^{f}_{C}\right)}^{T}
    \mm{\hat{\nu}}
    \mm{H}^{f}
    \dot{\vv{u}}^{f}.
  \end{split}
\end{equation}
With $g_{C} = 0$ we then have that $\dot{\mathcal{E}}^{f} \le 0$ and there is no energy
growth due to the characteristic boundary treatment; equality is obtained in the
Neumann case.

\subsection{Computational Interface}
For computational interfaces (e.g., interfaces between blocks in the domains
that have been introduced to mesh to either a material interface and/or needed
in the mesh generation) continuity of displacement and traction need to be
enforced. That is, across the interface it is required that
\begin{equation}
  \label{eqn:computational:interface}
  \begin{split}
    u^{-} &= u^{+},\\
    n_{i}^{-}C_{ij}^{-}\partial_{j}u^{-} &=
    -n_{i}^{+}C_{ij}^{+}\partial_{j}u^{+}.
  \end{split}
\end{equation}
Here the superscript $\pm$ denotes the value on either side of the interface
with the unit normal $\vec{n}^{\pm}$ is taken to be outward to each side of the
interface, i.e., $\vec{n}^{-} = -\vec{n}^{+}$. The standard approach to
enforcing this is to choose the numerical flux to be the average of the values
on the two sides of the interface,
\begin{equation}
  \label{eqn:disc:computational:interface}
  \begin{split}
    \vv{u}^{*f^{-}} &= \frac{1}{2}\left(\vv{u}^{f^{-}} + \vv{u}^{f^{+}}\right),\\
    \vv{\hat{\tau}}^{*f^{-}} &= \frac{1}{2}\left(\vv{\hat{\tau}}^{f^{-}} -
    \vv{\hat{\tau}}^{f^{+}}\right);
  \end{split}
\end{equation}
the minus sign in $\vv{\hat{\tau}}^{*f^{-}}$ is due to the unit normals being equal and
opposite. Here the two blocks connected across the interface are $B^{\pm}$
through faces $f^{\pm}$.

The face energy rate of change~\eqref{eqn:disc:faceenergy} for computational
interfaces is then
\begin{equation}
  \begin{split}
    \dot{\mathcal{E}}^{f^{\pm}} =&\;
    \frac{1}{2}
    {\left(\vv{\hat{\tau}}^{f^{\pm}} - \vv{\hat{\tau}}^{f^{\mp}}\right)}^{T}
    \mm{H}^{f}
    \dot{\vv{u}}^{f^{\pm}}
    +
    \frac{1}{2}
    {\left(\vv{\hat{\tau}}^{f^{\pm}}\right)}^{T}
    \mm{H}^{f}
    {\left(
    \dot{\vv{u}}^{f^{\mp}}
    - \dot{\vv{u}}^{f^{\pm}}
    \right)}\\
    =&\;
    -
    \frac{1}{2}
    {\left(\vv{\hat{\tau}}^{f^{\mp}}\right)}^{T}
    \mm{H}^{f}
    \dot{\vv{u}}^{f^{\pm}}
    +
    \frac{1}{2}
    {\left(\vv{\hat{\tau}}^{f^{\pm}}\right)}^{T}
    \mm{H}^{f}
    \dot{\vv{u}}^{f^{\mp}}.
  \end{split}
\end{equation}
Adding the two sides of the interface together gives
\begin{equation}
  \dot{\mathcal{E}}^{f} = \dot{\mathcal{E}}^{f^{+}} + \dot{\mathcal{E}}^{f^{-}} = 0,
\end{equation}
and energy stability results.

\subsection{Nonlinear Interface Condition}
The approach \citet{Duru2019} for nonlinear interfaces is to define the sliding
velocity $V^{\pm f}$ directly from the particle velocities on the grid and then
the traction $\tau^{f}$ is defined directly from the nonlinear function so the
numerical fluxes are
\begin{equation}
  \begin{split}
    \vv{u}^{*f^{\pm}} &= \vv{u}^{f^{\pm}},\\
    \vv{\hat{\tau}}^{*f^{\pm}} &= \hat{F}\left(\vv{V}^{f^{\pm}}\right),
    ~
    \vv{V}^{f^{\pm}} = \left(\dot{\vv{u}}^{f^{\mp}} - \dot{\vv{u}}^{f^{\pm}}\right).
  \end{split}
\end{equation}
The face energy rate of change~\eqref{eqn:disc:faceenergy} for a nonlinear
interface is then
\begin{equation}
  \begin{split}
    \dot{\mathcal{E}}^{f^{\pm}} =&\;
    {\left(\hat{F}\left(\vv{V}^{f^{\pm}}\right)\right)}^{T}
    \mm{H}^{f}
    \dot{\vv{u}}^{f^{\pm}}.
  \end{split}
\end{equation}
Adding the two sides of the interface together gives
\begin{equation}
  \begin{split}
    \dot{\mathcal{E}}^{f} = \dot{\mathcal{E}}^{f^{+}} + \dot{\mathcal{E}}^{f^{-}}
    &=
    {\left(\hat{F}\left(\vv{V}^{f^{+}}\right)\right)}^{T}
    \mm{H}^{f}
    \dot{\vv{u}}^{f^{+}}
    +
    {\left(\hat{F}\left(\vv{V}^{f^{-}}\right)\right)}^{T}
    \mm{H}^{f}
    \dot{\vv{u}}^{f^{-}}\\
    &=
    {\left(\hat{F}\left(\vv{V}^{f^{+}}\right)\right)}^{T}
    \mm{H}^{f}
    \dot{\vv{u}}^{f^{+}}
    -
    {\left(\hat{F}\left(\vv{V}^{f^{+}}\right)\right)}^{T}
    \mm{H}^{f}
    \dot{\vv{u}}^{f^{-}}\\
    &=
    -
    {\left(\hat{F}\left(\vv{V}^{f^{+}}\right)\right)}^{T}
    \mm{H}^{f}
    \vv{V}^{f^{+}}\\
    &\le
    0,
  \end{split}
\end{equation}
where we have used that $\vv{V}^{f^{-}} = -\vv{V}^{f^{+}}$ and the fact that $V
\hat{F}(V) \ge 0$.

\section{Nonlinear Interface Root Finding Problem}\label{app:rootfinding}
In general, evaluating $\mathcal{Q}^{\pm}$ for a nonlinear condition
$\tau^{\pm} = F\left(V^{\pm}\right)$ requires solving a nonlinear root finding
problem.
In particular, using the characteristic variables $w^{\pm}$ a root finding
problem for $V^{\pm}$ is solved after which $\mathcal{Q}^{\pm}$ can be
determined.

Recall that force balance, $\tau^{-} = -\tau^{+}$, and the fact that
$V^{-} = -V^{+}$ implies that $\tau^{-} = -F\left(V^{+}\right)$. Using this we
can compute the $Z^{\pm}$ weighted-average
\begin{equation}
  \frac{Z^{-}\tau^{+} - Z^{+}\tau^{-}}{Z^{+} + Z^{-}}
  =
  F\left(V^{+}\right).
\end{equation}
Expressing $\tau^{\pm}$ in terms of $\mathcal{Q}^{\pm}$ and $w^{\pm}$,
see~\eqref{eqn:tau:Q}, then gives
\begin{equation}
  \label{eqn:fric:root:Q}
  \frac{
    Z^{-}\mathcal{Q}^{+} - Z^{-}w^{+}
    - Z^{+}\mathcal{Q}^{-} + Z^{+}w^{-}
  }{2(Z^{+} + Z^{-})}
  =
  F\left(V^{+}\right).
\end{equation}
The sliding velocity $V^{+}$ can be written in terms of the characteristic
variables using~\eqref{eqn:dotu:Q}:
\begin{equation}
  V^{+} = \dot{u}^{-} - \dot{u}^{+} =
  \frac{\mathcal{Q}^{-} + w^{-}}{2Z^{-}} -
  \frac{\mathcal{Q}^{+} + w^{+}}{2Z^{+}}
  =
  \frac{Z^{+} \mathcal{Q}^{-} + Z^{+} w^{-} - Z^{-}\mathcal{Q}^{+} - Z^{-} w^{+}}{2Z^{-}Z^{+}}.
\end{equation}
Using this, we can rewrite~\eqref{eqn:fric:root:Q} as
\begin{equation}
  \frac{Z^{+}Z^{-}
  }{(Z^{+} + Z^{-})}
  V^{+}
  +
  \frac{
      Z^{+}w^{-}
    - Z^{-}w^{+}
  }{(Z^{+} + Z^{-})}
  =
  F\left(V^{+}\right).
\end{equation}
This expression can be more compactly written by defining
\begin{align}
  \tau^{+}_{l} = \frac{ Z^{+}w^{-} - Z^{-}w^{+} }{(Z^{+} + Z^{-})},
\end{align}
which depends only on the characteristic variables propagating into the
interface and is the traction that would result if the interface were a
computational interface; seen by using~\eqref{eqn:computational:Q}
in~\eqref{eqn:tau:Q}.
We can now write the final form of the root finding problem as
\begin{equation}
  \label{eqn:fric:root}
  \eta
  V^{+}
  +
  \tau_{l}^{+}
  =
  F\left(V^{+}\right),
\end{equation}
where $\eta = Z^{+}Z^{-} / (Z^{+} + Z^{-})$ is known as the radiation damping
coefficient.
Once this nonlinear system is solved for $V^{+}$ all other quantities can be
determined using~\eqref{eqn:Q}.
When numerically solving~\eqref{eqn:fric:root} it is useful to realize that
$\sgn\left(V^{+}\right) = \sgn\left(\tau_{l}^{+}\right)$ and that the root can
be bracketed: $\left|V^{+}\right| \in \left[0,
F^{-1}\left(\tau_{l}^{+}\right)\right]$.

\bibliographystyle{spmpscinat}
\bibliography{refs}

\end{document}